\documentclass[11pt]{amsart}

\usepackage[usenames,dvipsnames]{xcolor}
\usepackage{tikz}
\usetikzlibrary{arrows.meta,decorations.markings,calc}
\usepackage{amssymb}
\usepackage{amsthm}
\usepackage{amsmath}
\usepackage{graphicx}
\usepackage{amsaddr}
\usepackage{comment}
\usepackage[margin = 1in]{geometry}
\usepackage{enumerate}
\usepackage{array}
\usepackage{float}
\usepackage[toc,page]{appendix}
\usepackage{lineno}
\usepackage{color}
\usepackage{bbm}
\usepackage{hyperref}
\usepackage[version=4]{mhchem}
\usepackage{siunitx}

\definecolor{ink}{HTML}{16243A}
\definecolor{acc}{HTML}{1F5FA8}
\definecolor{org}{HTML}{B04E17}
\definecolor{mut}{HTML}{5B677C}
\definecolor{cov}{HTML}{9DB3D1}
\definecolor{covo}{HTML}{E8A06A}

\definecolor{mygreen}{rgb}{0.1,0.75,0.2}

 \newtheorem{thm}{Theorem}[section]
 \newtheorem{cor}[thm]{Corollary}
 \newtheorem{lem}[thm]{Lemma}
 \newtheorem{prop}[thm]{Proposition}

 \newtheorem{defn}[thm]{Definition}

 \newtheorem{rem}[thm]{Remark}

 \newtheorem{exm}[thm]{Example}

 \numberwithin{equation}{section}

\newcommand{\pt}{\partial}
\newcommand{\epss}{\varepsilon}

\newcommand{\bR}{\mathbb{R}}

\newcommand{\sA}{\mathcal{A}}

\newcommand{\vg}{ \gamma}

\newcommand{\vb}{\vec{b}}

\newcommand{\vl}{\vec{\ell}}

\begin{document}

\title[Weak KAM]{Selection principle for Weak KAM solutions on $\bR^n$ via large deviation principle for diffusions}

\author[Y. Gao]{Yuan Gao}
\address{Department of Mathematics, Purdue University, West Lafayette, IN}
\email{gao662@purdue.edu}

\author[Y. Han]{Yuxi Han}
\address{Department of Mathematics, Purdue University, West Lafayette, IN}
\email{han891@purdue.edu}



\date{\today}

\begin{abstract}
    Stationary Hamilton--Jacobi equations generally admit many weak KAM solutions. The question is how to select one weak KAM solution that reflects a meaningful physical principle. In this paper, we propose a selection principle for Hamilton-Jacobi equations arising from the large deviation principle of general drift-diffusion processes on $\bR^n$ based on the exponential asymptotics of invariant measures. We  first characterize the Aubry set as the finitely many equilibria of the drift, using an action estimate from the orthogonal decomposition of the drift. For any boundary data on the Aubry set compatible with the Peierls barrier, we construct the unique coercive weak KAM solution with these boundary data, which is also the maximal locally Lipschitz viscosity solution. The Freidlin--Wentzell graph formula then provides the boundary data on the Aubry set that select the unique weak KAM solution. The selected solution is exactly the rate function in the large deviation principle for invariant measures in the zero noise limit. We   also show that this rate function is the vanishing viscosity limit, to which the viscous approximations converge  locally uniformly. Finally, under additional regularity assumptions, we   show that the long-time and zero-noise limits commute for suitably prepared initial densities.
\end{abstract}

\maketitle

\tableofcontents

\section{Introduction}

Classical KAM (Kolmogorov-Arnold-Moser) theory addresses the persistence of invariant quasi-periodic solutions in perturbed integrable Hamiltonian systems. It primarily relies on the existence of classical solutions for a perturbed integrable system, often described in terms of action-angle variables, and their connection to classical stationary solutions of the Hamilton-Jacobi equation (HJE) \cite{Kolmogorov1954, Moser1962}. In contrast, weak KAM theory, pioneered by   \textsc{Fathi} \cite{Fathi1997, Fathi2008} and \textsc{Ma\~n\'e} \cite{mane1996generic}, \textsc{E} \cite{weinan1999aubry}, investigates stationary solutions to the HJE in the broader "viscosity" sense; see viscosity solution introduced by \textsc{Crandall and Lions} \cite{crandall1983viscosity}. Weak KAM framework focuses on the long-time convergence of dynamical solutions to these stationary viscosity solutions and analyzes the global structure of the dynamics. Concurrently, to study the invariant behavior of dynamical systems, Aubry and Mather independently developed what is now known as Aubry-Mather theory, introducing fundamental concepts such as invariant sets (Mather sets, Aubry sets) and invariant measures \cite{Aubry1983, Mather1982}. These concepts are now understood to be deeply connected to the construction and variational representation of weak KAM solutions.

The variational representation of stationary solutions to the Hamilton--Jacobi equation uses the Mañé potential to characterize least-action paths and boundary data on invariant sets. A fundamental difficulty is that stationary solutions are generally nonunique, raising the question of which solution is dynamically relevant. Two standard selection mechanisms are vanishing viscosity and vanishing discount; as illustrated in Remark~\ref{rem}, these procedures may select different limiting solutions. In the stochastic setting considered here, vanishing viscosity is the natural choice because it corresponds to the zero-noise limit of the underlying dynamics. But even this selection mechanism is subtle: for mechanical Hamiltonians on the torus, vanishing-viscosity approximations need not converge uniformly \cite{LiuTranYu2026}, while in the joint vanishing-discount and vanishing-viscosity regime the limiting solution may depend on the relative scaling of the two parameters \cite{MitakeNiTran2026}. Related selection problems have also been investigated through physical solutions of the HJE \cite{AnantharamanEtAl2005} and through the Freidlin--Wentzell large deviation principle (LDP) for invariant measures \cite{SIMA23}; see also the variational framework of weak KAM theory \cite{Fathi2008}. The central question is therefore: when the viscous approximations do converge uniformly, which weak KAM solution do they select, and how can this selected solution be characterized in terms of the underlying dynamics?

In this paper, we answer this question firmly by proposing a selection principle for weak KAM solutions of a class of HJEs arising from general drift-diffusion processes. The associated Hamiltonian and Lagrangian come from the  optimal control and LDP of the drift diffusion processes. We explore  the Hamiltonian structure, invariant sets  and constructs a unique weak KAM solution to the  stationary HJE. The selected solution not only represents the LDP rate function of invariant measures but also serves as the uniform limit from the vanishing viscosity formulation.

  This selection principle via LDP was first proposed in \cite{SIMA23}
for a one-dimensional irreversible diffusion process on   torus. Their work is based on an explicit one-dimensional  potential function which allows them to characterize the Peierls barrier and the projected Mather and Aubry set for one-dimensional dynamics explicitly.   Our work extends the previous work from the one-dimensional torus   to the full space in multiple dimensions, and to general drift-diffusion processes on the whole space. Novel characterization of invariant sets and the well-posedness of weak KAM solutions in the whole space are established in this work.
On the other hand, the same Freidlin--Wentzell variational construction of the
rate function via the invariant measure was used in both the current work and the previous one--dimensional case. However, more delicate arguments for the convergence of the invariant measures and the uniform convergence of the vanishing viscosity limit in the whole space  are needed.

Precisely, we consider the overdamped Langevin diffusion in $\mathbb R^d$,
\begin{equation}\label{LV}
d X_t = \vb(X_t)\, dt + \sqrt{2\epss}\, d B_t.
\end{equation}
Here $B_t$ is standard $d$-dimensional Brownian motion and
$\varepsilon>0$. We assume (A1)--(A3) in
Section~\ref{sec:assumptions}. In particular, the drift admits the
orthogonal decomposition
\[
 \vb=-\nabla U+\vl,\qquad \vl\cdot\nabla U=0,
\]
where $U$ is coercive and has finitely many critical points, which
coincide with the equilibria of $\vb$. Along the deterministic dynamics
$\dot x=\vb(x)$,
\[
 \frac{d}{dt}U(x(t))=-|\nabla U(x(t))|^2.
\]
Thus forward trajectories remain in compact sublevel sets of $U$ and
converge to equilibria. This structure allows us to identify the Aubry
set and construct weak KAM solutions in the whole space.

The probability density $\rho_\varepsilon$ of \eqref{LV} satisfies the
Fokker--Planck equation~\eqref{FP},
\[
 \partial_t\rho_\varepsilon
 =-\nabla\cdot(\rho_\varepsilon\vb)+\varepsilon\Delta\rho_\varepsilon.
\]
For sufficiently small $\varepsilon>0$, the diffusion has a unique
invariant probability measure
$\pi_\varepsilon(dx)=m_\varepsilon(x)\,dx$. The logarithmic transform
$u_\varepsilon=-\varepsilon\log m_\varepsilon$ converts the stationary
Fokker--Planck equation into the viscous HJE~\eqref{S_HJEeps},
\[
 \nabla u_\varepsilon\cdot(\nabla u_\varepsilon+\vb)
 =\varepsilon(\Delta u_\varepsilon+\nabla\cdot\vb).
\]
Its formal zero-noise limit is the stationary HJE~\eqref{eq:HJE0},
\[
 H(\nabla u,x)=\nabla u\cdot(\nabla u+\vb)=0,
 \qquad x\in\mathbb R^d.
\]
This equation generally admits several weak KAM solutions. We first
show that each admissible choice of values on the Aubry set determines
a unique coercive weak KAM solution through the Peierls barrier.

To select these values, we use the Freidlin--Wentzell graph formula
for the invariant-measure LDP
\cite[Chapter~6]{Freidlin_Wentzell_2012}. With the normalization
$\min I=0$, the rate function $I$ describes the stationary probabilities
at the exponential scale: for every Borel set $E$ satisfying
$\inf_{E^\circ}I=\inf_{\overline E}I$,
\[
 \lim_{\varepsilon\downarrow0}\varepsilon\log\pi_\varepsilon(E)
 =-\inf_{x\in E}I(x).
\]
The graph formula prescribes the relative values of $I$ at the
equilibria. We show that its Peierls-barrier extension is the unique
coercive weak KAM solution with these prescribed Aubry values.
Thus the LDP selects one member of the family of weak KAM solutions;
the normalization alone does not give uniqueness.

Our results have two parts. Theorem~\ref{thm:main-weak-kam} assumes (A1) and (A3), and characterizes the Aubry set and proves existence and uniqueness for prescribed admissible Aubry values. Theorem~\ref{thm:main-selection}, under (A1)--(A3), identifies the LDP-selected solution
and establishes convergence to equilibrium at fixed noise. Under the
additional regularity assumption (A4), it also gives locally uniform
convergence $u_\varepsilon\to I$ as $\varepsilon\downarrow0$ and shows
that the long-time and zero-noise limits commute for initial densities
proportional to $e^{-I/\varepsilon}$.
The following two theorems summarize these results.

\begin{thm}[Informal summary of Proposition~\ref{prop:aubry-characterization}
and Corollary~\ref{cor:unique-coercive-weak-kam}]
\label{thm:main-weak-kam}
Assume (A1) and (A3), and let $V$ and $h$ denote the Ma\~n\'e potential
and Peierls barrier. Then the following statements hold.
\begin{enumerate}[(i)]
\item The Aubry set is finite and satisfies
$\mathcal A=\{x:\nabla U(x)=0\}=\{x:\vb(x)=0\}$.
\item For any values $\{W_i\}_{x_i\in\mathcal A}$ satisfying
$W_i-W_j\leq h(x_j,x_i)$ for all $x_i,x_j\in\mathcal A$, the function
\[
 W(x)=\min_{x_i\in\mathcal A}\{W_i+h(x_i,x)\}
     =\min_{x_i\in\mathcal A}\{W_i+V(x_i,x)\}
\]
is the unique coercive weak KAM solution of \eqref{eq:HJE0} with
$W(x_i)=W_i$. It is also the maximal locally Lipschitz viscosity solution
with these values.
\item For the function $W$ in (ii) and every $x\in\mathbb R^d$, there
is a backward calibrated curve ending at $x$ and converging to a point
of $\mathcal A$ as time tends to $-\infty$.
\end{enumerate}
\end{thm}

Our setting differs from standard weak KAM theory in two key respects. First, the domain is $\mathbb{R}^d$ rather than the compact torus $\mathbb{T}^d$. Second, while the Hamiltonian
\[
H(p,x) = |p|^2 + p \cdot \vb(x)
\]
is coercive in $p$ for each fixed $x$, it need not be uniformly coercive as $x$ varies over $\mathbb{R}^d$. Together, these features place the problem outside the reach of standard weak KAM theory, and no explicit characterization of the Peierls barrier $h$ is available for drifts satisfying only (A1) and (A3). We therefore establish the requisite compactness and structural properties directly.

To handle the noncompact domain, we exploit the orthogonal decomposition of $\vb$ and the coercivity of $U$. Specifically, Proposition~\ref{prop:minimizer-regularity}(a) gives
\[
\nabla U(x)\cdot q\leq L(q,x),
\]
and integrating along a path $\gamma$ yields
\[
U(\gamma(t))-U(\gamma(0))\leq\int_0^t L(\dot\gamma(s),\gamma(s))\,ds.
\]
Since $L\geq 0$ and $U$ is coercive, any family of paths starting in a fixed compact set with uniformly bounded action remains in a common compact set, independent of the time interval. This provides the compactness needed to obtain minimizers and analyze the Peierls barrier on $\mathbb{R}^d$.

The same estimates yield $h(x,x)=0$ at equilibria and $h(x,x)>0$ elsewhere, giving the Aubry-set characterization in Proposition~\ref{prop:aubry-characterization}. For admissible values on $\mathcal{A}$, the bound $W(x)\geq U(x)-C$ follows from the same inequality, so the Peierls-barrier extension $W$ is coercive. Compactness of minimizing paths then produces backward calibrated curves converging to equilibria (Proposition~\ref{prop:calibrated}). Comparing solutions along these curves extends their ordering from $\mathcal{A}$ to all of $\mathbb{R}^d$ and establishes uniqueness for the prescribed Aubry values (Proposition~\ref{AB}).

\begin{thm}[Informal summary of Proposition~\ref{prop:ldp}
and Theorem~\ref{thm:exchange}]
\label{thm:main-selection}
Assume (A1)--(A3). Define $W_0$ on $\mathcal A$ by the minimum rooted-graph
costs in \eqref{eqn:disbdy}, and extend it to $\mathbb R^d$ by the
formula in Theorem~\ref{thm:main-weak-kam}. Then the following statements hold.
\begin{enumerate}[(i)]
\item The function
\[
 I(x)=W_0(x)-\min_{x_j\in\mathcal A}W_0(x_j),
 \qquad \min_{\mathbb R^d}I=0,
\]
is the rate function supplied by the Freidlin--Wentzell
LDP for the invariant measures $(\pi_\varepsilon)$ of \eqref{LV}.
It is also the unique coercive weak KAM solution with the Aubry values
$I(x_i)=W_0(x_i)-\min_{x_j\in\mathcal A}W_0(x_j)$.

\item For every fixed $0<\varepsilon<1/(dN)$, where $N$ is as in (A2),
the diffusion has a unique invariant probability measure
$\pi_\varepsilon(dx)=m_\varepsilon(x)\,dx$, with $m_\varepsilon>0$.
For any initial probability density in $L^1(\mathbb R^d)$, the
probability solution $\rho_\varepsilon$ of \eqref{FP} satisfies
$\rho_\varepsilon(\cdot,t)\to m_\varepsilon$ in $L^1(\mathbb R^d)$
and pointwise as $t\to\infty$. Consequently,
$\psi_\varepsilon:=-\varepsilon\log\rho_\varepsilon$ converges
pointwise to $u_\varepsilon:=-\varepsilon\log m_\varepsilon$.

\item If (A4) holds, then $u_\varepsilon\to I$ locally uniformly in
$\mathbb R^d$ as $\varepsilon\downarrow0$.

\item Under (A4), choose the initial density
$\rho_\varepsilon(x,0)=Z_\varepsilon^{-1}e^{-I(x)/\varepsilon}$, where
$Z_\varepsilon=\int_{\mathbb R^d}e^{-I(x)/\varepsilon}\,dx$.
For every fixed $t\geq0$,
\[
 \psi_\varepsilon(\cdot,t)\longrightarrow I
 \quad\text{as }\varepsilon\downarrow0,
\]
locally uniformly in $x\in\mathbb R^d$.

\item For the initial family in (iv), the long-time and zero-noise
limits commute:
\[
 \lim_{\varepsilon\downarrow0}\lim_{t\to\infty}
 \psi_\varepsilon(x,t)
 =I(x)=
 \lim_{t\to\infty}\lim_{\varepsilon\downarrow0}
 \psi_\varepsilon(x,t),
 \qquad x\in\mathbb R^d.
\]
\end{enumerate}
\end{thm}

The exchange of limits in (v) uses the initial family
specified in (iv). For general initial probability densities, the two
iterated limits need not agree; see
Remark~\ref{rem:initial-profile-commutation}.

The proof combines the Freidlin--Wentzell formula with
Theorem~\ref{thm:main-weak-kam}. Local elliptic estimates upgrade the
LDP for $\pi_\varepsilon$ to locally uniform convergence
$u_\varepsilon\to I$; see
Corollary~\ref{cor:density-convergence-strong-regularity}.
Recurrence and transition kernel estimates give convergence at fixed
noise in Proposition~\ref{prop:fixed-epsilon-density-convergence}.
For the finite-time zero-noise limit, one could consider the upper and lower half-relaxed limits of $\psi_\varepsilon$ as $\varepsilon\downarrow0$. By stability, these limits are respectively a subsolution and a supersolution of the limiting equation \eqref{eq:limit-HJ}. A comparison principle would then imply that they coincide. However, the standard comparison argument for arbitrary sub- and supersolutions does not apply under our assumptions. Although $H(p,x)$ is coercive in $p$ for each fixed $x$, this coercivity need not be uniform in $x\in\mathbb R^d$, since no global bound on $\vb$ is assumed. Proposition~\ref{prop:coercive-comparison} gives comparison when the initial data are ordered and the supersolution is uniformly coercive at infinity. While the limiting initial profile is $I$, the required uniform coercivity of the supersolution half-relaxed limit is hard to establish here. We therefore prove the finite-time limit by comparing $\rho_\varepsilon$ with $m_\varepsilon$ on a large ball and using parabolic estimates to control the contribution from the initial tail; see Proposition~\ref{prop:finite-time-wkb}. The fixed-point identity $S_tI=I$ for the Lax--Oleinik operator then
yields the exchange of limits in Theorem~\ref{thm:exchange}.

The paper is organized as follows.
Section~\ref{sec:assumptions} states the assumptions and notation.
Sections~\ref{sec:preliminaries} and \ref{sec:aubry} develop the action
potentials and characterize the Aubry set.
Section~\ref{sec:weakKAM} constructs weak KAM solutions with prescribed
Aubry values and proves comparison and uniqueness.
Section~\ref{sec:ldp} relates the Freidlin--Wentzell graph formula to
these solutions and establishes logarithmic invariant-density
asymptotics. Section~\ref{sec:selection} proves convergence to
equilibrium, the finite-time zero-noise limit, and the exchange of
limits.

\bigskip

\section{Notation and assumptions}\label{sec:assumptions}

\subsection{Assumptions on $\vb$} Let $\vb:\mathbb R^d\to\mathbb R^d$ be a continuous vector field. We consider the following assumptions on $\vb$.

\noindent
(A1) $\vb$ admits the orthogonal decomposition
\begin{equation}\label{eq:b-decomposition}
    \vb(x)=-\nabla U(x)+\vl(x),
    \qquad
    \vl(x)\cdot \nabla U(x)=0,
    \qquad x\in\mathbb R^d,
\end{equation}
where $U\in C^1(\mathbb R^d)$ satisfies $U(x) \to \infty$ as $|x| \to \infty$, and $\vl \in C(\mathbb{R}^d; \mathbb{R}^d)$. Moreover, $\vb$ has finitely many zero points (with at least one such point) and satisfies
\[
    \nabla U(x) = 0  \quad \Longrightarrow \quad \vb(x)= 0 .
\]
Equivalently, the critical points of $U$ coincide with a subset of the zero set of $\vb$, that is,
\[
\left\{x:\nabla U(x) = 0\right\}=\left\{x:\vb(x)=0\right\}.
\]
In particular, $\vl$ vanishes at every critical point of $U$.

\noindent
(A2) There exists a sufficiently large $N$ such that
\[
\vb(x) \cdot x < -\frac{1}{N}, \qquad |x| >N.
\]

\noindent
(A3) $U \in C^2(\mathbb{R}^d)$ and $\vl \in C^1(\mathbb{R}^d; \mathbb{R}^d)$.

\noindent
(A4) $U \in C^3(\mathbb{R}^d)$ and $\vl \in C^2(\mathbb{R}^d; \mathbb{R}^d)$.

\vspace{3mm}

Note that the decomposition in \eqref{eq:b-decomposition} is distinct from the Helmholtz decomposition. The latter requires the vector field $\vl$ to satisfy $\nabla\cdot \vl=0$, while no such condition is assumed here. In addition, the condition $\vl\cdot\nabla U=0$ is not part of the Helmholtz decomposition.

Assumption (A2) is a small-noise Lyapunov condition. Indeed, for the generator
\[
\mathcal L_\varepsilon=\vb\cdot\nabla+\varepsilon\Delta
\]
and the Lyapunov function $\Phi(x)=|x|^2$, assumption (A2) gives
\[
\mathcal L_\varepsilon\Phi(x)
=2\vb(x)\cdot x+2d\varepsilon
<-\frac{2}{N}+2d\varepsilon,
\qquad |x|>N.
\]
Thus, whenever $0<\varepsilon<1/(dN)$, the generator has a strictly negative
drift outside $B(0,N)$. Together with the regularity needed for well-posedness
(for example, (A3)) and the nondegeneracy of the diffusion, the standard
Foster--Lyapunov and elliptic arguments yield a unique invariant probability
measure $\pi_\varepsilon$ with a strictly positive density $m_\varepsilon$;
see \cite[p.~110]{Freidlin_Wentzell_2012}. For more general sufficient
conditions, we refer to \cite[Corollary~2.4.2 and
Theorem~4.1.6]{Bogachev2015}.

\subsection{Notation}\label{sec:notation}

We collect the notation used throughout the paper. We write $|x|$ and
$x\cdot y$ for the Euclidean norm and inner product, and
$B_r(x)=\{y\in\mathbb R^d:|y-x|<r\}$, with $B_r=B_r(0)$.
The symbols $\mathrm{AC}$, $\mathrm{Lip}$, $\mathrm{USC}$, and
$\mathrm{LSC}$ denote absolutely continuous curves, Lipschitz functions,
upper semicontinuous functions, and lower semicontinuous functions,
respectively. In \eqref{LV}, $\varepsilon>0$ is the noise intensity,
$X_t$ is the Langevin diffusion, and $B_t$ is standard Brownian motion in
$\mathbb R^d$.

\begin{table}[H]
\centering
\small
\renewcommand{\arraystretch}{1.08}
\begin{tabular}{@{}p{0.20\textwidth}>{\raggedright\arraybackslash}p{0.75\textwidth}@{}}
\hline
\textbf{Symbol} & \textbf{Meaning} \\
\hline
$\vb$, $U$, $\vl$ & Drift, coercive potential, and component orthogonal to
$\nabla U$ in the decomposition \eqref{eq:b-decomposition}. \\
$\mathcal L_\varepsilon$ & Diffusion generator
$\vb\cdot\nabla+\varepsilon\Delta$. \\
$\rho_\varepsilon(x,t)$ & Time-dependent probability density solving the
Fokker--Planck equation \eqref{FP}. \\
$\pi_\varepsilon$, $m_\varepsilon$ & Invariant probability measure and its
density: $\pi_\varepsilon(dx)=m_\varepsilon(x)\,dx$. \\
$\psi_\varepsilon$, $u_\varepsilon$ & Logarithmic densities:
$\psi_\varepsilon=-\varepsilon\log\rho_\varepsilon$ and
$u_\varepsilon=-\varepsilon\log m_\varepsilon$. \\
$H(p,x)$ & Hamiltonian $p\cdot(p+\vb(x))$, where $p$ is the momentum. \\
$L(s,x)$ & Lagrangian $\tfrac14|s-\vb(x)|^2$, where $s$ is the velocity. \\
$c^*$ & Ma\~n\'e critical value; $c^*=0$ in the present setting. \\
$A_T[\gamma]$ & Action $\int_0^T L(\dot\gamma(t),\gamma(t))\,dt$ of an
absolutely continuous curve $\gamma$. \\
$h_T(x,y)$ & Minimum action from $x$ to $y$ in time $T$. \\
$V(x,y)$ & Ma\~n\'e potential $\inf_{T>0}h_T(x,y)$. \\
$h(x,y)$ & Peierls barrier $\liminf_{T\to\infty}h_T(x,y)$. \\
$\mathcal A$, $x_i$ & Aubry set $\{x:h(x,x)=0\}$ and its points. Here
$\mathcal A=\{x:\vb(x)=0\}$; see Proposition~\ref{prop:aubry-characterization}. \\
$\mathcal M_0$ & Projected Mather set; it coincides with $\mathcal A$ when
$\vb\in C^2$, as in Proposition~\ref{prop:calibrated}. \\
$\alpha(\gamma)$, $\omega(\eta)$ & Limit sets of curves as time tends to
$-\infty$ and $+\infty$, respectively. \\
$w\prec L$ & Domination of $w$ by the action of $L$; see
Definition~\ref{defn:weakKAM_neg}. \\
$W$ & Weak KAM extension of admissible values on $\mathcal A$, defined by
\mbox{$W(x)=\min_{y\in\mathcal A}\{W(y)+h(y,x)\}$} in \eqref{eqn:Wextension}. \\
$S_t$ & Backward Lax--Oleinik operator:
$(S_tf)(y)=\inf_{x\in\mathbb R^d}\{f(x)+h_t(x,y)\}$. \\
$G(E)$ & Acyclic directed graphs on a finite set $S$, with one outgoing
edge from each vertex in $S\setminus E$ and none from $E$.
In \eqref{eqn:disbdy}, $S=\mathcal A$ and $G\{x_i\}=G(\{x_i\})$. \\
$W_0$ & Boundary values selected by the Freidlin--Wentzell graph costs
in \eqref{eqn:disbdy}, together with their weak KAM extension to
$\mathbb R^d$. \\
$I$ & Selected rate function
$I(x)=W_0(x)-\min_{x_i\in\mathcal A}W_0(x_i)$, normalized by \mbox{$\min I=0$}. \\
$Z_\varepsilon$ & Normalizing constant
$\int_{\mathbb R^d}e^{-I(x)/\varepsilon}\,dx$ for the initial density
\mbox{$\rho_\varepsilon(x,0)=Z_\varepsilon^{-1}e^{-I(x)/\varepsilon}$}. \\
\hline
\end{tabular}
\end{table}

\bigskip

\section{Preliminaries}\label{sec:preliminaries}
In this section, we introduce several fundamental concepts that are essential for the subsequent construction of the weak KAM solution. While some of these results are standard on compact domains such as $\mathbb{T}^d$, our analysis is carried out on the noncompact space $\mathbb{R}^d$. Hence, we need to carefully verify their validity in this setting.

\subsection{Associated Hamiltonian and Lagrangian}
In the equations below, we suppress the spatial and time
arguments of functions.
Unless indicated otherwise, $\nabla$ and $\Delta$ act on the spatial variable $x$.

The density $\rho_\varepsilon(x,t)$ of $X_t$ satisfying \eqref{LV} solves the
Fokker--Planck equation
\begin{equation}\label{FP}
\pt_t \rho_\varepsilon = -\nabla \cdot \left(\rho_\varepsilon \vb\right) + \epss\Delta \rho_\varepsilon,
\qquad (x,t)\in\mathbb R^d\times(0,\infty).
\end{equation}
Applying the WKB transformation $\rho_\varepsilon = e^{-\psi_\varepsilon/\varepsilon}$, and substituting this expression into \eqref{FP}, yields the equation satisfied by $\psi_\varepsilon$:
\begin{equation}\label{eqn:eplog}
    \partial_t \psi_\varepsilon  +\nabla \psi_\varepsilon \cdot \left(\nabla \psi_\varepsilon + \vb\right)= \varepsilon \left(\nabla \cdot\vb\right) + \varepsilon\Delta\psi_\varepsilon,
    \qquad (x,t)\in\mathbb R^d\times(0,\infty).
\end{equation}
 Taking $\epss\to 0$ formally, we obtain the Hamilton--Jacobi equation (HJE)
 \begin{equation}\label{eq:limit-HJ}
 \pt_t \psi + H(\nabla \psi, x) = 0,
 \qquad (x,t)\in\mathbb R^d\times(0,\infty),
 \end{equation}
 with a quadratic Hamiltonian $H: \mathbb{R}^d \times \mathbb{R}^d \to \mathbb{R}$ given by
 \begin{equation} \label{eqn:hamiltonian}
 H\left(p,x\right) = p \cdot \left(p+\vb(x)\right).
 \end{equation}
 Then the corresponding Lagrangian $L: \mathbb{R}^d \times \mathbb{R}^d \to \mathbb{R}$, which is the Legendre transform of $H$, is given by
 \begin{equation}\label{eqn:Lagrangian}
 L(s,x) = \sup_{p \in \mathbb{R}^d} \left\{s \cdot p - H\left(p,x\right)\right\}
 = \frac{1}{4}\left|s-\vb(x)\right|^2,
 \end{equation}
 where the maximizer $p^*$ solves
 \begin{equation}
 s = \nabla_p H\left(p^*, x\right) = 2p^* +\vb(x).
 \end{equation}
Note that $H$ and $L$ are strictly convex in $p$ and $s$, respectively.

The density $m_\varepsilon(x)$ of the invariant measure $\pi_\varepsilon$
satisfies the stationary Fokker--Planck equation
\begin{equation}
\nabla \cdot (\vb m_\epss) = \epss \Delta m_\epss,
\qquad x\in\mathbb R^d.
\end{equation}
Setting $m_\epss = e^{-u_\epss/\epss}$ and substituting it into the stationary Fokker--Planck equation yields
\begin{equation}\label{S_HJEeps}
 \nabla u_\epss \cdot \left(\nabla u_\epss +\vb\right) = \epss \left(\Delta u_\epss + \nabla \cdot \vb\right),
 \qquad x\in\mathbb R^d.
\end{equation}
 Formally letting $\varepsilon \to 0$, we obtain the stationary Hamilton--Jacobi equation
\begin{equation}\label{eq:HJE0}
 \nabla u \cdot \left(\nabla u +\vb\right)=0,
 \qquad x\in\mathbb R^d.
\end{equation}

One essential property of the Hamiltonian induced by an underlying stochastic process is
\begin{equation}\label{H_f1}
H(0, x) \equiv 0, \quad \text{ and thus }\quad  \nabla_x H(0, x) \equiv 0,
\qquad x\in\mathbb R^d.
\end{equation}
From this, immediate consequences for the corresponding Lagrangian are
\begin{equation}
L(s,x) = \sup_{p\in \mathbb{R}^d} \left(s \cdot p - H\left(p,x\right)\right) \geq -H(0,x) =0,
\end{equation}
and $L(s,x) = 0$ if and only if $s = \nabla_p H(0, x)$.

\subsection{Ma\~n\'e critical value, Ma\~n\'e potential, and Peierls Barrier}

Throughout this subsection, we assume (A1), with additional assumptions stated explicitly when needed.

We now introduce the Ma\~n\'e critical value associated with the Hamiltonian $H$ defined in \eqref{eqn:hamiltonian}. Since the spatial domain is $\mathbb{R}^d$ rather than a compact domain such as $\mathbb{T}^d$, some standard results from weak KAM theory require additional justification.

\begin{defn} Assume (A1). Define the Ma\~n\'e critical value as
    \begin{equation}
        c^* = \inf\left\{ c \in \mathbb{R} \;\middle|\; \exists\, u \in \mathrm{Lip}(\bR^d) \text{ such that } H(Du(x),x) \leq c \text{ a.e. in } \bR^d \right\},
    \end{equation}
    where $H$ is given by \eqref{eqn:hamiltonian}.
\end{defn}

\begin{prop}
    $c^* = 0$.
\end{prop}

\begin{proof}
For $u \equiv0$, we have $H(Du,x)=H(0,x) \equiv 0$ in $\mathbb{R}^d$. Hence, $c^\ast \leq 0$.

Now suppose, by contradiction, that $c^\ast<0$.
By the definition of $c^\ast$, there exist a constant $c<0$ and a function
$u\in \operatorname{Lip}(\mathbb{R}^d)$ such that
\[
    0>c\geq H(Du,x)
    =|Du|^2+Du\cdot\vb(x)
    \geq -\frac{|\vb(x)|^2}{4},
    \qquad \text{for a.e. }x\in\mathbb{R}^d .
\]
Let $x_0$ be a zero point of $\vb$, namely, $\vb(x_0)=0$.
Since the set of differentiability
points of $u$ has full measure, there exists a sequence of differentiability points
$x_n\to x_0$ such that the above inequality holds at each $x_n$. Passing to the limit and
using the continuity of $\vb$, we obtain
\[
    0>c
    \geq \lim_{n\to\infty}\left(-\frac{|\vb(x_n)|^2}{4}\right)
    =-\frac{|\vb(x_0)|^2}{4}=0,
\]
which is a contradiction. Therefore, $c^\ast =0$.
\end{proof}

We now recall the definitions of the Ma\~n\'e potential and the Peierls barrier in the case \(c^*=0\).

\begin{defn}
Let \(T>0\). For any \(\gamma\in \mathrm{AC}([0,T];\mathbb{R}^d)\), define the action functional
\begin{equation}
    A_T[\gamma]
    := \int_0^T L(\dot{\gamma}(t),\gamma(t))\,dt
    = \frac{1}{4}\int_0^T \left|\dot{\gamma}(t)-\vb(\gamma(t))\right|^2\,dt .
\end{equation}
Define $h_T: \mathbb{R}^d \times \mathbb{R}^d \to \mathbb{R}$ by
\[
 h_T(x, y):=\inf _{\substack{\gamma \in \operatorname{AC}\left([0, T]; \,\bR^d\right) \\ \gamma(0)=x,\,\gamma(T)=y}}   A_T[\gamma].
\]
The Ma\~n\'e potential associated with the Lagrangian \(L\) in \eqref{eqn:Lagrangian} is defined by
\begin{equation}\label{M_bm}
V(x,y)
:= \inf_{T>0} h_T(x,y).
\end{equation}
The Peierls barrier $h: \bR^d \times \bR^d \rightarrow \mathbb{R}$ is defined by
\[
h\left(x, y\right):=\liminf_{T \rightarrow+\infty} h_T(x, y).
\]
\end{defn}

By definition, $V(x,y)\leq h(x,y)$ for all $x,y\in\mathbb{R}^d$.
Equality holds when either endpoint belongs to the Aubry set; see
Proposition~\ref{prop:peierls-mane-aubry} and, for the equality on the
Aubry set in the compact setting, \cite[Remark~6.28]{sorrentino2010mather}.
Since $L\geq0$, an absolutely continuous curve has zero action if and only
if it follows the deterministic drift:
\begin{equation}
    \dot{\gamma}(t)=\vb(\gamma(t)),\qquad \gamma(0)=x,
    \qquad t\in[0,T].
\end{equation}

We now establish properties of $h$, $h_T$, and minimizers of $h_T$.
Standard weak KAM theory on compact manifolds uses $C^2$ Tonelli
Lagrangians, with a positive definite velocity Hessian and superlinear
growth in velocity; see \cite[Definition~2.1]{sorrentino2010mather}.
Our $L$ is quadratic and strictly convex in $s$, and superlinear in $s$
locally uniformly in $x$. However, (A1) guarantees only joint
continuity, while (A3) guarantees $C^1$ regularity; $C^2$ regularity is not
assumed. Moreover, the spatial domain is noncompact, global uniform
superlinearity in velocity is not guaranteed, and $L$ is not even coercive.
We use the coercive potential $U$ and the action estimate in
Proposition~\ref{prop:minimizer-regularity}(a) to obtain compactness of
minimizing sequences.

\begin{prop}\label{prop:minimizer-regularity}
Assume (A1). Let $T>0$, $x,y \in \mathbb{R}^d$. Then the following statements hold.

\noindent
(a) For any $\gamma \in \mathrm{AC}\left([0,T];\mathbb{R}^d\right)$, we have
\[
\frac{d}{dt}U(\gamma(t)) \leq \frac{1}{4} \left|\dot{\gamma}(t)-\vb(\gamma(t))\right|^2, \qquad \text{for a.e. }t\in [0,T].
\]

\noindent
(b) There exists a constant $C=C\left(x,y,T\right)>0$ such that $U(y)-U(x)\leq h_T(x,y) \leq C$.
Moreover, $U(y)-U(x)\leq h(x,y)<\infty$ and $U(y)-U(x)\leq V(x,y)<\infty$.

\noindent
(c) There exists \(\gamma \in \mathrm{AC}([0,T];\mathbb{R}^d)\) such that
\[
h_T(x,y)=\frac{1}{4}\int_0^T \left|\dot{\gamma}(t)-\vb(\gamma(t))\right|^2\,dt,
\qquad \gamma(0)=x,\quad \gamma(T)=y.
\]
Moreover, there exists a compact set \(K=K(x,y,T)\subset \mathbb{R}^d\) such that every minimizer \(\gamma\) of \(h_T(x,y)\) satisfies
$\gamma([0,T])\subset K$. In addition, if we assume (A3), then every minimizer $\gamma\in C^2([0,T];\mathbb{R}^d)$.

\noindent
(d) $h_T$ is continuous on $\mathbb{R}^d \times \mathbb{R}^d$ and $h$ is locally Lipschitz on $\mathbb{R}^d \times \mathbb{R}^d$.

\end{prop}

\begin{proof}

To prove part (a), we compute for any $x,v\in\mathbb R^d$,
\begin{equation}\label{eq:U-calibration-identity}
\begin{aligned}
 L(v,x)-\nabla U(x)\cdot v
 &=\frac14|v+\nabla U(x)-\vl(x)|^2-\nabla U(x)\cdot v\\
 &=\frac14|v-\nabla U(x)-\vl(x)|^2,
\end{aligned}
\end{equation}
where we used (A1) in the last equality. Equality
$L(v,x)=\nabla U(x)\cdot v$ holds if and only if
$v=\nabla U(x)+\vl(x)$. Hence, for any
$\gamma\in\mathrm{AC}([0,T];\mathbb R^d)$ and almost every
$t\in[0,T]$, we have
\[
\begin{aligned}
 \frac{d}{dt}U(\gamma(t))
 &=\nabla U(\gamma(t))\cdot\dot\gamma(t)\\
 &\leq L(\dot\gamma(t),\gamma(t))
 =\frac14|\dot\gamma(t)-\vb(\gamma(t))|^2.
\end{aligned}
\]
Now we prove part (b). For any
$\gamma\in\mathrm{AC}([0,T];\mathbb R^d)$ with $\gamma(0)=x$ and
$\gamma(T)=y$, integrating \eqref{eq:U-calibration-identity} gives
\[
\begin{aligned}
 A_T[\gamma]
 &=U(y)-U(x)
 +\frac14\int_0^T
 |\dot\gamma(t)-\nabla U(\gamma(t))-\vl(\gamma(t))|^2\,dt\\
 &\geq U(y)-U(x),
\end{aligned}
\]
which implies that
\[
 U(y)-U(x)\leq h_T(x,y),\qquad
 U(y)-U(x)\leq h(x,y),\qquad
 U(y)-U(x)\leq V(x,y).
\]

On the other hand, consider the straight line connecting $x$ and $y$ using time $T$
\[
\eta(t)=x+\frac{t}{T}(y-x), \qquad t\in[0,T].
\]
Then
\[
\begin{aligned}
h_T(x,y) &\le \int_0^T \frac14\left|\frac{y-x}{T}-\vb(\eta(s))\right|^2 ds\le \frac14\int_0^T \left(2\frac{|y-x|^2}{T^2}+2\left|\vb(\eta(s))\right|^2\right)ds \leq C,
\end{aligned}
\]
for some constant $C=C \left(|x|,|y|,T,\vb\right)>0$.

Next, we prove part (c). For any $n\in \mathbb{N}$, there exists $\gamma_n \in \mathrm{AC}([0,T]; \mathbb{R}^d)$ such that
\[
A_T[\gamma_n]=\frac{1}{4}\int_0^T \left|\dot{\gamma}_n(t)-\vb(\gamma_n(t))\right|^2\,dt\leq h_T(x,y)+\frac{1}{n}, \quad \gamma_n(0)=x,\ \gamma_n(T)=y.
\]
By part (a) and (b), we know
\[
U(\gamma_n(t))-U(x)\le \frac{1}{4}\int_0^t \left|\dot{\gamma}_n(s)-\vb(\gamma_n(s))\right|^2\,ds\leq h_T(x,y)+1\leq C,
\]
for some constant $C>0$. Then $\gamma_n([0,T]) \subset \left\{z\in \mathbb{R}^d: U(z) \leq U(x)+C\right\}:=K$. Since $U$ is coercive, $K$ is compact. Hence,
\[
\begin{aligned}
\int_0^T \left|\dot{\gamma}_n(t)\right|^2\,dt
&\le 2\int_0^T \left|\dot{\gamma}_n(t)-\vb(\gamma_n(t))\right|^2\,dt+2\int_0^T \left|\vb(\gamma_n(t))\right|^2\,dt\\
&\le 8 A_T[\gamma_n]+C \le 8(h_T(x,y)+1) +C \leq C.
\end{aligned}
\]
Thus $(\gamma_n)$ is bounded in $W^{1,2}([0,T];\mathbb R^d)$. After passing to a subsequence,
\[
\gamma_n \rightharpoonup  \gamma \text{ weakly in } W^{1,2}([0,T]; \mathbb{R}^d), \qquad \gamma_n \to \gamma \text{ uniformly on } [0,T].
\]
Uniform convergence preserves the endpoints and, by
continuity of $\vb$ on $K$, gives $\vb(\gamma_n)\to\vb(\gamma)$ uniformly.
Hence $\dot\gamma_n-\vb(\gamma_n)\rightharpoonup
\dot\gamma-\vb(\gamma)$ in $L^2(0,T;\mathbb R^d)$.
By weak lower semicontinuity of the $L^2$ norm,
\[
\frac{1}{4}\int_0^T \left|\dot{\gamma}-\vb(\gamma)\right|^2\,dt \leq \liminf_{n \to \infty} \frac{1}{4} \int_0^T \left|\dot{\gamma}_n-\vb(\gamma_n)\right|^2\,dt=h_T(x,y).
\]
Therefore, $\gamma \in \mathrm{AC}([0,T];\mathbb{R}^d)$ is a minimizer of $h_T(x,y)$.

Consider any minimizer $\gamma$ of $h_T(x,y)$. By a similar argument as above, we have
\[
\gamma([0,T]) \subset K, \qquad \gamma \in W^{1,2}([0,T];\mathbb{R}^d).
\]
For regularity, assume (A3).
Define $p(t):=\frac12\left(\dot\gamma(t)-\vb(\gamma(t))\right)$.
Then $p\in L^2([0,T];\mathbb R^d)$. Let
$\eta\in C^1([0,T];\mathbb R^d)$ with $\eta(0)=\eta(T)=0$. Since $\gamma$
is a minimizer of $h_T(x,y)$, we have
\[
\left.\frac{d}{d\varepsilon}\right|_{\varepsilon=0}
A_T[\gamma+\varepsilon\eta]=0 \quad \Longrightarrow \quad 0
=
\int_0^T
p(t)\cdot
\left(\dot\eta(t)-D\vb(\gamma(t))\eta(t)\right)\,dt,
\]
which implies
\[
\int_0^T p(t)\cdot \dot\eta(t)\,dt
=
\int_0^T
\bigl(D\vb(\gamma(t))^T p(t)\bigr)\cdot \eta(t)\,dt .
\]
Thus, in the sense of distributions,
\begin{equation}\label{eqn:pdot}
\dot p(t)
=
-D\vb(\gamma(t))^T p(t).
\end{equation}
Since $\vb\in C^1$ and $\gamma([0,T])$ is compact, $D\vb(\gamma(t))$ is bounded on
$[0,T]$. Since $p\in L^2([0,T];\mathbb R^d)$, we obtain $\dot p\in L^2([0,T];\mathbb R^d)$.
Therefore, $p\in W^{1,2}([0,T];\mathbb R^d)$. By the Sobolev embedding, $W^{1,2}([0,T]) \subset
C^{0,1/2}([0,T])$. Hence $p$ has a continuous representative, and \eqref{eqn:pdot} implies
$p\in C^1([0,T];\mathbb R^d)$. Now
\begin{equation} \label{eqn:gammadot}
\dot\gamma(t)=\vb(\gamma(t))+2p(t), \qquad \text{for a.e.} \, t\in[0,T],
\end{equation}
and hence
$\gamma\in C^1([0,T];\mathbb R^d)$. Using again \eqref{eqn:gammadot}, we see that $\gamma\in C^2([0,T];\mathbb R^d)$, which completes the proof of part (c).

Now we prove part (d), using only (A1). Let \(R>0\), and let
\(x,x',y,y'\in\overline{B(0,R)}\). The straight line connecting
$x$ and $y$ gives the uniform bound
\[
h_T(x,y)\leq C_0
:=\frac{T}{4}
\left(\frac{2R}{T}+\sup_{z\in\overline{B(0,R)}}|\vb(z)|\right)^2.
\]
If $\gamma$ is any minimizer for $h_T(x,y)$, then part (a) gives
\[
U(\gamma(t))
\leq U(x)+\frac14\int_0^t
|\dot\gamma(s)-\vb(\gamma(s))|^2\,ds
\leq \max_{\overline{B(0,R)}}U+C_0.
\]
Consequently, every such minimizer is contained in the compact set
\[
K:=\left\{z\in\mathbb R^d:
U(z)\leq\max_{\overline{B(0,R)}}U+C_0\right\},
\]
which is independent of the particular endpoints in $\overline{B(0,R)}$.
Moreover,
\begin{equation}\label{eq:uniform-L2-action}
\|\dot\gamma-\vb(\gamma)\|_{L^2(0,T)}^2
=4h_T(x,y)\leq4C_0.
\end{equation}

\noindent
(i) Define $\tilde{\gamma}:[0, T]\to \mathbb{R}^d$ by
\[
   \tilde{\gamma}(t) = \left\{\begin{aligned}
   & \gamma(t)+\frac{T-4t}{T}\left(x'-x\right), \qquad \qquad\,\,\,  \text{ if } t\in \left[0, \frac{T}{4}\right], \\
   &\gamma(t), \qquad \qquad \qquad \qquad \qquad  \qquad \text{ if } t \in \left[\frac{T}{4}, \frac{3T}{4}\right],\\
   &\gamma(t)+\frac{4t-3T}{T}\left(y'-y\right), \qquad  \text{ if } t\in \left[\frac{3}{4}T, T\right].
   \end{aligned}
   \right.
\]
Set
\[
\delta:=|x'-x|+|y'-y|.
\]
Then $|\tilde\gamma(t)-\gamma(t)|\leq \delta\leq4R$ and
\begin{equation}\label{eq:endpoint-perturbation-L2}
\|\dot{\tilde\gamma}-\dot\gamma\|_{L^2(0,T)}^2
=\frac4T\left(|x'-x|^2+|y'-y|^2\right)
\leq\frac{4\delta^2}{T}.
\end{equation}
Both $\gamma$ and $\tilde\gamma$ take values in the compact set
$K':=K+\overline{B(0,4R)}$. Let
\[
\omega(r):=\sup_{\substack{z,z'\in K'\\|z-z'|\leq r}}
|\vb(z)-\vb(z')|
\]
be the modulus of continuity of $\vb$ on $K'$. It follows from
\eqref{eq:endpoint-perturbation-L2} that
\begin{equation}\label{eq:endpoint-error-L2}
\begin{aligned}
\left\|
\bigl(\dot{\tilde\gamma}-\vb(\tilde\gamma)\bigr)
-\bigl(\dot\gamma-\vb(\gamma)\bigr)
\right\|_{L^2(0,T)} &\leq
\|\dot{\tilde\gamma}-\dot\gamma\|_{L^2(0,T)}
+\|\vb(\tilde\gamma)-\vb(\gamma)\|_{L^2(0,T)}\\
&\leq \frac{2\delta}{\sqrt T}+\sqrt T\,\omega(\delta)=:E(\delta),
\end{aligned}
\end{equation}
where $E(\delta)\to0$ as $\delta\to0$. Using the difference-of-squares identity,
Cauchy--Schwarz, \eqref{eq:endpoint-error-L2}, and
\eqref{eq:uniform-L2-action}, we obtain
\[
\begin{aligned}
&h_T(x',y')-h_T(x,y)\\
&\leq\frac14\int_0^T
\left(
|\dot{\tilde\gamma}-\vb(\tilde\gamma)|^2
-|\dot\gamma-\vb(\gamma)|^2
\right)\,dt\\
&\leq\frac14
\left\|
\bigl(\dot{\tilde\gamma}-\vb(\tilde\gamma)\bigr)
-\bigl(\dot\gamma-\vb(\gamma)\bigr)
\right\|_{L^2(0,T)}\left(
\left\|
\bigl(\dot{\tilde\gamma}-\vb(\tilde\gamma)\bigr)
-\bigl(\dot\gamma-\vb(\gamma)\bigr)
\right\|_{L^2(0,T)}
+2\|\dot\gamma-\vb(\gamma)\|_{L^2(0,T)}
\right)\\
&\leq \sqrt{C_0}\,E(\delta)+\frac14E(\delta)^2.
\end{aligned}
\]
The constants and compact sets above are uniform for all endpoints in
$\overline{B(0,R)}$. Exchanging $(x,y)$ with $(x',y')$ consequently gives
\[
|h_T(x,y)-h_T(x',y')|
\leq \sqrt{C_0}\,E(\delta)+\frac14E(\delta)^2.
\]
Thus $h_T$ is continuous.

\noindent
(ii) Let $a:= |x'-x|$ and $c:=|y'-y|$.
For each $T>0$, choose a minimizer $\gamma$ of $h_T(x,y)$.
Without loss of generality, assume $a,c>0$.
Consider the straight line connecting $x'$ to $x$ given by
\[
\eta(s): = x' +\frac{s(x-x')}{a}, \qquad s \in [0,a].
\]
Since $\eta \left([0, a]\right) \subset \overline{B(0,R)}$ and $\vb$ is continuous,
\begin{equation}\label{eqn:eta}
A_{a}[\eta]= \frac{1}{4} \int_0^{a} \left|\frac{x-x'}{|x-x'|}-\vb\left(\eta(s)\right)\right|^2 \, ds \leq Ma,
\end{equation}
where $M:=\left(1+\sup_{x\in \overline{B(0,R)}} \left|\vb(x)\right|\right)^2$.
Similarly, consider the straight line connecting $y$ to $y'$ given by
\[
\xi(s): = y +\frac{s(y'-y)}{c}, \qquad s \in [0, c],
\]
and
\begin{equation} \label{eqn:xi}
A_{c}[\xi]= \frac{1}{4} \int_0^{c} \left|\frac{y'-y}{|y'-y|}-\vb\left(\xi(s)\right)\right|^2 \, ds \leq Mc.
\end{equation}
Now consider $\tilde{\gamma}:[0, T+a+c] \to \mathbb{R}^d$ defined by
\[
\tilde{\gamma}(s)=\left\{\begin{aligned}
&\eta(s), \qquad \qquad \quad \,\,\,\,\, s \in [0,a],\\
&\gamma(s-a), \qquad \qquad s\in[a, T+a],\\
&\xi(s-T-a), \qquad \, s\in[T+a, T+a+c].
\end{aligned}\right.
\]
Then $\tilde{\gamma}(0)=x'$ and $\tilde{\gamma}(T+a+c)=y'$, and by \eqref{eqn:eta} and \eqref{eqn:xi},
\[
    h_{T+a+c}(x',y') \leq A_{T+a+c} [\tilde{\gamma}] \leq h_T(x,y)+M(a+c)=h_T(x,y) +M\left(|x-x'|+|y-y'|\right).
\]
Taking the liminf as $T\to\infty$, we obtain
\[
\begin{aligned}
h(x',y')= \liminf_{T \to \infty} h_{T+a+c}(x',y') &\leq \liminf_{T \to \infty} h_T (x,y) +M\left(|x-x'|+|y-y'|\right)\\&=h (x,y) +M\left(|x-x'|+|y-y'|\right).
\end{aligned}
\]
Exchanging $(x,y)$ with $(x',y')$ gives
\[
\left|h(x,y)-h(x',y')\right|\leq M\left(|x-x'|+|y-y'|\right)
\]
Thus, \(h\) is locally Lipschitz on
\(\mathbb{R}^d\times\mathbb{R}^d\).
\end{proof}

\begin{rem}
In general, $h$ is not globally Lipschitz on $\mathbb R^d \times \mathbb R^d$. If $U(x)=|x|^2$ and $\vb(x)=-\nabla U(x)=-2x$, then (see Proposition~\ref{prop:pbchar})
\[
h(0, x)=|x|^2,
\]
and $h$ is not globally Lipschitz on $\mathbb R^d\times \mathbb R^d$.
\end{rem}

\subsection{Characterization of Peierls Barrier for the case $\vb = -\nabla U$}

 In this section, we characterize the explicit formula for Peierls Barrier in the case $\vb = -\nabla U$. This gradient case gives a concrete example, but our main results do not rely on this characterization.

\begin{prop}\label{prop:pbchar}
Assume (A1). Suppose further that \(\vec{\ell}\equiv 0\) and that \(U:\mathbb{R}^d\to\mathbb{R}\) has a unique critical point \(x_1\), which is also its unique global minimum. Then for any $y \in \mathbb{R}^d$,
\[
h(x_1,y)=U(y)-U(x_1).
\]
\end{prop}

\begin{proof}
\textbf{Step 1}: $h(x_1,y) \geq U(y) - U(x_1)$.

\noindent
For any $\gamma \in \mathrm{AC}([0,T];\mathbb{R}^d)$ with $\gamma(0)=x_1$ and $\gamma(T)=y$, we have
\begin{equation}
\begin{aligned}
 &\int_0^T L(\dot{\gamma}(t),\gamma(t)) \, dt
 =\int_0^T \frac{1}{4}|\dot{\vg}(t)+ \nabla U(\gamma(t))|^2 \,dt \\
  =& \int_0^T \left(
  \frac14|\dot{\vg}(t)-\nabla U(\gamma(t))|^2
  +\dot{\vg}(t)\cdot\nabla U(\gamma(t))\right)\,dt\\
  \geq& \int_0^T \frac{d}{d t} U(\vg(t)) \,dt
  = U(y) - U(x_1).
\end{aligned}
\end{equation}
Taking the infimum over all admissible paths, and then passing to the liminf as $T\to \infty$, yields $h(x_1, y) \geq U(y) - U(x_1)$.

\textbf{Step 2}: $h(x_1, y) \leq U(y) - U(x_1)$.

\noindent
If $y=x_1$, the constant path gives the desired upper bound. Hence, assume
$y\neq x_1$ and consider the initial value problem
\begin{equation}
\left\{
\begin{aligned}
    &\dot{\gamma}(t)=-\nabla U(\gamma(t)),\\
    &\gamma(0)=y.
\end{aligned}
\right.
\end{equation}
Since $U\in C^1$, the vector field $-\nabla U$ is continuous. Peano's
existence theorem therefore gives at least one local solution $\gamma$. Along any such solution,
\begin{equation}\label{eqn:Udecrease}
\frac{d}{dt} U\bigl(\gamma(t)\bigr)
= \nabla U(\gamma(t)) \cdot \dot{\gamma}(t)
= -|\nabla U(\gamma(t))|^2 \leq 0,
\end{equation}
so $U(\gamma(t)) \leq U(y)$ for all $t$ in its interval of existence. By the coercivity assumption, there exists $R>0$ such that $U(x) > U(y)$ whenever $|x|>R$. This implies $|\gamma(t)| \leq R$. Thus the solution remains in the compact set $\{U\leq U(y)\}$, on which $\nabla U$ is bounded. The standard continuation argument for ODEs with continuous vector fields then extends this chosen solution to all positive times, so
$\gamma:[0,\infty)\to\mathbb{R}^d$.
Set $M:=1+\max_{z\in\overline{B(0,R)}}|\nabla U(z)|$.
Since $\gamma([0,\infty))\subset\overline{B(0,R)}$, we have
\begin{equation}\label{eqn:gavelbd}
|\dot{\gamma}(t)| = |\nabla U(\gamma(t))| \leq M
\qquad\text{for all }t\ge 0.
\end{equation}
By \eqref{eqn:Udecrease} and the fact that $U(\gamma(t))\ge U(x_1)$, there exists $C_0\in\mathbb{R}$ such that
\begin{equation}\label{eqn:Uconverge}
\lim_{t\to\infty} U(\gamma(t))=C_0.
\end{equation}
Let $\{t_n\}\subset [0,\infty)$ with $t_n\to\infty$. Since $\gamma([0,\infty))\subset \overline{B(0,R)}$, there exist a subsequence $\{t_{n_k}\}$ and a point $x^\ast\in \overline{B(0,R)}$ such that
$\gamma(t_{n_k})\to x^\ast$. By continuity of $U$, we have $U(x^\ast)=C_0$. We claim that $x^\ast=x_1$. Suppose not. Since $x_1$ is the unique critical point of $U$, we have $\nabla U(x^\ast)\neq 0$. By continuity of $\nabla U$, there exist $\delta>0$ and $\epss_0>0$ such that
\[
|\nabla U(x)|\ge \epss_0
\qquad\text{for all }x\in B(x^\ast,\delta).
\]
Since $\gamma(t_{n_k})\to x^\ast$, there exists $K_0\in\mathbb{N}$ such that $\gamma(t_{n_k})\in B\left(x^\ast,\frac{\delta}{2}\right)$ for all $k\ge K_0$. For any $k\ge K_0$, by \eqref{eqn:gavelbd}, for $\tau := \frac{\delta}{4M}$, we have
\[
|\gamma(t)-\gamma(t_{n_k})|\le M(t-t_{n_k})< \frac{\delta}{2}
\qquad\text{for all }t\in [t_{n_k},\,t_{n_k}+\tau],
\]
and hence $\gamma(t)\in B(x^\ast,\delta)$ for all $t\in [t_{n_k},\,t_{n_k}+\tau]$.
Therefore, by \eqref{eqn:Udecrease},
\[
U(\gamma(t_{n_k}+\tau)) - U(\gamma(t_{n_k}))
= -\int_{t_{n_k}}^{t_{n_k}+\tau} |\nabla U(\gamma(t))|^2\,dt
\le -\epss_0^2 \tau.
\]
Thus
\[
U(\gamma(t_{n_k}+\tau))
\le U(\gamma(t_{n_k})) - \epss_0^2 \tau.
\]
Since $\gamma(t_{n_k})\to x^\ast$, we have $U(\gamma(t_{n_k}))\to C_0$. On the other hand, $t_{n_k}+\tau\to\infty$, so \eqref{eqn:Uconverge} also gives $U(\gamma(t_{n_k}+\tau))\to C_0$. This contradicts the preceding uniform decrease by $\epss_0^2\tau$. Therefore every accumulation point of $\gamma(t)$ as $t\to\infty$ is $x_1$. Since the trajectory is precompact, it follows that $\lim_{t\to \infty} \gamma(t)=x_1$.

Let $\epss>0$. Choose $T_\epss>0$ such that
$|\gamma(T)-x_1|<\epss$ for all $T\geq T_\epss$.
For any $T \geq T_\epss$, define $\tilde{\gamma}:[0,T+\epss] \to \mathbb{R}^d$ by
\[
\tilde{\gamma}(s):=\begin{cases}
x_1+\frac{s}{\epss}(\gamma(T)-x_1), & s\in [0,\epss],\\
\gamma(T+\epss-s), & s \in [\epss, T+\epss].
\end{cases}
\]
Note that $\tilde{\gamma}(0)=x_1$, $\tilde{\gamma}(T+\epss)=y$.
Since $x_1,\gamma(T)\in\overline{B(0,R)}$, the connecting
segment also lies in this ball, so $|\nabla U(\tilde\gamma(s))|\leq M$
for $s\in[0,\epss]$.
Hence, by \eqref{eqn:Udecrease} and
$U(\gamma(T))\geq U(x_1)$, for any $T\geq T_\epss$, we have
\[
\begin{aligned}
h_{T+\epss}(x_1, y) &\leq \int_0^{T+\epss} L\left(\dot{\tilde{\gamma}}(t), \tilde{\gamma}(t)\right)\,dt\\
&= \int_0^\epss \frac{1}{4}\left|\frac{\gamma(T)-x_1
}{\epss}+\nabla U(\tilde{\gamma}(t))\right|^2\,dt
+\int_0^T |\nabla U(\gamma(t))|^2\,dt\\
&\leq \frac{(1+M)^2}{4}\epss + U(y)-U\left(\gamma(T)\right)\\
& \leq \frac{(1+M)^2}{4}\epss +U(y)-U(x_1).
\end{aligned}
\]
For a fixed $\epss>0$, we obtain
\[
h(x_1,y)=\liminf_{T\to \infty}h_{T+\epss}(x_1,y)
\leq U(y)-U(x_1)+\frac{(1+M)^2}{4}\epss.
\]
Letting $\epss \to 0$, we conclude that
$h(x_1, y)\leq U(y)-U(x_1)$.
\end{proof}

\bigskip

\section{Characterization of Aubry Set}\label{sec:aubry}
In this section, we characterize the Aubry set under assumption (A1). In particular, when $\vb$ admits the orthogonal decomposition specified in (A1), the Aubry set coincides with the set of critical points of $U$.

We first recall the definition of the Aubry set and then characterize it in the proposition that follows.

\begin{defn}
Assume (A1). The Aubry set $\mathcal{A} \subset \mathbb{R}^d$ is defined by
\[
\mathcal{A}:=\left\{x\in \mathbb{R}^d:h(x,x)=0\right\}.
\]
\end{defn}

\begin{prop}
\label{prop:aubry-characterization}
Assume (A1). Then the Aubry set coincides with the set of critical points of \(U\), namely,
\begin{equation}
    \mathcal{A} = \{ x \in \bR^d : \nabla U(x) = 0 \}= \{ x \in \bR^d : \vb(x) = 0 \}.
\end{equation}
\end{prop}

\begin{proof}
    Let $B:=\{ x \in \bR^d : \nabla U(x) = 0 \}$. We proceed in two steps to prove $\mathcal{A} = B$.

Step 1. We show $B \subseteq \mathcal{A}$. For any $x \in B$, we have $\nabla U(x) = 0$, hence $\vb(x) = 0$ by (A1). Let $T >0$ and consider the constant path $\eta(t) \equiv x$ for $t \in [0,T]$. Then
\begin{equation}
0\leq h_T(x,x)\leq\int_0^TL(\dot{\eta}(t), \eta(t)) \,dt = \frac{1}{4} \int_0^T \left| 0 - \vb(x) \right|^2 \, dt= 0,
\end{equation}
and $h(x, x)= \liminf_{T \to  \infty} h_T(x,x)= 0$. Hence, $x \in \mathcal{A}$ and $B \subseteq \mathcal{A}$.

Step 2. We now prove that \( \mathcal{A} \subset B \). Suppose, for contradiction, there exists \( y \in \mathcal{A} \) such that \( \vec{b}(y) \neq 0 .\) Then $\nabla U(y)\neq 0$. By continuity of $\nabla U$ and $\vb$, there exist constants $\delta,r,C_r>0$ such that
\begin{equation}\label{eqn:bulbd}
0<\delta<\left|\nabla U(x)\right|\leq\left|\vb(x)\right|\leq C_r
\quad \text{for all }x\in B(y,r),
\end{equation}
where $|\vb|\geq|\nabla U|$ follows from the orthogonal decomposition in
(A1). Since $y \in \mathcal{A}$, by Proposition~\ref{prop:minimizer-regularity} and the definition of Aubry set, there exists a sequence of paths $\left\{\gamma_k \right\}\subset \mathrm{AC}\left([0,t_k];\bR^d\right)$ with $t_k\to\infty$ such that $\gamma_k(0)=\gamma_k(t_k)=y$ and
\begin{equation}\label{eqn:limvl}
\lim_{k \to \infty}\int_0^{t_k}\left|\dot{\gamma}_k(s)-\vb\left(\gamma_k(s)\right)\right|^2 \,ds =0.
\end{equation}
Since $\gamma_k\left(0\right)=\gamma_k({t_k}) =y$ and $\nabla U \cdot \vl=0$, by the fundamental theorem of calculus,
\[
\begin{aligned}
0&=U(\gamma_k(t_k))-U(\gamma_k(0))=\int_0^{t_k}\nabla U(\gamma_k(s))\cdot \dot{\gamma}_k(s)\, ds\\
&= \int_0^{t_k}\nabla U(\gamma_k(s))\cdot \left(\dot{\gamma}_k(s)-\vb(\gamma_k(s))\right)\, ds - \int_0^{t_k} \left|\nabla U (\gamma_k(s))\right|^2 ds,
\end{aligned}
\]
which implies
\[
\begin{aligned}
\int_0^{t_k} \left|\nabla U (\gamma_k(s))\right|^2 ds&=\int_0^{t_k}\nabla U(\gamma_k(s))\cdot \left(\dot{\gamma}_k(s)-\vb(\gamma_k(s))\right)\, ds\\
&\leq \left(\int_0^{t_k} \left|\nabla U (\gamma_k(s))\right|^2 \, ds \right)^\frac{1}{2} \left(\int_0^{t_k}\left|\dot{\gamma}_k(s)-\vb(\gamma_k(s))\right|^2 \,ds\right)^\frac{1}{2},
\end{aligned}
\]
where the last inequality follows from H\"{o}lder's inequality.
Combining with~\eqref{eqn:limvl}, we obtain
\begin{equation}\label{eqn:bgul2bd}
\int_0^{t_k} \left|\nabla U (\gamma_k(s))\right|^2 \, ds \leq \int_0^{t_k}\left|\dot{\gamma}_k(s)-\vb(\gamma_k(s))\right|^2 \,ds \to 0, \qquad \text {as}\,\, k \to \infty.
\end{equation}
For each $k$, define
\[
\tau_k:=\inf\left\{s\in[0,t_k]:\gamma_k(s)\in\partial B(y,r)\right\},
\]
with the convention that $\tau_k=t_k$ if this set is empty. We claim that
$\tau_k\geq c>0$ for all sufficiently large $k$, where $c$ is independent
of $k$.

If $\gamma_k$ does not reach $\partial B(y,r)$, then $\tau_k=t_k$, and the
claim follows from $t_k\to\infty$. If $\gamma_k$ reaches the boundary, then
\eqref{eqn:limvl} implies that, for all sufficiently large $k$,
\[
\begin{aligned}
r
&=|\gamma_k(\tau_k)-\gamma_k(0)|
\leq \int_0^{\tau_k}|\dot\gamma_k(s)|\,ds\\
&\leq
\left(\int_0^{t_k}
|\dot\gamma_k(s)-\vb(\gamma_k(s))|^2\,ds\right)^{1/2}
\tau_k^{1/2}
+\int_0^{\tau_k}|\vb(\gamma_k(s))|\,ds\\
&\leq \tau_k^{1/2}+C_r\tau_k.
\end{aligned}
\]
Since $s^{1/2}+C_rs\to0$ as $s\downarrow0$, this inequality gives a
constant $c>0$, independent of $k$, such that $\tau_k\geq c$ for all
sufficiently large $k$.

By the definition of $\tau_k$, the path remains in $B(y,r)$ for
$0\leq s<\tau_k$. Using \eqref{eqn:bulbd}, we therefore obtain
\[
\int_0^{t_k} \left|\nabla U(\gamma_k(s))\right|^2 \, ds \geq \int_0^{\tau_k} \left|\nabla U(\gamma_k(s))\right|^2 \, ds \geq \delta^2\tau_k \geq \delta^2c >0,
\]
which contradicts \eqref{eqn:bgul2bd}, and thus completes the proof.
\end{proof}

We next show that the Ma\~n\'e potential and the Peierls
barrier coincide when one endpoint belongs to $\mathcal A$. For equality
on the Aubry set in the compact setting, see
\cite[Remark~6.28]{sorrentino2010mather}.

\begin{prop}\label{prop:peierls-mane-aubry}
Assume (A1). Let $x \in \mathcal{A}$ and $y \in \mathbb{R}^d$. Then $h(x, y)=V(x,y)$ and $h(y, x)=V(y,x)$.
\end{prop}

\begin{proof}
By Proposition~\ref{prop:aubry-characterization}, $\vb(x)=0$, so $L(0,x)=0$.
Fix $0<s<T$. By Proposition~\ref{prop:minimizer-regularity}(c), choose
minimizers $\gamma,\sigma\in\mathrm{AC}([0,s];\mathbb R^d)$ for
$h_s(x,y)$ and $h_s(y,x)$, respectively. Define
\[
\eta_T(t)=
\begin{cases}
x, & 0\leq t\leq T-s,\\
\gamma(t-T+s), & T-s\leq t\leq T,
\end{cases}
\qquad
\xi_T(t)=
\begin{cases}
\sigma(t), & 0\leq t\leq s,\\
x, & s\leq t\leq T.
\end{cases}
\]
These curves are absolutely continuous on $[0,T]$, with
$\eta_T(0)=x$, $\eta_T(T)=y$, $\xi_T(0)=y$, and $\xi_T(T)=x$. Thus
\[
\begin{aligned}
h_T(x,y)&\leq A_T[\eta_T]
=A_s[\gamma]+(T-s)L(0,x)=h_s(x,y),\\
h_T(y,x)&\leq A_T[\xi_T]
=A_s[\sigma]+(T-s)L(0,x)=h_s(y,x).
\end{aligned}
\]
Therefore, both $T\mapsto h_T(x,y)$ and $T\mapsto h_T(y,x)$ are
nonincreasing and bounded below by zero. Their limits as $T\to\infty$
equal their infima over $T>0$, so the definitions of $h$ and $V$ give
\[
\begin{aligned}
h(x,y)&=\lim_{T\to\infty}h_T(x,y)=\inf_{T>0}h_T(x,y)=V(x,y),\\
h(y,x)&=\lim_{T\to\infty}h_T(y,x)=\inf_{T>0}h_T(y,x)=V(y,x).
\end{aligned}
\]
\end{proof}

\bigskip

\section{Construction of Weak KAM solutions}\label{sec:weakKAM}
In this section, we develop a general construction of weak KAM solutions to \eqref{eq:HJE0}. We first show that coercive weak KAM solutions are uniquely determined by their values on the Aubry set. Next, using the Peierls barrier, we construct a global energy landscape $W:\mathbb{R}^d\to\mathbb{R}$ as in \eqref{eqn:Wextension}, with boundary values on the Aubry set satisfying the discrete weak KAM problem, and prove that \(W\) is a weak KAM solution of negative type. We then establish its invariance under the Lax--Oleinik semigroup and uniqueness among coercive weak KAM solutions with the prescribed values on the Aubry set. Finally, we compare $W$ with sub- and supersolutions of the time-dependent HJE \eqref{eq:limit-HJ}.

\subsection{Weak KAM solutions, calibrated curves, and Aubry set}
We begin by recalling the definition of a weak KAM solution of negative type.

\begin{defn}[Weak KAM solution of negative type]
\label{defn:weakKAM_neg}
A function $u \in C(\bR^d)$ is a weak KAM solution of negative type to the stationary HJE \eqref{eq:HJE0} if:
\begin{enumerate}[(I)]
    \item $u$ is \emph{dominated by $L$} ($u \prec L$): for any absolutely continuous curve $\gamma:[a,b]\to\bR^d$,
    \begin{equation}
        u(\gamma(b)) - u(\gamma(a)) \leq \int_a^b L(\dot{\gamma}(t),\gamma(t))\,d t.
    \end{equation}
    \item For every $x \in \bR^d$, there exists a continuous piecewise $C^1$ curve $\gamma:(-\infty,0]\to\bR^d$ with $\gamma(0)=x$ such that for all $a < b \leq 0$,
    \begin{equation}
        u(\gamma(b)) - u(\gamma(a)) = \int_a^b L(\dot{\gamma}(t),\gamma(t))\,dt.
    \end{equation}
    Such a $\gamma$ is called a \emph{calibrated curve} (backward characteristic).
\end{enumerate}
\end{defn}

Throughout this paper, the term weak KAM solution refers to a weak KAM solution of negative type.

We first prove that $U$ is a weak KAM solution.

\begin{prop}\label{prop:U-weak-KAM}
Assume (A1) and (A3). Then $U$ is a coercive weak KAM
solution of negative type to the stationary HJE \eqref{eq:HJE0}.
\end{prop}

\begin{proof}
We proceed in two steps to prove that $U$ is a weak KAM solution.

Step 1. We show $U\prec L$. For any
$\eta\in\mathrm{AC}([a,b];\bR^d)$, integrating
\eqref{eq:U-calibration-identity} gives
\[
 U(\eta(b))-U(\eta(a))
 \leq\int_a^b L(\dot\eta(t),\eta(t))\,dt.
\]
Hence, $U\prec L$.

Step 2. For any $x\in\bR^d$, let $\gamma$ solve
\begin{equation}\label{eq:U-calibrated-ode}
 \dot\gamma(t)=\nabla U(\gamma(t))+\vl(\gamma(t)),
 \qquad \gamma(0)=x.
\end{equation}
By (A3), this equation has a unique local solution. Since
$\nabla U\cdot\vl=0$, we have
\[
\begin{aligned}
 \frac{d}{dt}U(\gamma(t))
 &=\nabla U(\gamma(t))\cdot
   \bigl(\nabla U(\gamma(t))+\vl(\gamma(t))\bigr)\\
 &=|\nabla U(\gamma(t))|^2\geq0.
\end{aligned}
\]
Hence, $U(\gamma(t))\leq U(x)$ for every $t\leq0$ in the interval
of existence. By (A1), the set
\[
 K:=\{y\in\bR^d:U(y)\leq U(x)\}
\]
is compact. Since $\gamma(t)\in K$ for all such $t$, the solution
extends to $(-\infty,0]$. By \eqref{eq:U-calibrated-ode} and
\eqref{eq:U-calibration-identity}, for any $a<b\leq0$, we have
\[
\begin{aligned}
 U(\gamma(b))-U(\gamma(a))
 &=\int_a^b\nabla U(\gamma(t))\cdot\dot\gamma(t)\,dt\\
 &=\int_a^b L(\dot\gamma(t),\gamma(t))\,dt.
\end{aligned}
\]
Thus $\gamma$ is a backward calibrated curve for $U$ with $\gamma(0)=x$.
By Definition~\ref{defn:weakKAM_neg}, $U$ is a weak KAM solution of
negative type, and $U$ is coercive by (A1).
\end{proof}

In the next lemma, we establish the connection between calibrated curves and the Aubry set. In particular, we show that for any calibrated curve of a weak KAM solution that remains in a compact set, its \(\alpha\)-limit set is contained in the Aubry set.

\begin{defn}
For $\gamma: (-\infty, 0] \to \mathbb{R}^d$, define the $\alpha$-limit set of $\gamma$ by
\[
\alpha(\gamma):=\left\{z:\exists \,t_n \to -\infty, \gamma(t_n) \to z\right\}.
\]

For $\eta: [0, \infty) \to \mathbb{R}^d$, define the $\omega$-limit set of $\eta$ by
\[
\omega(\eta):=\left\{z:\exists \,t_n \to \infty, \eta(t_n) \to z\right\}.
\]
\end{defn}

\begin{lem}\label{lem:alphaset}
    Suppose $w \in C(\mathbb{R}^d)$ is a weak KAM solution of negative type to the stationary HJE, and $\gamma:(-\infty, 0] \to \mathbb{R}^d$ is a calibrated curve for $w$ such that $\gamma((-\infty, 0]) \subset K$ for some compact set $K \subset \mathbb{R}^d$. Then $\alpha(\gamma) \subset \mathcal{A}$.
\end{lem}

\begin{proof}
    Let $z \in \alpha(\gamma)$. Then there exists a sequence $\left\{t_n\right\}$ with $t_n \to -\infty$ such that $\lim_{n \to \infty} \gamma(t_n)=z$. Choose a subsequence $\{t_{n_k}\}$ such that $t_{n_{k+1}} -t_{n_{k}} <-k$, and define $T_k=t_{n_k}-t_{n_{k+1}}$, so that $T_k\to\infty$. Then,
\[
h_{T_k}(\gamma(t_{n_{k+1}}),\gamma(t_{n_k}))\leq \int_{t_{n_{k+1}}}^{t_{n_k}} L(\dot{\gamma}(\tau), \gamma(\tau)) d\tau=w(\gamma(t_{n_k}))-w(\gamma(t_{n_{k+1}})).
\]
Let $\xi_k:[t_{n_{k+1}}, t_{n_k}] \to \mathbb{R}^d$ be a minimizer for $h_{T_k}(\gamma(t_{n_{k+1}}),\gamma(t_{n_k}))$, that is,
\[
h_{T_k}(\gamma(t_{n_{k+1}}),\gamma(t_{n_k}))= \frac{1}{4}\int_{t_{n_{k+1}}}^{t_{n_k}} \left|\dot{\xi_k}(s)-\vb(\xi_k(s))\right|^2\,ds, \quad \xi_k(t_{n_{k+1}})=\gamma(t_{n_{k+1}}), \quad \xi_k(t_{n_k})=\gamma(t_{n_{k}}).
\]
Fix $0<\delta\leq1$. For $k \in \mathbb{N}$ sufficiently large, we have
\begin{equation}\label{eqn:gammayclose}
|\gamma(t_{n_k}) - z| \leq \delta.
\end{equation}
Define $\tilde{\xi}_k:[-T_k-2\delta, 0] \to \mathbb{R}^d$ by
\begin{equation}
    \tilde{\xi}_k(s):=\left\{\begin{aligned}&z + \frac{s-(-T_k-2\delta)}{\delta}(\gamma(t_{n_{k+1}})-z), \qquad \quad\qquad -T_k-2\delta \leq s \leq -T_k-\delta, \\
 &\xi_k(s+t_{n_k}+\delta),\qquad \qquad\qquad\qquad\qquad\qquad\qquad -T_k-\delta \leq s \leq -\delta,\\
    & \gamma (t_{n_k}) +\frac{s+\delta}{\delta} (z-\gamma(t_{n_{k}})), \qquad \qquad \qquad \qquad \,\,\, -\delta \leq s \leq 0.
    \end{aligned}
    \right.
\end{equation}
Then $\tilde{\xi}_k$ is an admissible path for $h_{T_k+2\delta}(z, z)$. Hence, for $k$ sufficiently large,
\[
\begin{aligned}
 h_{T_k+2\delta}(z, z) &\leq \int_{-T_k-2\delta}^0 \frac{1}{4} \left|\dot{\tilde{\xi}}_k(s)-\vb(\tilde{\xi}_k(s))\right|^2 \,ds\leq C\delta+ \frac{1}{4}\int_{t_{n_{k+1}}}^{t_{n_k}} \left|\dot{\xi_k}(s)-\vb(\xi_k(s))\right|^2\,ds\\
 &\leq C\delta +h_{T_k}(\gamma(t_{n_{k+1}}),\gamma(t_{n_k}))\leq C\delta+w(\gamma(t_{n_k}))-w(\gamma(t_{n_{k+1}})),
\end{aligned}
\]
where the second inequality follows from \eqref{eqn:gammayclose}, the fact that $\gamma ((-\infty,0
]) \subset K$, and the continuity of $b$.
Taking liminf as $k \to \infty$, we obtain
\[
0\leq h(z, z)\leq \liminf_{k\to \infty} h_{T_k+2\delta}(z, z)\leq C\delta .
\]
Letting $\delta\downarrow0$, we have $h(z, z)=0$, and hence $z \in \mathcal{A}$.
\end{proof}

We now prove a comparison principle for weak KAM solutions to
\eqref{eq:HJE0}, assuming coercivity of the upper solution.

\begin{prop}\label{AB}
Let $V_1, V_2$ be two weak KAM solutions of negative type to \eqref{eq:HJE0}. Assume that $\lim_{|x| \to \infty} V_2(x)=+\infty$. If $V_1 \leq V_2$ on $\sA$, then $V_1 \leq V_2$ on $\bR^d$.
\end{prop}

\begin{proof}
Fix $x\in\mathbb R^d$. Since $V_2$ is a backward weak KAM solution, there exists a backward calibrated curve $\gamma:(-\infty,0]\to\mathbb R^d$ such that $\gamma(0)=x$ and for all $t<0$,
\begin{equation}\label{eqn:v1L}
V_2(x)-V_2(\gamma(t))
=
\int_t^0 L(\dot\gamma(s),\gamma(s))\,ds.
\end{equation}
Since $L\ge 0$, we have $V_2(\gamma(t))\le V_2(x)$ for all $t<0$. Therefore,
\[
\gamma((-\infty,0])\subset \{z\in\mathbb R^d:\ V_2(z)\le V_2(x)\}=:K.
\]
The set $K$ is compact by the continuity and coercivity of
$V_2$. By Lemma~\ref{lem:alphaset}, $\alpha(\gamma)\subset\mathcal A$.

On the other hand, since $V_1$ is dominated by $L$, for all $t<0$,
\[
V_1(x)-V_1(\gamma(t))
\le
\int_t^0 L(\dot\gamma(s),\gamma(s))\,ds.
\]
Combining domination of $V_1$ with calibration of $V_2$
in \eqref{eqn:v1L}, we obtain
\begin{equation}\label{eqn:V1V2}
V_1(x)-V_2(x)
\le
V_1(\gamma(t))-V_2(\gamma(t)),
\qquad\forall\, t<0.
\end{equation}
Since $\gamma((-\infty,0])\subset K$, there exists a sequence $t_n\to-\infty$ and a point $a\in\alpha(\gamma)\subset \mathcal A$ such that
$\gamma(t_n)\to a$. Passing to the limit along $t_n$ in \eqref{eqn:V1V2} and using the continuity of $V_1,V_2$, we obtain
\[
V_1(x)-V_2(x)
\le
V_1(a)-V_2(a)
\le 0, \qquad \text{for all } x\in \mathbb{R}^d,
\]
since $V_1\le V_2$ on $\mathcal{A}$ by assumption. Hence, $V_1\le V_2$.
\end{proof}

\subsection{Construction of Weak KAM Solutions via the Aubry Set and the Peierls Barrier}

We now present a general construction of weak KAM solutions to \eqref{eq:HJE0} based on the Aubry set and the Peierls barrier. Specifically, we prescribe the values of the solution on the Aubry set by requiring them to satisfy the discrete weak KAM problem, and then extend these values to the whole space \(\mathbb{R}^d\) using the Peierls barrier.

\begin{defn}
Let $W:\mathcal{A} \to \mathbb{R}$ satisfy
\begin{equation}\label{eqn:dwkam}
W(y)-W(z) \leq h(z,y)=V(z,y),
\end{equation}
for any $y,z \in \mathcal{A}$. Its extension to $\mathbb{R}^d$ is defined by
\begin{equation}\label{eqn:Wextension}
W(x) := \min_{y \in \mathcal{A}} \left\{{ W(y) + h(y, x) }\right\}=\min_{y \in \mathcal{A}} \left\{{ W(y) + V(y, x) }\right\}, \qquad x \in \mathbb{R}^d.
\end{equation}
We say $W$ satisfies discrete weak KAM problem on $\sA$ if \eqref{eqn:dwkam} holds.
\end{defn}

\begin{rem}\label{rem:wlip}
Note that \eqref{eqn:dwkam} is equivalent to the condition
\[
W(y) := \min_{z \in \mathcal{A}} \left\{{ W(z) + h(z, y) }\right\}, \qquad y \in \mathcal{A}.
\]

Since $\mathcal{A}$ consists of finitely many points and $h$ is locally Lipschitz by Proposition~\ref{prop:minimizer-regularity}, the extension $W$ defined by \eqref{eqn:Wextension} is locally Lipschitz on $\mathbb{R}^d$.
\end{rem}

We now prove several properties of the function $W$ defined by \eqref{eqn:Wextension}.

\begin{prop}
\label{prop:calibrated}
Assume (A1) and (A3). Let \(\{W_y\}_{y\in\mathcal{A}}\subset\mathbb{R}\) satisfy
\[
    W_y-W_z\leq h(z,y), \qquad y,z\in\mathcal{A}.
\]
Let \(W:\mathbb{R}^d\to\mathbb{R}\) be the function defined by the Peierls barrier construction in \eqref{eqn:Wextension}, with boundary values satisfying the discrete weak KAM problem, namely,
\[
    W(y)=W_y, \qquad y\in\mathcal{A}.
\]
Then:
\begin{enumerate}[(i)]
    \item \(W\) is coercive, i.e., $\lim_{|x|\to\infty}W(x)=+\infty$.
    \item $W$ is dominated by $L$, that is, for every absolutely continuous curve $\gamma:[a,b]\to\mathbb{R}^d$,
    \[
    W(\gamma(b))-W(\gamma(a))\leq\int_a^b L(\dot\gamma(t),\gamma(t))\,dt.
    \]
    Moreover, for every $x\in\mathbb{R}^d$, there exists a calibrated curve $\gamma\in C^2((-infty,0];\mathbb{R}^d)$ for $W$ with $\gamma(0)=x$, which converges to a zero of $\vb$ as time tends to $-\infty$.
    \item Suppose, in addition, that \(\vb\in C^2\). Then the projected Mather set \(\mathcal{M}_0\) coincides with the Aubry set \(\mathcal{A}\).
\end{enumerate}
\end{prop}

\begin{proof}
(i) Let $x \in \mathbb{R}^d$, $C_1:=\min_{y\in\mathcal{A}}W(y)$, and $C_2:=\max_{y\in\mathcal{A}}U(y)$.
By the definition of $W$ and Proposition~\ref{prop:minimizer-regularity}, for some $y_i \in \mathcal{A}$,
\[
W(x) =W(y_i) +h (y_i, x) \geq C_1 +U(x)-U(y_i) \geq C_1 -C_2+U(x).
\]
Since $U(x)\to\infty$ as $|x|\to\infty$, we have $\lim_{|x|\to \infty} W(x) = \infty$.

(ii) We first prove that $W$ is dominated by $L$. Given $x,z\in\mathbb R^d$,
choose $y\in\mathcal A$ such that
\[
W(x)=W_y+V(y,x).
\]
The triangle inequality for the Ma\~n\'e potential gives
\[
W(z)
\leq W_y+V(y,z)
\leq W_y+V(y,x)+V(x,z)
=W(x)+V(x,z).
\]
If $\gamma:[a,b]\to\mathbb R^d$ is absolutely continuous with
$\gamma(a)=x$ and $\gamma(b)=z$, then
\[
V(x,z)\leq
\int_a^b L(\dot\gamma(t),\gamma(t))\,dt.
\]
Consequently,
\[
W(\gamma(b))-W(\gamma(a))
\leq \int_a^b L(\dot\gamma(t),\gamma(t))\,dt,
\]
so $W\prec L$.

Let $x \in \mathbb{R}^d$. Choose $y\in \mathcal{A}$ such that
\begin{equation}\label{eqn:wmane}
W(x)=W(y)+h(y, x)=W(y)+V(y, x).
\end{equation}
By the definition of $h(y, x)$ and Proposition~\ref{prop:minimizer-regularity}, there exists a sequence $T_n \to \infty$ and $C^2$ curves $\gamma_n:[-T_n, 0] \to \mathbb{R}^d$ such that
\[
h_{T_n}(y, x)=\int_{-T_n}^0L(\dot{\gamma}_n(t), \gamma_n(t))\,dt \to h(y, x), \qquad \gamma_n(-T_n)=y, \quad \gamma_n(0)=x.
\]
For $T_n>1$, consider the path $\xi_n:[-T_n,0] \to \mathbb{R}^d$ defined by
\[
\xi_n(t)=
\begin{cases}
y, & -T_n \le t \le -1,\\
y+(t+1)(x-y), & -1 \le t \le 0,
\end{cases}
\]
which is an admissible path for $h_{T_n}(y, x)$. Since $y\in\mathcal A$, we have $b(y)=0$. Hence,
\begin{equation}\label{eqn:gamnbd}
\begin{aligned}
\int_{-T_n}^0L(\dot{\gamma}_n(t), \gamma_n(t))\,dt
&=h_{T_n}(y, x)
\leq \int_{-T_n}^0L(\dot{\xi}_n(t),\xi_n(t))\,dt \\
&=\frac{1}{4}\int_{-T_n}^0 \left|\dot{\xi}_n(t)-\vb(\xi_n(t))\right|^2\,dt\\
&\leq  \frac{1}{2} \int_{-T_n}^0 \left|\dot{\xi}_n(t)\right|^2\,dt+\frac{1}{2} \int_{-T_n}^0 |\vb(\xi_n(t))|^2dt\\
&\leq \frac{1}{2}|x-y|^2+\frac{1}{2}\int_{-1}^0 | b(\xi_n(t))|^2dt+\frac{1}{2}\int_{-T_n}^{-1}|b(y)|^2\,dt\le C_0,
\end{aligned}
\end{equation}
where $C_0>0$ is independent of $n$. Here we use that $\vb$ is continuous and $\vb(y)=0$.

By Proposition~\ref{prop:minimizer-regularity}, we have
\[
\begin{aligned}
U(\gamma_n(t))-U(y)
&\le \int_{-T_n}^t \frac14|\dot\gamma_n(s)-b(\gamma_n(s))|^2\,ds\\
&=\int_{-T_n}^t L(\dot\gamma_n(s),\gamma_n(s))\,ds
\le \int_{-T_n}^0 L(\dot\gamma_n(s), \gamma_n(s))\,ds
\le C_0.
\end{aligned}
\]
Therefore,
\[
U(\gamma_n(t))\le U(y)+C_0=:M_0
\qquad\text{for all }t\in[-T_n,0].
\]
Since $U(z)\to\infty$ as $|z|\to\infty$, the set $
K:=\{z\in\mathbb R^d:\ U(z)\le M_0\}$ is compact, and $
\gamma_n([-T_n,0])\subset K$ for all $n$ sufficiently large.

Fix $T>0$. For all sufficiently large $n$, $[-T,0]\subset [-T_n,0]$. Since $K$ is compact and $b$ is continuous, there exists $M_1:=\sup_{z\in K}|b(z)|<\infty$.
Then
\[
\begin{aligned}
\int_{-T}^0|\dot{\gamma}_n(t)|^2\,dt
&\le 2\int_{-T}^0|\dot{\gamma}_n(t)-b(\gamma_n(t))|^2\,dt+2\int_{-T}^0|b(\gamma_n(t))|^2\,dt\\
&\le 8\int_{-T}^0L(\dot{\gamma}_n(t), \gamma_n(t))\,dt+2TM_1^2\le 8C_0+2TM_1^2.
\end{aligned}
\]
Hence for any $-T\le s<t\le0$, by Cauchy--Schwarz,
\[
|\gamma_n(t)-\gamma_n(s)| \le \int_s^t|\dot{\gamma}_n(\tau)|\,d\tau \le |t-s|^{1/2}\left(\int_s^t|\dot{\gamma}_n(\tau)|^2\,d\tau\right)^{1/2} \leq |t-s|^{1/2}(8C_0+2TM_1^2)^{1/2}.
\]
Thus $\{\gamma_n\}$ is equicontinuous on $[-T,0]$. Since $\gamma_n([-T,0])\subset K$, the family is also uniformly bounded on $[-T,0]$. The same bounds give boundedness in $W^{1,2}([-T,0];\mathbb R^d)$. By weak compactness and Arzelà--Ascoli, for each $T>0$ there exists a subsequence converging weakly in $W^{1,2}$ and uniformly on $[-T,0]$. Using a diagonal argument on the intervals $[-k,0]$, $k\in\mathbb N$, we may extract a subsequence, still denoted $\gamma_n$, and obtain a continuous curve $\gamma:(-\infty,0]\to K$ such that for every $T>0$,
\[
\gamma_n \rightharpoonup  \gamma \text{ weakly in } W^{1,2}([-T,0]; \mathbb{R}^d), \qquad \gamma_n \to \gamma \text{ uniformly on } [-T,0].
\]
Since $\gamma_n(0)=x$ for all $n \in \mathbb{N}$, we have $\gamma(0)=x$.

Fix $s<t<0$. For all sufficiently large $n \in \mathbb{N}$, one has $T_n>-s$. Since $\gamma_n$ is a minimizer on $[-T_n,0]$,
\begin{align}
h_{T_n}(y, x)
&=
h_{T_n+s}(y,\gamma_n(s))
+
h_{t-s}(\gamma_n(s),\gamma_n(t))
+
h_{-t}(\gamma_n(t),x).
\label{eqn:split}
\end{align}
Since $W(x)-W(y)=h(y, x)$ and $
h_{T_n}(y, x)\to h(y, x)$, we obtain
\begin{align}
&\Bigl[h_{T_n+s}(y,\gamma_n(s))-\bigl(W(\gamma_n(s))-W(y)\bigr)\Bigr]
\notag+
\Bigl[h_{t-s}(\gamma_n(s),\gamma_n(t))-\bigl(W(\gamma_n(t))-W(\gamma_n(s))\bigr)\Bigr]
\notag\\
&\qquad \qquad\qquad\qquad\qquad\qquad\qquad\qquad+
\Bigl[h_{-t}(\gamma_n(t),x)-\bigl(W(x)-W(\gamma_n(t))\bigr)\Bigr]
\longrightarrow 0.
\label{eqn:sum}
\end{align}
By domination of $W$, we have
\[
W(\gamma_n(s))-W(y)\le h_{T_n+s}(y,\gamma_n(s)),
\]
\[
W(\gamma_n(t))-W(\gamma_n(s))\le h_{t-s}(\gamma_n(s),\gamma_n(t)),
\]
\[
W(x)-W(\gamma_n(t))\le h_{-t}(\gamma_n(t), x).
\]
Therefore each of the three terms in \eqref{eqn:sum} must converge to $0$. In particular,
\begin{align}
h_{t-s}(\gamma_n(s), \gamma_n(t))
-
\bigl(W(\gamma_n(t))-W(\gamma_n(s))\bigr)
\longrightarrow 0.
\label{eqn:ts}
\end{align}

Now, since $\gamma_n(s)\to \gamma(s)$ and $\gamma_n(t)\to \gamma(t)$, and since the value function $(a,b)\mapsto h_{t-s}(a,b)$ is continuous, we have
\[
h_{t-s}(\gamma_n(s), \gamma_n(t))
\to
h_{t-s}(\gamma(s),\gamma(t)).
\]
Since $W$ is continuous by Remark \ref{rem:wlip}, we have
\[
W(\gamma_n(t))-W(\gamma_n(s))
\to
W(\gamma(t))-W(\gamma(s)).
\]
Passing to the limit in \eqref{eqn:ts}, we get
\begin{align}
h_{t-s}(\gamma(s),\gamma(t))
=
W(\gamma(t))-W(\gamma(s)).
\label{eq:W-equals-fixed-time-action}
\end{align}

On the other hand, by definition of $h_{t-s}$,
\begin{align}
h_{t-s}(\gamma(s),\gamma(t))
\le
\int_s^t L(\dot\gamma(\tau),\gamma(\tau))\,d\tau.
\label{eq:h-less-than-action}
\end{align}
We now prove the reverse inequality. Since each $\gamma_n|_{[s,t]}$ is minimizing,
\[
\int_s^t L(\dot\gamma_n(\tau),\gamma_n(\tau))\,d\tau
=
h_{t-s}(\gamma_n(s),\gamma_n(t)).
\]
Hence
\begin{align}
\int_s^t L(\dot\gamma_n(\tau),\gamma_n(\tau))\,d\tau
\to
h_{t-s}(\gamma(s), \gamma(t)).
\label{eq:actions-converge-to-fixed-time-value}
\end{align}

By the compactness argument, $\gamma_n\to\gamma$ uniformly on $[s,t]$, while $\dot\gamma_n\rightharpoonup \dot\gamma$ weakly in $L^2([s,t])$. Since $b$ is continuous,
\[
b(\gamma_n)\to b(\gamma)\qquad \text{strongly in }L^\infty([s,t]).
\]
Therefore
\[
\dot\gamma_n-b(\gamma_n)\rightharpoonup \dot\gamma-b(\gamma)
\qquad \text{weakly in }L^2([s,t]).
\]
By weak lower semicontinuity of the $L^2$ norm,
\begin{align}
\int_s^t L(\dot\gamma(\tau),\gamma(\tau))\,d\tau
\le
\liminf_{n\to\infty}
\int_s^t L(\dot\gamma_n(\tau),\gamma_n(\tau))\,d\tau.
\label{eq:lsc}
\end{align}
Combining \eqref{eq:actions-converge-to-fixed-time-value} and \eqref{eq:lsc}, we obtain
\begin{align}
\int_s^t L(\dot\gamma(\tau),\gamma(\tau))\,d\tau
\le
h_{t-s}(\gamma(s),\gamma(t)).
\label{eq:action-less-than-h}
\end{align}
Together with \eqref{eq:h-less-than-action}, this yields
\[
\int_s^t L(\dot\gamma(\tau),\gamma(\tau))\,d\tau
=
h_{t-s}(\gamma(s),\gamma(t)).
\]
Since $\gamma$ is a minimizer of $h_{t-s}(\gamma(s), \gamma(t))$, by Proposition~\ref{prop:minimizer-regularity}, $\gamma \in C^2([s,t];\mathbb{R}^d)$. Using \eqref{eq:W-equals-fixed-time-action}, we conclude that
\[
W(\gamma(t))-W(\gamma(s))
=
\int_s^t L(\dot\gamma(\tau),\gamma(\tau))\,d\tau.
\]
Since $s<t\le 0$ are arbitrary, $\gamma \in C^2((-\infty, 0];\mathbb{R}^d)$ is a backward calibrated curve for $W$.

Now since $\gamma((-\infty,0]) \subset K$, there exists a sequence $\left\{t_n\right\}$ with $t_n \to -\infty$ such that
\[
\lim_{n \to \infty} \gamma(t_n)=y^\ast
\]
for some $y^\ast \in K$. Hence, $y^\ast \in \alpha(\gamma) \subset \mathcal{A}$ by Lemma \ref{lem:alphaset}. By standard ODE theory, $\alpha(\gamma)$ is connected. Since $\mathcal{A}$ consists of finitely many isolated points, we must have $\alpha(\gamma)=\left\{y^\ast\right\}$. Therefore, $\lim_{t\to -\infty} \gamma(t)=y^\ast$.

\noindent
(iii) Let $\mathcal P_{\min}$ denote the probability
measures invariant under the Euler--Lagrange flow, supported on complete
trajectories, and minimizing $\int_{\mathbb R^{2d}}L(v,x)\,d\mu(x,v)$
among such invariant probability measures. The projected Mather set
is defined by
\[
 \mathcal M_0
 :=\operatorname{pr}_x\!\left(
 \overline{\bigcup_{\mu\in\mathcal P_{\min}}\operatorname{supp}\mu}
 \right),
 \qquad \operatorname{pr}_x(x,v)=x,
\]
where the closure is taken in $\mathbb R^d\times\mathbb R^d$.
We first prove $\mathcal A\subset\mathcal M_0$. Let $y \in \mathcal{A}$. Then $b(y)=0$, and the constant trajectory
\[
x(t)\equiv y, \qquad v(t)\equiv 0,
\]
is an orbit of the Euler-Lagrange flow. Define $\mu:=\delta_{(y,0)}$. Then $\mu$ is invariant and satisfies
\[
\int L\, d\mu=0.
\]
Since $L \geq 0$, $\mu$ is a minimizing invariant measure. By definition, the support of every minimizing invariant measure is contained in the Mather set, and thus $y \in \mathcal{M}_0$. Therefore, $\mathcal{A} \subset \mathcal{M}_0$.

For the reverse inclusion, let $\mu$ be any minimizing invariant probability measure. The minimum action is zero, so $\int L\,d\mu=0$ and $v=\vb(q)$ for $\mu$-almost every $(q,v)$. By invariance, along $\mu$-almost every trajectory $(q(t),v(t))$, we have $\dot q(t)=\vb(q(t))$ and, by (A1),
\[
\frac{d}{dt}U(q(t))=-|\nabla U(q(t))|^2\le0.
\]
Poincar\'e recurrence gives times $t_n\to\infty$ with $q(t_n)\to q(0)$ for $\mu$-almost every initial point. Monotonicity and continuity of $U$ then imply $U(q(t))=U(q(0))$ for every $t\ge0$, so the preceding identity gives $\nabla U(q(0))=0$. By continuity and Proposition~\ref{prop:aubry-characterization}, $\operatorname{supp}\mu\subset\mathcal A\times\{0\}$. Since $\mathcal A$ is closed, this proves $\mathcal M_0\subset\mathcal A$.
\end{proof}

We next verify that $W$ is a viscosity solution of \eqref{eq:HJE0} and is maximal among locally Lipschitz viscosity solutions with the same values on the Aubry set.
\begin{cor}
\label{cor:weakKAM}
Under the assumptions of Proposition~\ref{prop:calibrated}, let $W$ be the function constructed there. Then $W$ is a weak KAM solution of negative type to \eqref{eq:HJE0} and the maximal locally Lipschitz viscosity solution satisfying the prescribed values $W(y)=W_y$ for $y\in\mathcal A$.
\end{cor}

\begin{proof}
By Remark~\ref{rem:wlip}, $W$ is locally Lipschitz. Proposition~\ref{prop:calibrated}(ii) gives both the domination property $W\prec L$ and a backward calibrated curve through every point. Thus $W$ satisfies the two defining conditions of a weak KAM solution of negative type.

For completeness, we verify directly that $W$ is a viscosity solution of
\eqref{eq:HJE0}. Suppose that $\phi\in C^1$ touches $W$ from above at
$x_0$. For any $q\in\mathbb R^d$ and all sufficiently small $t>0$, consider the straight segment
$\gamma(s)=x_0+s q$, $-t\leq s\leq0$. Since $W-\phi$ has a local maximum
at $x_0$, we obtain
\[
\phi(x_0)-\phi(x_0-tq)
\leq W(x_0)-W(x_0-tq)
\leq \int_{-t}^0L(q,x_0+s q)\,ds.
\]
Dividing by $t$ and passing to the limit as $t\downarrow0$,
\[
D\phi(x_0)\cdot q\leq L(q,x_0),
\qquad q\in\mathbb R^d.
\]
Taking the supremum over $q$ gives
\[
H(D\phi(x_0),x_0)
=\sup_{q\in\mathbb R^d}\{D\phi(x_0)\cdot q-L(q,x_0)\}
\leq0.
\]

Now suppose that $\phi\in C^1$ touches $W$ from below at $x_0$, and let
$\gamma:(-\infty,0]\to\mathbb R^d$ be the backward calibrated curve with
$\gamma(0)=x_0$. For small
$t>0$,
\[
\int_{-t}^0L(\dot\gamma(s),\gamma(s))\,ds
=W(x_0)-W(\gamma(-t))
\leq \phi(x_0)-\phi(\gamma(-t)).
\]
The calibrated curve is $C^2$ by the proof of
Proposition~\ref{prop:calibrated}. Dividing by $t$ and letting $t\downarrow0$,
with $q_0:=\dot\gamma(0)$, yields
\[
L(q_0,x_0)\leq D\phi(x_0)\cdot q_0.
\]
Therefore
\[
H(D\phi(x_0),x_0)
=\sup_{q\in\mathbb R^d}\{D\phi(x_0)\cdot q-L(q,x_0)\}
\geq0.
\]
This proves the viscosity supersolution property and hence the viscosity
solution property.

Finally, let $v$ be any locally Lipschitz viscosity solution of
\eqref{eq:HJE0} such that $v(y)=W_y$ on $\mathcal A$. We show explicitly
that $v$ is dominated by $L$. Let
$\gamma:[a,b]\to\mathbb R^d$ be absolutely continuous and set
\[
K_1:=\{z\in\mathbb R^d:
\operatorname{dist}(z,\gamma([a,b]))\leq1\},
\]
which is compact. Since $v$ is locally Lipschitz, there exists a constant $M>0$ such
that $|Dv|\leq M$ a.e. on $K_1$, and
\begin{equation}\label{eq:viscosity-ae-subsolution}
H(Dv(z),z)\leq0
\qquad\text{for a.e. }z\in K_1.
\end{equation}

Let $\varrho_\delta$ be a standard mollifier supported in
$B_\delta(0)$ and define $v_\delta=v*\varrho_\delta$. For
$x\in\gamma([a,b])$ and $0<\delta<1$, convexity of $H$ in its first
variable and Jensen's inequality give
\[
\begin{aligned}
H(Dv_\delta(x),x)
&=H\left(\int_{\mathbb R^d}\varrho_\delta(x-z)Dv(z)\,dz,x\right)\\
&\leq\int_{\mathbb R^d}\varrho_\delta(x-z)H(Dv(z),x)\,dz.
\end{aligned}
\]
For the Hamiltonian $H(p,x)=|p|^2+p\cdot\vb(x)$,
\[
H(Dv(z),x)-H(Dv(z),z)
=Dv(z)\cdot\bigl(\vb(x)-\vb(z)\bigr).
\]
Define the modulus of continuity of $\vb$ on $K_1$ by
\[
\omega_{\vb,K_1}(\delta)
:=\sup_{\substack{x,z\in K_1\\|x-z|\leq\delta}}
|\vb(x)-\vb(z)|.
\]
Using \eqref{eq:viscosity-ae-subsolution} and the support of
$\varrho_\delta$, we obtain
\begin{equation}\label{eq:mollified-subsolution-error}
H(Dv_\delta(x),x)
\leq M\omega_{\vb,K_1}(\delta)=:r_\delta,
\qquad x\in\gamma([a,b]).
\end{equation}

Now since $v_\delta$ is smooth, for almost every $t\in[a,b]$,
\[
\begin{aligned}
\frac{d}{dt}v_\delta(\gamma(t))
&=Dv_\delta(\gamma(t))\cdot\dot\gamma(t) \leq H(Dv_\delta(\gamma(t)),\gamma(t))
+L(\dot\gamma(t),\gamma(t))\leq r_\delta+L(\dot\gamma(t),\gamma(t)).
\end{aligned}
\]
Integrating both sides, we obtain
\[
v_\delta(\gamma(b))-v_\delta(\gamma(a))
\leq \int_a^bL(\dot\gamma(t),\gamma(t))\,dt+(b-a)r_\delta.
\]
Since $v_\delta\to v$ uniformly on compact sets, letting
$\delta\downarrow0$ gives
\[
v(\gamma(b))-v(\gamma(a))
\leq \int_a^bL(\dot\gamma(t),\gamma(t))\,dt.
\]
Thus $v\prec L$. Taking the infimum over all times and paths from
$y\in\mathcal A$ to $x\in\mathbb R^d$ gives
\[
v(x)-v(y)\leq V(y,x)=h(y,x).
\]
Using $v(y)=W_y$ and taking the minimum over $y\in\mathcal A$, we conclude
that
\[
v(x)\leq\min_{y\in\mathcal A}\{W_y+h(y,x)\}=W(x).
\]
Thus $W$ is the maximal locally Lipschitz viscosity solution with the
prescribed boundary data.
\end{proof}

The domination and calibration properties also imply that weak KAM solutions are fixed points of the backward Lax--Oleinik semigroup.

\begin{cor}
\label{cor:laxoleinik}
Any weak KAM solution $w:\mathbb{R}^d \to \mathbb{R}$ is invariant under the
backward Lax--Oleinik semigroup: for every $t>0$ and $y\in\mathbb{R}^d$,
\begin{equation}
    w(y)=(S_t w)(y)
    =\inf_{x\in\mathbb{R}^d}\bigl\{w(x)+h_t(x,y)\bigr\}.
\end{equation}
Consequently, every weak KAM solution satisfies the free-time identity
\begin{equation}
    w(y)=\inf_{x\in\mathbb{R}^d}\bigl\{w(x)+V(x,y)\bigr\}.
\end{equation}
In particular, the function $W$ defined in \eqref{eqn:Wextension} is invariant
under $S_t$ and satisfies
\begin{equation}
    W(y)=\min_{x_j\in\mathcal{A}}\bigl\{W(x_j)+h(x_j,y)\bigr\}.
\end{equation}
\end{cor}

\begin{proof}
Fix $t>0$ and $y\in\mathbb{R}^d$. Since $w$ is dominated by $L$, every
absolutely continuous curve $\eta:[0,t]\to\mathbb{R}^d$ with
$\eta(0)=x$ and $\eta(t)=y$ satisfies
\[
    w(y)-w(x)
    \leq \int_0^t L(\dot\eta(s),\eta(s))\,ds.
\]
Taking first the infimum over all such curves and then the infimum over
$x\in\mathbb{R}^d$ gives
\[
    w(y)\leq \inf_{x\in\mathbb{R}^d}
    \bigl\{w(x)+h_t(x,y)\bigr\}=(S_t w)(y).
\]

For the reverse inequality, let
$\gamma:(-\infty,0]\to\mathbb{R}^d$ be a backward calibrated curve for
$w$ with $\gamma(0)=y$, and set $x_t:=\gamma(-t)$. Calibration on
$[-t,0]$ yields
\[
    w(y)-w(x_t)
    =\int_{-t}^0 L(\dot\gamma(s),\gamma(s))\,ds.
\]
The restriction of $\gamma$ to $[-t,0]$ is an admissible competitor in
the definition of $h_t(x_t,y)$. Therefore,
\[
\begin{aligned}
    (S_t w)(y)
    &\leq w(x_t)+h_t(x_t,y)\\
    &\leq w(x_t)+\int_{-t}^0
    L(\dot\gamma(s),\gamma(s))\,ds
    =w(y).
\end{aligned}
\]
Combining the two inequalities proves that $S_t w=w$ for every $t>0$.

It remains to derive the free-time identity. The domination argument above
shows that, for every $x\in\mathbb{R}^d$ and $T>0$,
\[
    w(y)\leq w(x)+h_T(x,y).
\]
Taking the infimum over $T>0$ and then over $x$ gives
\[
    w(y)\leq\inf_{x\in\mathbb{R}^d}
    \bigl\{w(x)+V(x,y)\bigr\}.
\]
On the other hand, $V(x,y)\leq h_t(x,y)$ for every $t>0$, so invariance
implies
\[
    \inf_{x\in\mathbb{R}^d}\bigl\{w(x)+V(x,y)\bigr\}
    \leq\inf_{x\in\mathbb{R}^d}\bigl\{w(x)+h_t(x,y)\bigr\}
    =w(y).
\]
This proves the free-time identity. Finally, $W$ is a weak KAM solution by
Corollary~\ref{cor:weakKAM}, while its representation over $\mathcal A$ is
precisely \eqref{eqn:Wextension}.
\end{proof}
For the extension $W$, coercivity allows us to apply Proposition~\ref{AB} in both directions and obtain uniqueness with prescribed values on $\mathcal A$.

\begin{cor}
\label{cor:unique-coercive-weak-kam}
Under the assumptions of Proposition~\ref{prop:calibrated}, the function $W$ constructed there is the unique coercive weak KAM solution of \eqref{eq:HJE0} satisfying $W(y)=W_y$ for $y\in\mathcal A$.
\end{cor}

\begin{proof}
By Proposition~\ref{prop:calibrated}(i) and Corollary~\ref{cor:weakKAM}, $W$ is a coercive weak KAM solution with the prescribed boundary values. Let $w$ be another coercive weak KAM solution with $w(y)=W_y$ on $\mathcal A$. Since $w=W$ on $\mathcal A$ and both functions are coercive, Proposition~\ref{AB}, applied in both directions, gives $w\le W$ and $W\le w$ on $\mathbb R^d$. Hence $w=W$.
\end{proof}

\bigskip

\section{Determine boundary values on Aubry set via LDP}\label{sec:ldp}
In Section~\ref{sec:weakKAM}, we showed that compatible values on the Aubry set uniquely determine a coercive weak KAM solution through the Peierls barrier. We now use the large deviation principle to select these values up to an additive constant, following the framework of \cite{Freidlin_Wentzell_2012}.

To introduce the prescribed boundary values on the Aubry set based on the large deviation principle, we introduce the definition of $\{x_i\}$-graph.

\begin{defn}
    Let $S$ be a finite set, and $E \subset S$. A graph consisting of arrows $m \to n \,(m \in S\setminus E,\ n \in S,\ n \neq m)$ is called an $E$-graph if it satisfies the following conditions:

    \noindent
    (i) every point $m \in S \setminus E$ is the initial point of exactly one arrow;

    \noindent
    (ii) there are no directed cycles in the graph.

    Let $G(E)$ denote the set of $E$-graphs, and we use $g$ to denote a generic $E$-graph.
\end{defn}

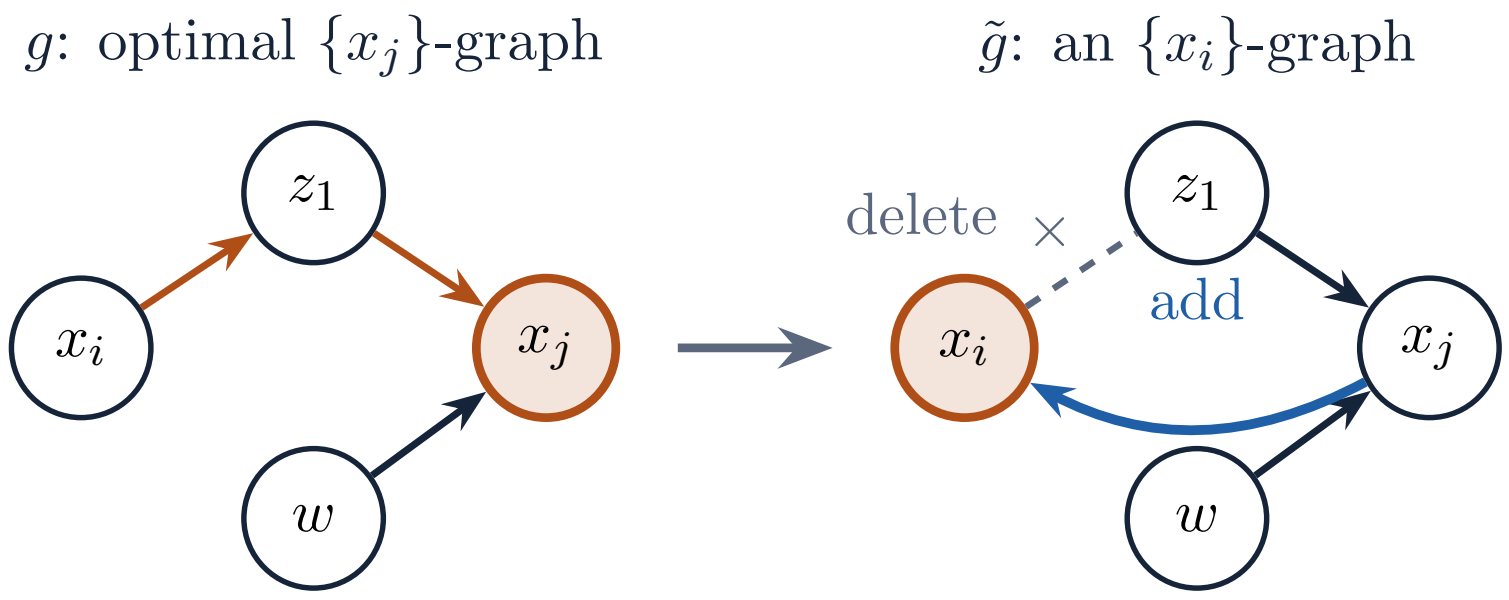
\begin{figure}[htbp]
\centering
\resizebox{0.6\linewidth}{!}{\begin{tikzpicture}[x=1cm,y=1cm,
  v/.style={circle,draw=ink,fill=white,line width=1pt,minimum size=0.9cm,inner sep=0pt},
  e/.style={-{Stealth[length=8pt]},line width=1.3pt,ink}]
\node[ink] at (1.5,1.95) {$g$: optimal $\{x_j\}$-graph};
\node[v] (i) at (0,0) {$x_i$};
\node[v] (z) at (1.5,1.0) {$z_1$};
\node[v,draw=org,fill=org!15,line width=1.6pt] (j) at (3.0,0) {$x_j$};
\node[v] (w) at (1.5,-1.1) {$w$};
\draw[e,org] (i) -- (z);
\draw[e,org] (z) -- (j);
\draw[e] (w) -- (j);
\draw[-{Stealth[length=10pt]},line width=1.6pt,mut] (3.85,0) -- (4.85,0);
\begin{scope}[xshift=5.7cm]
\node[ink] at (1.5,1.95) {$\tilde g$: an $\{x_i\}$-graph};
\node[v,draw=org,fill=org!15,line width=1.6pt] (i2) at (0,0) {$x_i$};
\node[v] (z2) at (1.5,1.0) {$z_1$};
\node[v] (j2) at (3.0,0) {$x_j$};
\node[v] (w2) at (1.5,-1.1) {$w$};
\draw[mut,dashed,line width=1.2pt] (i2) -- (z2);
\node[mut] at (0.55,0.75) {$\times$};
\draw[e] (z2) -- (j2);
\draw[e] (w2) -- (j2);
\draw[e,acc,line width=1.8pt] (j2) to[bend left=28] (i2);
\node[acc] at (1.5,0.32) {add};
\node[mut,anchor=south east] at (0.35,0.6) {delete};
\end{scope}
\end{tikzpicture}}
\caption{Left: an optimal $\{x_j\}$-graph $g$. There is only one arrow out of every other vertex, all leading to the root $\{x_i\}$. Right: the $\{x_i\}$-graph $\tilde g$ obtained from $g$ by deleting the arrow $x_i \to z_1$ and adding the arrow $x_j \to x_i$.}
\label{fig:egraph}
\end{figure}

Under (i), finiteness of $S$ makes (ii) equivalent to the following reachability condition (see \cite[p.~158]{Freidlin_Wentzell_2012}):
(ii') for any point $m \in S \setminus E$, there exists a sequence of arrows leading from it to some point $n \in E$.
Figure~\ref{fig:egraph} illustrates an $\{x_j\}$-graph and an $\{x_i\}$-graph.

In the following proposition, we first present several important properties of critical points. These properties imply that the Freidlin--Wentzell classes are singleton sets consisting of critical points of $\vb$ (see \cite[p.~150]{Freidlin_Wentzell_2012}). We then determine the boundary values on the Aubry set by taking the infimum over all directed graphs rooted at a given point, where the total cost is obtained by summing the costs associated with the edges of the graph (see \cite[Chapter 6, Lemmas 4.1 and 4.3(c)]{Freidlin_Wentzell_2012}).
\begin{prop}
Assume (A1) and (A3). Then the following statements hold.

\noindent
(i) Every $\omega$-limit set of the dynamical system $\dot{x}_t=\vb(x_t)$ is contained in $\mathcal{A}$.

\noindent
(ii) If $x,y\in\mathcal{A}$ and $x\neq y$, then either $V(x,y)\neq 0$ or $V(y,x)\neq 0$.

\noindent
(iii) Write $\mathcal A=\{x_1,\ldots,x_N\}$ and define $W_0:\mathcal A\to\mathbb R$ by
\begin{equation}\label{eqn:disbdy}
W_0(x_i):=\min_{g\in G\{x_i\}}\sum_{(x_m\to x_n)\in g}V(x_m,x_n),
\end{equation}
where $G\{x_i\}=G(\{x_i\})$ is the family of $\{x_i\}$-graphs on $S=\mathcal A$, and each edge $x_m\to x_n$ has cost $V(x_m,x_n)$. Then $W_0$ satisfies the discrete weak KAM relation \eqref{eqn:dwkam}.
\end{prop}

\begin{proof}
We first prove part (i). Let $\eta(t)$ solve $\dot{\eta}=\vb(\eta)$.
Then
\begin{equation}\label{eqn:Uly}
\frac{d}{dt}U(\eta(t))=\nabla U(\eta(t)) \cdot \vb(\eta(t))=-\left|\nabla U(\eta(t))\right|^2,
\end{equation}
which implies
\[
U(\eta(t))\leq U(\eta(0)).
\]
Since $U(x) \to \infty$ as $|x| \to \infty$, $\eta([0, \infty)) \subset K:=\left\{y \in \mathbb{R}^d : U(y) \leq U(\eta(0))\right\}$ where $K$ is compact. Thus, $\omega(\eta)$ is nonempty and compact. By LaSalle's invariance principle, $\omega(\eta) \subset \mathcal{A}$. Since $\omega(\eta)$ is connected and $\mathcal{A}$ consists of isolated finitely many critical points of $U$, $\omega(\eta)=\{x_i\}$ for some $x_i \in \mathcal{A}$.

Now for part (ii), suppose $x_i, x_j \in \mathcal{A}$ and $x_i \neq x_j$. By Proposition~\ref{prop:minimizer-regularity}, we know
\[
V(x,y) \geq U(y)-U(x), \qquad x, y\in \mathbb{R}^d.
\]
Thus if $U(x_i) \neq U(x_j)$, then $V(x_i,x_j)>0$ or $V(x_j,x_i)>0$. It remains to consider
the case
\[
U(x_i)=U(x_j)=:c, \qquad x_i\neq x_j.
\]
We claim that in this case both $V(x_i,x_j)$ and $V(x_j,x_i)$ are
strictly positive. We prove the first assertion; the other follows by
interchanging $x_i$ and $x_j$.

We proceed by contradiction and suppose that $V(x_i,x_j)=0$. Then there exist
$T_n>0$ and paths
$\gamma_n\in\operatorname{AC}([0,T_n];\mathbb R^d)$ joining $x_i$ to
$x_j$ such that
\[
A_n:=A_{T_n}[\gamma_n]
=\frac14\int_0^{T_n}
\left|\dot\gamma_n(t)-\vb(\gamma_n(t))\right|^2\,dt
\longrightarrow0.
\]
By Proposition~\ref{prop:minimizer-regularity}(a), for every
$t\in[0,T_n]$,
\[
U(\gamma_n(t))-c
\leq A_t[\gamma_n]
\leq A_n.
\]
Thus, for all sufficiently large $n$, the image of $\gamma_n$ is contained
in the compact set
\[
K:=\{x\in\mathbb R^d:U(x)\leq c+1\}.
\]
Since $\nabla U\cdot\vb=-|\nabla U|^2$, we have
\[
\begin{aligned}
\int_0^{T_n}|\nabla U(\gamma_n(t))|^2\,dt
&=\int_0^{T_n}\nabla U(\gamma_n(t))\cdot \left(\dot\gamma_n(t)-\vb(\gamma_n(t))\right)\,dt\\
&\leq
\left(\int_0^{T_n}|\nabla U(\gamma_n(t))|^2\,dt\right)^{1/2}
\left(\int_0^{T_n} \left|\dot\gamma_n(t)-\vb(\gamma_n(t))\right|^2\,dt\right)^{1/2}.
\end{aligned}
\]
Consequently,
\begin{equation}\label{eq:small-action-gradient}
\int_0^{T_n}|\nabla U(\gamma_n(t))|^2\,dt
\leq \int_0^{T_n}\left|\dot\gamma_n(t)-\vb(\gamma_n(t))\right|^2\,dt=4A_n.
\end{equation}

Write $\mathcal A=\{x_1,\ldots,x_N\}$. Since these points are distinct
and finite, choose $r>0$ so that the closed balls
$\overline{B_r(x_k)}$ are pairwise disjoint, and set
\[
d_0:=\min_{k\neq l}
\operatorname{dist}\bigl(\overline{B_r(x_k)},
\overline{B_r(x_l)}\bigr)>0.
\]
Because the critical set of $U$ is precisely $\mathcal A$, compactness
gives constants
\[
\kappa:=\min_{K\setminus\bigcup_{k=1}^N B_r(x_k)}|\nabla U|>0,
\qquad
M:=\max_K|\vb|<\infty.
\]
Every path from $x_i$ to $x_j$ must contain a subpath, defined on some
interval $[s_n,t_n]$, which joins the boundary of one of these balls to
the boundary of a different ball and otherwise stays outside their union.
In particular,
\[
|\gamma_n(t_n)-\gamma_n(s_n)|\geq d_0.
\]
Writing $\tau_n:=t_n-s_n$, \eqref{eq:small-action-gradient} gives
\[
\kappa^2\tau_n
\leq\int_{s_n}^{t_n}|\nabla U(\gamma_n(t))|^2\,dt
\leq4A_n,
\qquad
\tau_n\leq\frac{4A_n}{\kappa^2}.
\]
On the other hand, since $\dot\gamma_n=\vb(\gamma_n)+\left(\dot\gamma_n(t)-\vb(\gamma_n(t))\right)$,
\[
\begin{aligned}
d_0
&\leq\int_{s_n}^{t_n}|\dot\gamma_n(t)|\,dt\leq M\tau_n+
\tau_n^{1/2}
\left(\int_{s_n}^{t_n}\left|\dot\gamma_n(t)-\vb(\gamma_n(t))\right|^2\,dt\right)^{1/2}\\
&\leq
\left(\frac{4M}{\kappa^2}+\frac{4}{\kappa}\right)A_n,
\end{aligned}
\]
which contradicts the fact that $A_n\to0$. Hence $V(x_i,x_j)>0$. Reversing the
roles of the endpoints gives $V(x_j,x_i)>0$.

Finally, we prove part (iii). If $x_i=x_j$, the required inequality is immediate since $V(x_i,x_i)=0$. Assume $x_i\neq x_j$ and choose a minimizing $\{x_j\}$-graph $g$, which exists because $G\{x_j\}$ is finite and nonempty. Then
\[
W_0(x_j)=\sum_{(x_m\to x_n)\in g} V(x_m,x_n).
\]
There exists a directed path going from $x_i$ to $x_j$, that is,
\[
x_i=z_0 \rightarrow z_1 \rightarrow \cdots \rightarrow z_r=x_j,
\]
where $\left\{z_0, \cdots, z_r\right\} \subset \mathcal{A}$.
Now add the arrow $x_j \rightarrow x_i$, which creates a cycle
\[
x_j \rightarrow x_i=z_0 \rightarrow z_1 \rightarrow \cdots \rightarrow z_r=x_j.
\]
Delete the original arrow leaving $x_i$, namely, $x_i \to z_1$.
The resulting graph $\tilde{g}$ is a $\{x_i\}$-graph; see illustration in Figure ~\ref{fig:egraph}. Therefore,
\[
W_0(x_i) \leq \sum_{(x_m\to x_n)\in\tilde g}V(x_m,x_n)
= W_0(x_j) + V(x_j,x_i)-V(x_i,z_1)
\leq W_0(x_j)+V(x_j,x_i),
\]
which proves that $W_0$ satisfies the discrete weak KAM problem \eqref{eqn:dwkam}.
\end{proof}

We illustrate the graph construction \eqref{eqn:disbdy} with a two-dimensional double-well potential. By computing the transition costs and minimizing rooted graphs explicitly, the example shows how the normalized values of $W_0$ recover the potential levels at the two minima and the saddle.

\begin{exm}[A two-dimensional double-well potential]
Consider the reversible diffusion on $\mathbb R^2$ with
\[
U(x_1,x_2)=\frac14(x_1^2-1)^2+\frac12x_2^2,
\qquad
\vb(x_1,x_2)=-\nabla U(x_1,x_2)
=\bigl(x_1-x_1^3,-x_2\bigr).
\]
Writing $e_1=(1,0)$, the only three critical points of $U$ are the two
minima $x_-=-e_1$, $x_+=e_1$, and the index-one saddle $x_0=0$.
Consequently,
\[
\mathcal A=\{x_-,x_0,x_+\}=\{-e_1,0,e_1\},
\]
with
\[
U(x_-)=U(x_+)=0,
\qquad
U(x_0)=\frac14.
\]
For this gradient drift, the cost of moving from either minimum to the saddle
is the corresponding increase of $U$, whereas the Ma\~n\'e cost from the
saddle into either well is zero. More precisely,
\[
\begin{array}{c|ccc}
V(x,y) & y=x_- & y=x_0 & y=x_+\\ \hline
x=x_- & 0 & \frac14 & \frac14\\
x=x_0 & 0 & 0 & 0\\
x=x_+ & \frac14 & \frac14 & 0
\end{array}
\]
To justify these values in the above cost table, Proposition~\ref{prop:minimizer-regularity}(a) first
gives
\[
V(x_\pm,x_0)\geq U(x_0)-U(x_\pm)=\frac14.
\]
The reverse inequality is obtained, in the variational sense, by perturbing
the initial and terminal points slightly and following the reversed gradient
flow $\dot\gamma=\nabla U$ along the $x_1$-axis. Along such a segment,
\[
\frac14\left|\dot\gamma+\nabla U(\gamma)\right|^2
=|\nabla U(\gamma)|^2
=\frac{d}{dt}U(\gamma),
\]
so its action is the increase of $U$. The short connecting pieces at the two
ends have action tending to zero as the perturbations tend to zero. Hence
$V(x_\pm,x_0)=1/4$.

Similarly, $V(x_0,x_\pm)=0$. This statement is variational: a deterministic
trajectory starting exactly at $x_0$ remains there. Instead, connect $x_0$
to $\pm\delta e_1$ with action $O(\delta^2)$, follow the forward gradient
flow $\dot\gamma=-\nabla U$ toward $x_\pm$ at zero cost, and use a final
short connector whose action tends to zero. Letting $\delta\downarrow0$ and
the terminal gap tend to zero proves the claim.

Finally, every path from $x_-$ to $x_+$ must cross the line $\{x_1=0\}$,
on which $U\geq U(x_0)=1/4$. Applying
Proposition~\ref{prop:minimizer-regularity}(a) up to the first crossing yields
\[
V(x_-,x_+)\geq U(x_0)-U(x_-)=\frac14.
\]
Concatenating the preceding approximate uphill path from $x_-$ to $x_0$ with
the approximate zero-cost path from $x_0$ to $x_+$ gives the reverse
inequality. Thus $V(x_-,x_+)=1/4$; by symmetry,
$V(x_+,x_-)=1/4$. This proves the above cost table. The quadratic term in
$x_2$ makes motion away from the $x_1$-axis unnecessary and introduces no
additional critical points.

We now compute \eqref{eqn:disbdy}. For a graph rooted at $x_-$, both of the
other vertices must eventually lead to $x_-$. All possible $\left\{x_{-}\right\}$-graphs are
\[
x_+\longrightarrow x_0\longrightarrow x_- ,
\]
\[
x_0\longrightarrow x_+\longrightarrow x_- ,
\]
\[
x_0\longrightarrow x_- ,
\qquad
x_+\longrightarrow x_- .
\]
A minimizing choice is
\[
x_0\longrightarrow x_- ,
\qquad
x_+\longrightarrow x_- ,
\]
whose total cost is $0+1/4=1/4$. In fact, the other two graphs also have
cost $1/4$, so all three are minimizing. By symmetry, every graph rooted at
$x_+$ also has cost $1/4$. If the root is $x_0$, the three possible graphs
are
\[
x_-\longrightarrow x_0,
\qquad
x_+\longrightarrow x_0,
\]
\[
x_-\longrightarrow x_+\longrightarrow x_0,
\]
and
\[
x_+\longrightarrow x_-\longrightarrow x_0.
\]
Each contains two edges of cost $1/4$, and therefore has total cost $1/2$.
Thus
\[
W_0(x_-)=\frac14,
\qquad
W_0(x_0)=\frac12,
\qquad
W_0(x_+)=\frac14.
\]
After the normalization used in Proposition~\ref{prop:ldp},
\[
W_0(x)-\min_{z\in\mathcal A}W_0(z)
=
\begin{cases}
0, & x=x_-\text{ or }x=x_+,\\[2mm]
\dfrac14, & x=x_0,
\end{cases}
\qquad x\in\mathcal A.
\]
Since the invariant density in this reversible example is proportional to
$e^{-U/\varepsilon}$, these normalized graph values agree with
$U-\min_{\mathbb R^2}U$ on the Aubry set. Thus the rooted-graph formula
assigns the correct relative levels to the two wells and the saddle.
\end{exm}

The following proposition is a direct application of \cite[Chapter 6, Theorem 4.3]{Freidlin_Wentzell_2012}.

\begin{prop}\label{prop:ldp}
Assume (A1)--(A3). Extend $W_0$ from \eqref{eqn:disbdy} to $W_0:\mathbb{R}^d \to \mathbb{R}$ by
\[
W_0(x)=\min_{x_i \in\mathcal{A}} \left[W_0(x_i)+V(x_i, x)\right].
\]
Let $\gamma >0$. Then, for every sufficiently small neighborhood $E_\rho(x)$ of $x$, there exists $\varepsilon_0>0$ such that for all $0<\varepsilon \leq \varepsilon_0$,
\[
\exp \left\{-\varepsilon^{-1} \left(W_0(x)-\min_{i\in \mathcal{A}} W_0(x_i)+\gamma\right)\right\} \leq \pi_\varepsilon(E_\rho(x)) \leq \exp \left\{-\varepsilon^{-1} \left(W_0(x)-\min_{i\in \mathcal{A}} W_0(x_i)-\gamma\right)\right\}
\]
\end{prop}

Proposition~\ref{prop:ldp} estimates the invariant measure of small neighborhoods. The next result passes from these probability estimates to locally uniform logarithmic asymptotics of the density, using a uniform local gradient bound on $-\varepsilon\log m_\varepsilon$ to control its oscillation.

\begin{prop}[Density-level logarithmic asymptotics]
\label{prop:density-ldp}
Under the same assumptions as Proposition~\ref{prop:ldp}, let
\begin{equation}\label{eq:log-density-and-rate}
I(x):=W_0(x)-\min_{x_i\in\mathcal A}W_0(x_i),
\qquad
u_\varepsilon(x):=-\varepsilon\log m_\varepsilon(x), \qquad x\in \mathbb{R}^d.
\end{equation}

Suppose in addition that for every compact
set $K\subset\mathbb R^d$ there are $C_K>0$ and $\varepsilon_K>0$ such that
\begin{equation}\label{eq:local-gradient-density}
\|\nabla u_\varepsilon\|_{L^\infty(K)}\leq C_K,
\qquad 0<\varepsilon\leq\varepsilon_K.
\end{equation}
Then
\[
u_\varepsilon\longrightarrow I
\qquad\text{locally uniformly in }\mathbb R^d.
\]
In particular, $u_\varepsilon(x)\to I(x)$ for every $x$.
\end{prop}

\begin{proof}
By Remark~\ref{rem:wlip}, $I$ is locally Lipschitz and hence continuous. Fix $x\in\mathbb R^d$, and choose a compact neighborhood $K$ of $x$. For
$\rho>0$ sufficiently small, $\overline{B_\rho(x)}\subset K$. Set
\[
q_\varepsilon(\rho,x)
:=-\varepsilon\log\pi_\varepsilon(B_\rho(x)), \qquad
a_\varepsilon(\rho,x):=\inf_{B_\rho(x)}u_\varepsilon,
\qquad
b_\varepsilon(\rho,x):=\sup_{B_\rho(x)}u_\varepsilon.
\]
Since $m_\varepsilon=e^{-u_\varepsilon/\varepsilon}$, we have
\[
|B_\rho|e^{-b_\varepsilon(\rho,x)/\varepsilon}
\leq \pi_\varepsilon(B_\rho(x))
\leq |B_\rho|e^{-a_\varepsilon(\rho,x)/\varepsilon}.
\]
Taking $-\varepsilon\log$ yields
\[
a_\varepsilon(\rho,x)-\varepsilon\log|B_\rho|
\leq q_\varepsilon(\rho,x)
\leq b_\varepsilon(\rho,x)-\varepsilon\log|B_\rho|.
\]
Equivalently, the preceding inequality can be written as
\[
a_\varepsilon(\rho,x)
\leq q_\varepsilon(\rho,x)+\varepsilon\log|B_\rho|
\leq b_\varepsilon(\rho,x).
\]
Moreover, since $x\in B_\rho(x)$, the definition of
$a_\varepsilon(\rho,x)$ and $b_\varepsilon(\rho,x)$ gives
\[
a_\varepsilon(\rho,x)
\leq u_\varepsilon(x)
\leq b_\varepsilon(\rho,x).
\]
Thus both $u_\varepsilon(x)$ and
$q_\varepsilon(\rho,x)+\varepsilon\log|B_\rho|$ belong to the same interval
$[a_\varepsilon(\rho,x),b_\varepsilon(\rho,x)]$. Their distance is therefore
at most the length of this interval:
\begin{equation}\label{eq:density-ball-comparison}
\left|
u_\varepsilon(x)-q_\varepsilon(\rho,x)
-\varepsilon\log|B_\rho|
\right|
\leq \operatorname{osc}_{B_\rho(x)}u_\varepsilon
\leq 2C_K\rho,
\end{equation}
where the last inequality follows from
\eqref{eq:local-gradient-density}: for any $z_1,z_2\in B_\rho(x)$,
\[
|u_\varepsilon(z_1)-u_\varepsilon(z_2)|
\leq C_K|z_1-z_2|
\leq 2C_K\rho.
\]

We next derive the required fixed-ball bounds directly from the local
estimates in Proposition~\ref{prop:ldp}. Let $\eta>0$. Choose
$y_\eta\in B_\rho(x)$ such that
\[
I(y_\eta)\leq \inf_{B_\rho(x)}I+\eta.
\]
Since $B_\rho(x)$ is open, we may choose a sufficiently small neighborhood
$E_\eta(y_\eta)$ of $y_\eta$ such that
$E_\eta(y_\eta)\subset B_\rho(x)$ and the local lower estimate in
Proposition~\ref{prop:ldp}, with error $\eta$, applies. Thus, for all
sufficiently small $\varepsilon$,
\[
\pi_\varepsilon(B_\rho(x))
\geq \pi_\varepsilon\bigl(E_\eta(y_\eta)\bigr)
\geq \exp\left\{-\frac{I(y_\eta)+\eta}{\varepsilon}\right\}.
\]
It follows that
\[
\limsup_{\varepsilon\downarrow0}q_\varepsilon(\rho,x)
\leq I(y_\eta)+\eta
\leq \inf_{B_\rho(x)}I+2\eta.
\]
Letting $\eta\downarrow0$ gives
\begin{equation}\label{eq:fixed-ball-upper-q}
\limsup_{\varepsilon\downarrow0}q_\varepsilon(\rho,x)
\leq \inf_{B_\rho(x)}I.
\end{equation}

For the converse bound, set $F:=\overline{B_\rho(x)}$. For each $y\in F$,
choose a sufficiently small open neighborhood $E(y)$ for which the local
upper estimate, again with error $\eta$, holds. The sets $E(y)$ cover the
compact set $F$, so there are $y_1,\ldots,y_N\in F$ such that
\[
F\subset \bigcup_{j=1}^N E(y_j).
\]
For sufficiently small $\varepsilon$,
\[
\begin{aligned}
\pi_\varepsilon(F)
\leq \sum_{j=1}^N\pi_\varepsilon(E(y_j)) \leq \sum_{j=1}^N
\exp\left\{-\frac{I(y_j)-\eta}{\varepsilon}\right\} \leq N\exp\left\{-\frac{\inf_F I-\eta}{\varepsilon}\right\}.
\end{aligned}
\]
Since $B_\rho(x)\subset F$, taking $-\varepsilon\log$ gives
\[
q_\varepsilon(\rho,x)
\geq \inf_F I-\eta-\varepsilon\log N.
\]
First letting $\varepsilon\downarrow0$ and then $\eta\downarrow0$ yields
\begin{equation}\label{eq:fixed-ball-lower-q}
\liminf_{\varepsilon\downarrow0}q_\varepsilon(\rho,x)
\geq \inf_{\overline{B_\rho(x)}}I.
\end{equation}
Combining \eqref{eq:fixed-ball-upper-q} and
\eqref{eq:fixed-ball-lower-q}, we have
\[
\inf_{\overline{B_\rho(x)}}I
\leq \liminf_{\varepsilon\downarrow0}q_\varepsilon(\rho,x)
\leq \limsup_{\varepsilon\downarrow0}q_\varepsilon(\rho,x)
\leq \inf_{B_\rho(x)}I.
\]
Since $I$ is continuous and $B_\rho(x)$ is dense in
$\overline{B_\rho(x)}$,
$\inf_{\overline{B_\rho(x)}}I=\inf_{B_\rho(x)}I$.
Consequently, for every fixed $\rho>0$,
\begin{equation}\label{eq:fixed-ball-rate}
\lim_{\varepsilon\downarrow0}q_\varepsilon(\rho,x)
=\inf_{B_\rho(x)}I.
\end{equation}
Let
\[
\omega_I(\rho,x)
:=\sup_{y\in\overline{B_\rho(x)}}|I(y)-I(x)|.
\]
Since $x\in B_\rho(x)$,
\[
\left|\inf_{B_\rho(x)}I-I(x)\right|
\leq \omega_I(\rho,x).
\]
For fixed $\rho>0$, $\varepsilon\log|B_\rho|\to0$. Combining
\eqref{eq:density-ball-comparison} and \eqref{eq:fixed-ball-rate}, we obtain
\[
\begin{aligned}
\limsup_{\varepsilon\downarrow0}|u_\varepsilon(x)-I(x)|
\leq 2C_K\rho+
\left|\inf_{B_\rho(x)}I-I(x)\right|\leq 2C_K\rho+\omega_I(\rho,x).
\end{aligned}
\]
By continuity of $I$, $\omega_I(\rho,x)\to0$ as $\rho\downarrow0$.
Letting $\rho\downarrow0$ proves
$u_\varepsilon(x)\to I(x)$ for every $x$.

It remains to show that the convergence is locally uniform. Suppose instead that it is not uniform on some compact set $K$. Then there exist
$\eta>0$, a sequence $\varepsilon_n\downarrow0$, and $x_n\in K$ such that
\[
|u_{\varepsilon_n}(x_n)-I(x_n)|\geq\eta.
\]
After passing to a subsequence, $x_n\to x\in K$. Choose a compact
neighborhood $K'$ of $K$. For all sufficiently large $n$, the segment joining $x$ to $x_n$ lies in $K'$. The gradient estimate on $K'$, the pointwise
convergence at $x$, and the continuity of $I$ give
\begin{align*}
|u_{\varepsilon_n}(x_n)-I(x_n)|
\leq{} |u_{\varepsilon_n}(x_n)-u_{\varepsilon_n}(x)|
+|u_{\varepsilon_n}(x)-I(x)|+|I(x)-I(x_n)|\longrightarrow 0,
\end{align*}
which is a contradiction. Hence the convergence is locally uniform.
\end{proof}

\begin{rem}
    In the case $\vb=-\nabla U$, this conclusion follows directly
from Laplace principle:
\[
m_\varepsilon(x)=Z_\varepsilon^{-1}e^{-U(x)/\varepsilon},
\qquad
-\varepsilon\log m_\varepsilon(x)\longrightarrow U(x)-\min U.
\]
The same explicit invariant density holds for
$\vb=-\nabla U+\vl$ if, in addition to
$\vl\cdot\nabla U=0$, one assumes $\nabla\cdot\vl=0$.
\end{rem}

Proposition~\ref{prop:density-ldp} reduces the convergence of the logarithmic invariant densities to establishing the local gradient bound \eqref{eq:local-gradient-density}. The following corollary obtains this bound from a local Bernstein estimate under (A4).

\begin{cor}[Density convergence under strengthened regularity]
\label{cor:density-convergence-strong-regularity}
Assume (A1), (A2), and (A4).
Then, with $u_\varepsilon$ and $I$ defined in \eqref{eq:log-density-and-rate}, we have
\[
u_\varepsilon\longrightarrow I
\qquad\text{locally uniformly in }\mathbb R^d
\quad\text{as }\varepsilon\downarrow0.
\]
\end{cor}

\begin{proof}
By (A4), we have $\vb=-\nabla U+\vl\in C^2(\mathbb R^d;\mathbb R^d)$. Set $\varepsilon_0:=\min\{1,1/(2dN)\}$, where $N$ is as in (A2). For $0<\varepsilon\leq\varepsilon_0$, we have $\varepsilon<1/(dN)$, so (A2) ensures existence of a strictly positive invariant density $m_\varepsilon$, which satisfies
\[
\Delta m_\varepsilon
-\nabla\cdot\left(\frac{\vb}{\varepsilon}m_\varepsilon\right)=0
\qquad\text{in }\mathbb R^d.
\]
Apply \cite[Theorem~5.2]{MetafunePallaraRhandi2005} to the rescaled generator
\[
\widetilde L_\varepsilon
=\Delta+\frac{\vb}{\varepsilon}\cdot\nabla,
\]
which has the same invariant density as
$L_\varepsilon=\varepsilon\Delta+\vb\cdot\nabla$. If
$K\subset\mathbb R^d$ is compact and
\[
K_1:=\{y\in\mathbb R^d:\operatorname{dist}(y,K)\leq1\},
\]
then
\[
\begin{aligned}
|\nabla\log m_\varepsilon(x)|
&\leq C\sup_{\substack{y\in\mathbb R^d\\|y-x|\leq1}}
\left(
1+\frac{|\vb(y)|+|D\vb(y)|+|D^2\vb(y)|}{\varepsilon}
\right)\\
&\leq \frac{C_K}{\varepsilon},
\qquad x\in K,\quad 0<\varepsilon\leq\varepsilon_0,
\end{aligned}
\]
since $\vb\in C^2$ and hence the relevant derivatives are bounded on
$K_1$. Consequently,
\[
\|\nabla u_\varepsilon\|_{L^\infty(K)}
=\varepsilon\|\nabla\log m_\varepsilon\|_{L^\infty(K)}
\leq C_K.
\]
The conclusion now follows from Proposition~\ref{prop:density-ldp}.
\end{proof}

\bigskip

\section{Selection principle via large deviation principle}\label{sec:selection}

Weak KAM solutions to the stationary HJE are not unique: any boundary values $\{W_j\}$ satisfying the discrete weak KAM problem yield a distinct solution. A selection principle is needed to identify the physically meaningful one. We show that first taking the long time limit $t\to+\infty$ (ergodicity) and then the zero noise limit $\epss\to 0$ (large deviation principle for the invariant measure) selects a unique weak KAM solution, which coincides with the global energy landscape of the stochastic process.

\subsection{Exchange of limits and the selected solution}

We compare the long-time and zero-noise limits of the logarithmic density.
For each sufficiently small fixed noise level, the probability density
converges to the invariant density as $t\to+\infty$, by
Proposition~\ref{prop:fixed-epsilon-density-convergence}. Under the additional
assumption (A4), Proposition~\ref{prop:finite-time-wkb} gives the zero-noise
limit $S_tI=I$ at every fixed time for the initial density proportional to
$e^{-I/\varepsilon}$. Thus the two iterated limits agree, as shown in
Theorem~\ref{thm:exchange}.

\begin{prop}[Fixed-noise convergence to equilibrium]
\label{prop:fixed-epsilon-density-convergence}
Assume (A1)--(A3). Let $N$ be as in (A2), fix $0<\varepsilon<1/(dN)$, and let
$\rho_\varepsilon$ be the probability solution of the Fokker--Planck equation~\eqref{FP}
with an initial probability density
$\rho_\varepsilon(\cdot,0)\in L^1(\mathbb R^d)$. Then
\[
\rho_\varepsilon(\cdot,t)\longrightarrow m_\varepsilon
\qquad\text{in }L^1(\mathbb R^d)
\quad\text{as }t\to+\infty.
\]
Moreover, for every $x\in\mathbb R^d$,
\[
\rho_\varepsilon(x,t)\longrightarrow m_\varepsilon(x)
\qquad\text{as }t\to+\infty.
\]
\end{prop}

\begin{proof}
\emph{Step 1: existence, invariance, and uniqueness of the equilibrium.}
Let $X_t$ solve with initial condition $X_0=x\in\mathbb R^d$
\[
 dX_t=\vb(X_t)\,dt+\sqrt{2\varepsilon}\,dB_t,
 \qquad
 \mathcal L_\varepsilon=\vb\cdot\nabla+\varepsilon\Delta,
\]
and put $V(x)=|x|^2$.  Assumption (A2) and
$\varepsilon<1/(dN)$ give
\begin{equation}
\label{eq:fixed-epsilon-lyapunov}
 \mathcal L_\varepsilon V(x)
 =2\vb(x)\cdot x+2d\varepsilon
 <-\frac2N+2d\varepsilon<0,
 \qquad |x|>N.
\end{equation}
Hence \cite[Corollary~2.4.2]{Bogachev2015} yields a stationary probability density
$m_\varepsilon$.  Local elliptic regularity makes $m_\varepsilon$ continuous,
and the Harnack inequality for the stationary Fokker--Planck
equation shows that $m_\varepsilon$ is strictly positive on the connected
domain $\mathbb R^d$.  Write
$\pi_\varepsilon(dx)=m_\varepsilon(x)\,dx$.

By (A1) and (A3), $\vb\in C^1(\mathbb R^d;\mathbb R^d)$ and hence is locally Lipschitz. Standard SDE theory therefore
gives a pathwise unique maximal strong solution, defined up to its explosion
time.  It remains to show that this explosion time is almost surely infinite.
Since
$\mathcal L_\varepsilon V$ is bounded above on $\mathbb R^d$, set
$C_\varepsilon:=\max\{0,\sup_{\mathbb R^d}\mathcal L_\varepsilon V\}<\infty$.
For $R>|x|$, Dynkin's
formula, stopped at
$\tau_R:=\inf\{s\geq0:|X_s|\geq R\}$, gives
\[
 R^2\mathbb P_x(\tau_R\leq t)
 \leq \mathbb E_xV(X_{t\wedge\tau_R})
 \leq V(x)+C_\varepsilon t.
\]
Letting $R\to\infty$ proves nonexplosion.  Local pathwise uniqueness together
with nonexplosion therefore gives a unique global solution of the SDE.  Its
transition kernel is
\[
 P_t(x,A):=\mathbb P_x(X_t\in A)
\]
and the Markov property gives the semigroup law.  Because the solution exists
in $\mathbb R^d$ for every finite time, the semigroup is conservative:
$P_t(x,\mathbb R^d)=1$.

The Lyapunov bound
$\mathcal L_\varepsilon(1+V)\leq(1+C_\varepsilon)(1+V)$ and
\cite[Section~9.4]{Bogachev2015} give uniqueness of probability solutions
of \eqref{FP} with prescribed initial data. The constant family $\pi_\varepsilon$
and the family $\pi_\varepsilon P_t$ therefore coincide, so
\[
 \pi_\varepsilon P_t=\pi_\varepsilon,
 \qquad t\geq0.
\]
Finally, $V_1=1+|x|^2$ satisfies
$\mathcal L_\varepsilon V_1\leq CV_1$ on $\mathbb R^d$; hence condition
(iii) of \cite[Theorem~4.1.6]{Bogachev2015} shows that
$\pi_\varepsilon$ is the unique stationary probability measure.

\medskip
\noindent
\emph{Step 2: convergence in $L^1$.}
Since $\varepsilon I_d$ is uniformly elliptic and
$\vb\in C^1\subset L^p_{\mathrm{loc}}$ for every finite
$p$, the local hypotheses of standard parabolic theory apply.  In particular,
by \cite[Theorem~4.1]{BogachevKrylovRockner2001}, there exists a sub-probability kernel
$K_t(x,dy)=p_t(x,y)\,dy$ with $p_t$ jointly continuous, nonnegative, and
locally H\"older continuous for $t>0$.  The uniqueness of the nonexplosive
martingale problem established in Step~1 identifies this kernel with the SDE
transition kernel $\pi_\varepsilon$-almost everywhere in $x$.
For $\varphi\in C_c^\infty(\mathbb R^d)$, $P_t\varphi$ is continuous
by \cite[Theorems~3.4--3.5]{Khasminskii}, and $K_t\varphi$ is continuous
by \cite[Theorem~4.1]{BogachevKrylovRockner2001}.
Since $m_\varepsilon>0$, $\pi_\varepsilon$ has full support, so their
almost-everywhere equality extends to every $x$.
Testing against all such $\varphi$ then gives
\[
 K_t(x,dy)=P_t(x,dy)=p_t(x,y)\,dy,
 \qquad \int_{\mathbb R^d}p_t(x,y)\,dy=1.
\]
For fixed $x$, $p_t(x,y)$ is a nonnegative solution of the Fokker--Planck
equation~\eqref{FP} in the $(t,y)$ variables.  Since it has total mass one, it is not
identically zero; the parabolic strong maximum principle therefore gives
$p_t(x,y)>0$ for every $t>0$ and $x,y\in\mathbb R^d$.

We next verify the recurrence hypotheses needed for the continuous-time
ergodic theorem. Fix $h>0$ and consider the skeleton chain
$Y_n:=X_{nh}$, whose one-step transition kernel is $P_h$. The strict
positivity of the transition density immediately implies irreducibility:
for every Borel set $A\subset\mathbb R^d$ with positive Lebesgue measure,
\[
 P_h(x,A)=\int_A p_h(x,y)\,dy>0
 \qquad\text{for every }x\in\mathbb R^d.
\]
Thus the $h$-skeleton is Lebesgue-irreducible.

By Step~1, $\pi_\varepsilon P_h=\pi_\varepsilon$. Since
$m_\varepsilon>0$ everywhere, $\pi_\varepsilon$ is equivalent to Lebesgue
measure. Moreover, since $P_h(x,dy)=p_h(x,y)\,dy$, we have
\[
 P_h(x,\cdot)\ll\pi_\varepsilon
 \qquad\text{for every }x\in\mathbb R^d.
\]
Hence, if $G$ is a Borel set with $\pi_\varepsilon(G)=1$, then
$\pi_\varepsilon(G^c)=0$, and therefore
\[
 P_h(x,G^c)=0,
 \qquad\text{equivalently}\qquad
 \mathbb P_x(Y_1\in G)=1,
 \qquad x\in\mathbb R^d.
\]
Thus every set of full invariant measure is reached almost surely from every
initial state (in fact, after one skeleton step). Since the skeleton is
Lebesgue-irreducible and admits the invariant probability measure
$\pi_\varepsilon$, the Harris recurrence criterion
\cite[Theorem~6, equivalence of (i) and (v)]{RobertsRosenthal2006Harris}
shows that the skeleton chain is Harris recurrent.

We now transfer this property to the continuous-time process. In the
notation of \cite[Section~2.3]{MeynTweedie1993StabilityII}, a sampled kernel
has the form
\[
 K_a(x,A):=\int_0^\infty P_s(x,A)\,a(ds),
\]
where $a$ is a probability measure on $\mathbb R_+$. Choosing
$a=\delta_h$ gives $K_a=P_h$. Therefore
\cite[Theorem~2.2(i)]{MeynTweedie1993StabilityII} implies that the diffusion
itself is Harris recurrent. Since it admits the invariant probability
measure $\pi_\varepsilon$, it is positive Harris recurrent. Together with
the existence of the irreducible $h$-skeleton established above,
\cite[Theorem~6.1]{MeynTweedie1993StabilityII} gives
\begin{equation}
\label{eq:point-mass-tv-convergence}
 \|P_t(x,\cdot)-\pi_\varepsilon\|_{\mathrm{TV}}
 \longrightarrow0,
 \qquad x\in\mathbb R^d.
\end{equation}
Here and below we use the total variation distance convention
\[
 \|\mu-\nu\|_{\mathrm{TV}}
 :=\sup_{A\in\mathcal B(\mathbb R^d)}|\mu(A)-\nu(A)|
 =\frac12\sup_{|\varphi|\leq1}
 \left|\int\varphi\,d(\mu-\nu)\right|
\]
for probability measures $\mu$ and $\nu$, where the functions $\varphi$
are measurable. If $\mu$ and $\nu$ have densities $f$ and $g$, respectively, then
\[
 \|f-g\|_{L^1}=2\|\mu-\nu\|_{\mathrm{TV}}.
\]

Let $\mu_0^\varepsilon(dx)=\rho_\varepsilon(x,0)\,dx$. By the uniqueness
of probability solutions established in Step~1, the law
$\mu_0^\varepsilon P_t$ has density $\rho_\varepsilon(\cdot,t)$, and
by the definition of the push-forward law,
\[
 \mu_0^\varepsilon P_t-\pi_\varepsilon
 =\int_{\mathbb R^d}
 \big(P_t(y,\cdot)-\pi_\varepsilon\big)\,\mu_0^\varepsilon(dy).
\]
The triangle inequality for the total variation distance therefore gives
\[
 \|\mu_0^\varepsilon P_t-\pi_\varepsilon\|_{\mathrm{TV}}
 \leq\int_{\mathbb R^d}
 \|P_t(y,\cdot)-\pi_\varepsilon\|_{\mathrm{TV}}
 \,\mu_0^\varepsilon(dy).
\]
Here $\mu_0^\varepsilon$ is a probability measure because
$\rho_\varepsilon(\cdot,0)$ is an initial probability density.  Define
\[
 h_t(y):=\|P_t(y,\cdot)-\pi_\varepsilon\|_{\mathrm{TV}}.
\]
The function $h_t$ is measurable, $h_t(y)\to0$ for every $y$ by
\eqref{eq:point-mass-tv-convergence}, and $0\leq h_t(y)\leq1$.  The dominated convergence theorem then gives
\[
 \lim_{t\to\infty}\int_{\mathbb R^d}h_t(y)\,\mu_0^\varepsilon(dy)
 =\int_{\mathbb R^d}\lim_{t\to\infty}h_t(y)\,\mu_0^\varepsilon(dy)=0.
\]
Since $\mu_0^\varepsilon P_t$ and $\pi_\varepsilon$
have densities $\rho_\varepsilon(\cdot,t)$ and $m_\varepsilon$, respectively,
\[
 2\|\mu_0^\varepsilon P_t-\pi_\varepsilon\|_{\mathrm{TV}}
 =\sup_{|\varphi|\leq1}
 \left|\int_{\mathbb R^d}\varphi(x)
 \bigl(\rho_\varepsilon(x,t)-m_\varepsilon(x)\bigr)\,dx\right|
 =\|\rho_\varepsilon(\cdot,t)-m_\varepsilon\|_{L^1(\mathbb R^d)}.
\]
Consequently,
\[
 \|\rho_\varepsilon(\cdot,t)-m_\varepsilon\|_{L^1(\mathbb R^d)}
 \leq2\int_{\mathbb R^d}h_t(y)\,\mu_0^\varepsilon(dy)
 \longrightarrow0.
\]
\medskip

\noindent
\emph{Step 3: pointwise convergence.}
Set $f^\varepsilon:=\rho_\varepsilon-m_\varepsilon$. Since $m_\varepsilon$ is stationary, $f^\varepsilon$ satisfies
\[
\partial_t f^\varepsilon
=
\varepsilon\Delta f^\varepsilon
-\nabla\cdot\bigl(\vb f^\varepsilon\bigr)
\qquad\text{in }\mathbb R^d\times(0,\infty).
\]
Fix $x\in\mathbb R^d$, $r>0$, and $t\geq2$.  Throughout this step,
$\varepsilon>0$ is fixed.  For $1\leq a\leq2$, define the space--time cylinders
\[
 Q_a:=B_{ar}(x)\times(t-a/2,t+a/2),
 \qquad M(a):=\sup_{Q_a}|f^\varepsilon|,
 \qquad J_t:=\int_{Q_2}|f^\varepsilon(y,s)|\,dy\,ds.
\]
Note that $Q_a\subset Q_b$ whenever $a<b$.
Local continuity of $f^\varepsilon$ implies $M(3/2)<\infty$ for each
fixed $t$, since $\overline{Q_{3/2}}$ is a compact subset of
$\mathbb R^d\times(0,\infty)$.

For a point $(y,s)$ and
$0<h\leq 1$, write
\[
 D_h(y,s):=B_{h/2}(y)\times(s-h^2,s].
\]
This is the cylinder denoted by $Q(h)$ in
\cite[Theorem~B, p.~616]{Aronson1968Nonnegative}, with top center $(y,s)$.
We fix the upper scale bound in that theorem to be $\rho_0=1$.
Since the coefficients are bounded on $B_{2r}(x)$ and
independent of time, the constant below can be chosen independently of $t\geq2$.
Thus, whenever $0<h\leq\rho_0$ and $D_{3h}(y,s)\subset Q_2$, we have
\begin{equation}\label{eq:aronson-local-boundedness}
 |f^\varepsilon(y,s)|
 \leq C h^{-(d+2)/2}
 \left(\int_{D_{3h}(y,s)}|f^\varepsilon(z,\tau)|^2\,dz\,d\tau\right)^{1/2}.
\end{equation}

To apply this estimate uniformly on $Q_a$, let
$1\leq a<b\leq3/2$ and choose
\[
 h=c_r(b-a),\qquad c_r:=\min\{1,r/3,1/6\}.
\]
Since $b-a\leq1/2$ and $c_r\leq1/6$, we have
$0<h\leq1/12<\rho_0=1$; hence the scale bound is automatic.
Since $c_r\leq r/3$, $c_r\leq1/6$, and $0<b-a\leq1/2$,
\begin{align*}
 \frac{3h}{2}
 &=\frac32c_r(b-a)
 \leq\frac32\frac r3(b-a)
 =\frac{(b-a)r}{2}<(b-a)r,\\
 9h^2
 &=9c_r^2(b-a)^2
 \leq\frac{(b-a)^2}{4}
 \leq\frac{b-a}{8}<\frac{b-a}{2}.
\end{align*}
For $(y,s)\in Q_a$ and $(z,\tau)\in D_{3h}(y,s)$, it follows that
\[
 |z-x|\leq|z-y|+|y-x|
 <\frac{3h}{2}+ar<(b-a)r+ar=br,
\]
and
\[
 \tau>s-9h^2
 >t-\frac a2-\frac{b-a}{2}=t-\frac b2,
 \qquad
 \tau\leq s<t+\frac a2<t+\frac b2.
\]
Hence $D_{3h}(y,s)\subset Q_b$.
Taking the supremum over $(y,s)\in Q_a$ in
\eqref{eq:aronson-local-boundedness}, using $D_{3h}(y,s)\subset Q_b$
and $h=c_r(b-a)$, yields
\begin{equation}\label{eq:local-l2-boundedness}
 M(a)\leq C(b-a)^{-(d+2)/2}\|f^\varepsilon\|_{L^2(Q_b)}.
\end{equation}
Here and below, $C$ depends on $x,r,\varepsilon,d$ and the local
coefficient bounds, but not on $a,b,t$.

Note that
\[
 \|f^\varepsilon\|_{L^2(Q_b)}^2
 =\int_{Q_b}|f^\varepsilon(y,s)|^2\,dy\,ds
 \leq  M(b)J_t.
\]
Combining this with \eqref{eq:local-l2-boundedness} and
Young's inequality, for every $\theta>0$,
\begin{equation}\label{eq:nested-cylinder-interpolation}
 M(a)\leq C(b-a)^{-(d+2)/2}M(b)^{1/2}J_t^{1/2}
 \leq\theta M(b)+C_\theta(b-a)^{-(d+2)}J_t.
\end{equation}
The larger-cylinder supremum $M(b)$ still occurs on the right.  To
remove it, choose an increasing sequence of radii
\[
 a_j:=\frac32-2^{-j-1}\quad(j\geq0).
\]
Then $a_0=1$, $a_j\uparrow3/2$, and $a_{j+1}-a_j=2^{-j-2}$.
Applying \eqref{eq:nested-cylinder-interpolation} successively to
$Q_{a_j}\subset Q_{a_{j+1}}$ gives
\[
 M(1)\leq\theta^N M(a_N)
 +C_\theta J_t\sum_{j=0}^{N-1}\theta^j2^{(d+2)(j+2)}.
\]
Fix $\theta=2^{-(d+3)}$.  The sum is a convergent geometric series,
since $\theta2^{d+2}=1/2$, and
$\theta^N M(a_N)\leq\theta^N M(3/2)\to0$ as $N\to\infty$.
Consequently,
\[
 M(1)\leq C_{x,r,\varepsilon}J_t.
\]
This is the desired local $L^1$--$L^\infty$ estimate.  Notice that we only
used finiteness of $M(3/2)$ for each fixed $t$; no uniform-in-time
bound on this quantity was assumed.

\smallskip
\noindent\emph{Conclusion using Step~2.}
Since $(x,t)\in Q_1$, the preceding estimate implies
\begin{equation}\label{eq:pointwise-from-l1}
\begin{aligned}
 |f^\varepsilon(x,t)|
 &\leq C_{x,r,\varepsilon}
 \int_{t-1}^{t+1}\int_{B_{2r}(x)}|f^\varepsilon(y,s)|\,dy\,ds\\
 &\leq 2C_{x,r,\varepsilon}
 \sup_{s\geq t-1}\|f^\varepsilon(\cdot,s)\|_{L^1(\mathbb R^d)}
 \longrightarrow0,
\end{aligned}
\end{equation}
where the last limit follows from Step~2: for any $\delta>0$, the $L^1$ norm
is at most $\delta$ for every sufficiently large $s$, so its supremum
over $s\geq t-1$ also tends to zero.  Therefore
$\rho_\varepsilon(x,t)\to m_\varepsilon(x)$ for every $x\in\mathbb R^d$.

\end{proof}

The following lemma establishes the integrability of the initial profile $e^{-I/\varepsilon}$ for sufficiently small $\varepsilon>0$ and the logarithmic asymptotics of its normalizing constant.
\begin{lem}
\label{lem:wkb-normalization}
Assume (A1)--(A3), and let $I$ be defined in \eqref{eq:log-density-and-rate}. If $N$ is as in (A2), then
\[
 Z_\varepsilon:=\int_{\mathbb R^d}e^{-I(x)/\varepsilon}\,dx<\infty
 \qquad\text{for every }0<\varepsilon<\frac1{dN},
\]
and
\[
\lim_{\varepsilon\downarrow0} \varepsilon\log Z_\varepsilon = 0.
\]
\end{lem}

\begin{proof}
For $r>N$, set
\[
 q(r):=\min_{|x|=r}U(x).
\]
If $x_r$ is a minimizer of $U$ on $\partial B_r$, the Lagrange multiplier
rule gives
$\nabla U(x_r)=\alpha_r x_r$ for some constant $\alpha_r$.  Moreover, $\nabla U(x_r)\ne0$: otherwise
(A1) would give $\vb(x_r)=0$, contradicting (A2).  The orthogonality in
(A1) therefore implies
$\vl(x_r)\cdot x_r=0$.
Consequently, (A2) gives
\[
 \nabla U(x_r)\cdot x_r
 =-\vb(x_r)\cdot x_r>\frac1N.
\]
To account explicitly for the dependence of the constraint on $r$, write
\[
 q(r)=\min_{|y|=1}U(ry).
\]
Since
$U\in C^2$, the function $q$ is locally Lipschitz.  At almost every $r>N$
it is differentiable, and Danskin's theorem gives, for some minimizer $y_r$,
\[
 q'(r)=\partial_r U(ry_r)
      =\nabla U(ry_r)\cdot y_r
      =\nabla U(x_r)\cdot\frac{x_r}{r}
      >\frac1{Nr},
 \qquad x_r:=ry_r.
\]
Integrating this
inequality shows that, for any fixed $R>N$ and every $r\geq R$,
\begin{equation}
\label{eq:logarithmic-growth-U}
 q(r)\geq q(R)+\frac1N\log\frac rR.
\end{equation}

The coercivity proof in Proposition~\ref{prop:calibrated} gives a constant
$C_0$ such that $I(x)\geq U(x)-C_0$.  It follows from
\eqref{eq:logarithmic-growth-U}, upon fixing $R_0>N$, that there is a
constant $C_1$ such that
\begin{equation}
\label{eq:logarithmic-growth-I}
 I(x)\geq\frac1N\log|x|-C_1,
 \qquad |x|\geq R_0.
\end{equation}

For any $R\geq R_0$ and $0<\varepsilon<1/(dN)$,
\eqref{eq:logarithmic-growth-I} implies, for $|x|\geq R$,
\[
 e^{-I(x)/\varepsilon}
 \leq e^{C_1/\varepsilon}|x|^{-1/(N\varepsilon)}.
\]
Using polar coordinates, we therefore obtain
\begin{equation}
\label{eq:wkb-tail-bound}
\begin{aligned}
 \int_{|x|\geq R}e^{-I(x)/\varepsilon}\,dx
 \leq |\mathbb S^{d-1}|e^{C_1/\varepsilon}
       \int_R^\infty r^{d-1-1/(N\varepsilon)}\,dr =\frac{|\mathbb S^{d-1}|e^{C_1/\varepsilon}}
         {1/(N\varepsilon)-d}
   R^{\,d-1/(N\varepsilon)}<\infty.
\end{aligned}
\end{equation}
The integral on $B_R$ is finite by continuity of $I$, proving
$Z_\varepsilon<\infty$.

Since
$I\geq0$ and $I(x_*)=0$ for some $x_*\in\mathcal A$, for every $\delta>0$
there is $r_\delta>0$ such that $I\leq\delta$ on
$B(x_*,r_\delta)$.  Hence
\[
 Z_\varepsilon\geq |B(x_*,r_\delta)|e^{-\delta/\varepsilon},
\]
so $\liminf_{\varepsilon\downarrow0}\varepsilon\log Z_\varepsilon
\geq-\delta$.  Letting $\delta\downarrow0$ gives
$\liminf_{\varepsilon\downarrow0}\varepsilon\log Z_\varepsilon\geq0$.

For the upper bound, fix $R\geq R_0$.  Taking the logarithm of
\eqref{eq:wkb-tail-bound} and multiplying by $\varepsilon$ gives
\begin{align*}
 \varepsilon\log
 \int_{|x|\geq R}e^{-I(x)/\varepsilon}\,dx
 \leq C_1+\left(d\varepsilon-\frac1N\right)\log R +\varepsilon\log
 \left(\frac{|\mathbb S^{d-1}|}{1/(N\varepsilon)-d}\right).
\end{align*}
For fixed $R$, the last term tends to zero as $\varepsilon\downarrow0$;
hence
\[
 \limsup_{\varepsilon\downarrow0}
 \varepsilon\log
 \int_{|x|\geq R}e^{-I(x)/\varepsilon}\,dx
 \leq C_1-\frac1N\log R.
\]
On the other hand, $I\geq0$ gives
$\int_{B_R}e^{-I/\varepsilon}\,dx\leq |B_R|$, and therefore
\[
 \limsup_{\varepsilon\downarrow0}\varepsilon\log Z_\varepsilon
 \leq \max\left\{0,C_1-\frac1N\log R\right\}.
\]
Finally,
letting $R\to\infty$ gives
$\limsup_{\varepsilon\downarrow0}\varepsilon\log Z_\varepsilon\leq0$.
Together with the lower bound, this proves
$\varepsilon\log Z_\varepsilon\to0$.
\end{proof}

\begin{prop}[Finite-time zero-noise limit]
\label{prop:finite-time-wkb}
Assume (A1), (A2), and (A4).
Let $I$ be defined in \eqref{eq:log-density-and-rate}.
For $0<\varepsilon<1/(dN)$, let $\rho_\varepsilon$ be the probability solution of \eqref{FP} with initial density
\[
\rho_\varepsilon(x,0)
=Z_\varepsilon^{-1}e^{-I(x)/\varepsilon},
\qquad
Z_\varepsilon:=\int_{\mathbb R^d}e^{-I(x)/\varepsilon}\,dx.
\]
Then, for every fixed $t\geq0$,
\begin{equation}\label{eq:finite-time-wkb}
-\varepsilon\log\rho_\varepsilon(x,t)
\longrightarrow (S_tI)(x)=I(x)
\end{equation}
locally uniformly in $x\in\mathbb R^d$ as $\varepsilon\downarrow0$.
\end{prop}

\begin{proof}
For $t=0$, the identity
$-\varepsilon\log\rho_\varepsilon(x,0)=I(x)+\varepsilon\log Z_\varepsilon$
and Lemma~\ref{lem:wkb-normalization} give the assertion.

Fix $t>0$ and a compact set $K\subset\mathbb R^d$. We first establish a local smoothing estimate.  Since $t$ is fixed,
we may restrict to
$0<\varepsilon<\min\{1,1/(dN),t/4\}$.
For $f\geq0$ in $L^1(\mathbb R^d)$, let $r_\varepsilon(y,s)$ solve
\eqref{FP} with initial data $r_\varepsilon(y,0)=f(y)$.

Set
\[
 B_2:=B_2(0)=\{z\in\mathbb R^d:|z|<2\},
 \qquad K_2:=\{y:\operatorname{dist}(y,K)\leq2\}.
\]
Define the fixed space--time cylinder in the rescaled coordinates by
\[
 D_*:=B_2(0)\times(-4,0].
\]
For each fixed $x\in K$, define
\[
 v(z,s)=r_\varepsilon(x+\varepsilon z,t+\varepsilon s),
 \qquad (z,s)\in D_*.
\]
Under this change of variables, $D_*$ corresponds to
$B_{2\varepsilon}(x)\times(t-4\varepsilon,t]$ in the original
coordinates.  Since $\varepsilon<t/4$, this time interval lies after
the initial time, so $v$ is defined throughout $D_*$.

Define $a(y):=\nabla\cdot\vec b(y)$ for $y\in\mathbb R^d$.  Then
\[
 \partial_s v-\Delta_zv
 +\vec b(x+\varepsilon z)\cdot\nabla_zv
 +\varepsilon a(x+\varepsilon z)v=0.
\]
For $x\in K$, $z\in B_2$, and $\varepsilon<1$, we have
$x+\varepsilon z\in K_2$ and
$|\vec b(x+\varepsilon z)|+\varepsilon|a(x+\varepsilon z)|
\leq\sup_{K_2}(|\vec b|+|a|)<\infty$.
Thus the local parabolic estimates below have constants independent of
$x$, $t$, and $\varepsilon$.
For the cylinders $Q(\sigma):=B_{\sigma/2}(0)\times(-\sigma^2,0]$,
choose $\rho_0=1$ and $\rho=1/2$.  Since
$Q(3\rho)=B_{3/4}(0)\times(-9/4,0]\subset D_*$,
the local estimate in \cite[Theorem~B, p.~616]{Aronson1968Nonnegative}
yields
\[
 |v(0,0)|\leq C_K
 \left(\int_{-9/4}^{0}\int_{B_{3/4}(0)}v(z,s)^2\,dz\,ds\right)^{1/2}.
\]
Applying the interpolation and iteration in
\eqref{eq:local-l2-boundedness}--\eqref{eq:nested-cylinder-interpolation}
to the backward cylinders
$\widehat Q_\alpha:=B_\alpha(0)\times(-\alpha^2,0]$, $1\leq\alpha\leq3/2$,
gives $v(0,0)\leq C_K\int_{-4}^{0}\int_{B_2}v(z,s)\,dz\,ds$, since $v\geq0$.
Changing variables and using mass conservation, we obtain
\begin{equation}
\begin{aligned}
 r_\varepsilon(x,t)
 &\leq C_K\varepsilon^{-d-1}
 \int_{t-4\varepsilon}^{t}\int_{B_{2\varepsilon}(x)}
 r_\varepsilon(y,s)\,dy\,ds\\
 &\leq C_K\varepsilon^{-d}\|f\|_{L^1(\mathbb R^d)},
 \qquad x\in K.
\end{aligned}
\label{eq:wkb-local-density-bound}
\end{equation}

Now fix $\eta>0$ and put $M_K:=\sup_K I$.
By coercivity of $I$ and \eqref{eq:wkb-tail-bound}, together with
$\varepsilon\log Z_\varepsilon\to0$, we can choose $R$ so large that
$K\subset B_R$ and $J_R:=\inf_{\partial B_R}I>M_K+3\eta$, with
\begin{equation}
 \limsup_{\varepsilon\downarrow0}\varepsilon\log
 \int_{\mathbb R^d\setminus B_R}\rho_\varepsilon(y,0)\,dy
 <-M_K-3\eta.
 \label{eq:wkb-initial-tail-rate}
\end{equation}
Corollary~\ref{cor:density-convergence-strong-regularity} and
Lemma~\ref{lem:wkb-normalization} give
$\sup_{y\in\overline{B_R}}
| -\varepsilon\log m_\varepsilon(y)-I(y)|\leq\eta/2$ and
$|\varepsilon\log Z_\varepsilon|\leq\eta/2$
for sufficiently small $\varepsilon$.
Using $\rho_\varepsilon(y,0)=Z_\varepsilon^{-1}e^{-I(y)/\varepsilon}$,
we obtain, for $y\in\overline{B_R}$,
\[
 \left|\varepsilon\log\frac{\rho_\varepsilon(y,0)}{m_\varepsilon(y)}\right|
 =\left| -\varepsilon\log m_\varepsilon(y)-I(y)
 -\varepsilon\log Z_\varepsilon\right|\leq\eta.
\]
Dividing by $\varepsilon$ and exponentiating yields
\begin{equation}
 e^{-\eta/\varepsilon}m_\varepsilon(y)
 \leq\rho_\varepsilon(y,0)\leq e^{\eta/\varepsilon}m_\varepsilon(y),
 \qquad y\in\overline{B_R}.
 \label{eq:wkb-initial-density-comparison}
\end{equation}
Let $q_\varepsilon$ solve \eqref{FP} with initial data
$q_\varepsilon(y,0)=\rho_\varepsilon(y,0)
\mathbf1_{\mathbb R^d\setminus B_R}(y)$.
By linearity and stationarity, $e^{\eta/\varepsilon}m_\varepsilon
+q_\varepsilon$ solves the same equation as $\rho_\varepsilon$.
Its initial data dominate $\rho_\varepsilon(\cdot,0)$ on $B_R$ by
\eqref{eq:wkb-initial-density-comparison}, and outside $B_R$ by the
definition of $q_\varepsilon$. The nonnegative transition kernel
from the proof of Proposition~\ref{prop:fixed-epsilon-density-convergence}
preserves the nonnegativity of this initial difference and therefore gives
\[
 \rho_\varepsilon(\cdot,t)
 \leq e^{\eta/\varepsilon}m_\varepsilon+q_\varepsilon(\cdot,t).
\]
By \eqref{eq:wkb-local-density-bound} and the strict limsup bound
\eqref{eq:wkb-initial-tail-rate}, for sufficiently small $\varepsilon$,
\[
 \sup_K q_\varepsilon(\cdot,t)
 \leq C_K\varepsilon^{-d}
 \int_{\mathbb R^d\setminus B_R}\rho_\varepsilon(y,0)\,dy
 \leq C_K\varepsilon^{-d}e^{-(M_K+3\eta)/\varepsilon}.
\]
Since $\inf_K m_\varepsilon\geq e^{-(M_K+\eta)/\varepsilon}$ for
sufficiently small $\varepsilon$, we have
\[
 \sup_{x\in K}\frac{q_\varepsilon(x,t)}{m_\varepsilon(x)}
 \leq\frac{\sup_Kq_\varepsilon(\cdot,t)}{\inf_Km_\varepsilon}
 \leq C_K\varepsilon^{-d}e^{-2\eta/\varepsilon}
 \longrightarrow0.
\]
Thus the tail contribution is negligible relative to $m_\varepsilon$,
uniformly on $K$.  In particular, for sufficiently small $\varepsilon$,
$q_\varepsilon(x,t)\leq m_\varepsilon(x)$ for every $x\in K$.
The additive error in the preceding upper bound can therefore be
absorbed into the multiplicative estimate
\[
 \rho_\varepsilon(x,t)
 \leq\bigl(e^{\eta/\varepsilon}+1\bigr)m_\varepsilon(x)
 \leq2e^{\eta/\varepsilon}m_\varepsilon(x),\qquad x\in K,
\]
which permits a logarithmic comparison with $m_\varepsilon$.

For the lower bound, let $A_\varepsilon:=\sup_{\partial B_R}m_\varepsilon$
and $C_R:=\|a\|_{L^\infty(B_R)}$.
The function $h_\varepsilon(y,s):=e^{-\eta/\varepsilon}
\bigl(m_\varepsilon(y)-A_\varepsilon e^{C_Rs}\bigr)$
is a subsolution of \eqref{FP} in
$B_R\times(0,t]$ since
\[
 \partial_s h_\varepsilon-\varepsilon\Delta h_\varepsilon
 +\vec b\cdot\nabla h_\varepsilon+a(y)h_\varepsilon
 =-e^{-\eta/\varepsilon}A_\varepsilon e^{C_Rs}(C_R+a(y))\leq0.
\]
On $\partial B_R\times(0,t]$,
\[
 h_\varepsilon(y,s)
 \leq e^{-\eta/\varepsilon}A_\varepsilon(1-e^{C_Rs})
 \leq0\leq\rho_\varepsilon(y,s).
\]
Moreover,
\eqref{eq:wkb-initial-density-comparison} gives
$h_\varepsilon(\cdot,0)\leq\rho_\varepsilon(\cdot,0)$ in $B_R$.
Parabolic comparison therefore yields
\[
 \rho_\varepsilon(y,t)\geq e^{-\eta/\varepsilon}
 \bigl(m_\varepsilon(y)-A_\varepsilon e^{C_Rt}\bigr),
 \qquad y\in B_R.
\]
Since $u_\varepsilon=-\varepsilon\log m_\varepsilon$ converges uniformly
to $I$ on $\overline{B_R}$, we have
$I(y)-\eta\leq u_\varepsilon(y)\leq I(y)+\eta$ there for sufficiently
small $\varepsilon$.
Using $J_R=\inf_{\partial B_R}I$ and $M_K=\sup_K I$, we obtain
\[
 A_\varepsilon=\sup_{\partial B_R}m_\varepsilon
 \leq e^{-(J_R-\eta)/\varepsilon},
 \qquad
 \inf_Km_\varepsilon
 \geq e^{-(M_K+\eta)/\varepsilon}.
\]
Dividing these bounds gives
\[
 \frac{A_\varepsilon e^{C_Rt}}{\inf_K m_\varepsilon}
 \leq e^{C_Rt}
 \exp\left(-\frac{J_R-M_K-2\eta}{\varepsilon}\right)
 \longrightarrow0.
\]
The two errors $\eta$ account for the term $2\eta$ in the exponent.
The limit is zero because $J_R>M_K+3\eta$ and $R$ and $t$ are fixed.
In particular, for sufficiently small $\varepsilon$,
$A_\varepsilon e^{C_Rt}\leq\tfrac12\inf_Km_\varepsilon
\leq\tfrac12m_\varepsilon(x)$ for every $x\in K$.
Combining this with the preceding lower and upper bounds yields
\[
 \frac12 e^{-\eta/\varepsilon}m_\varepsilon(x)
 \leq\rho_\varepsilon(x,t)
 \leq2e^{\eta/\varepsilon}m_\varepsilon(x),
 \qquad x\in K.
\]
Taking logarithms gives
\[
 \sup_{x\in K}
 \left| -\varepsilon\log\rho_\varepsilon(x,t)-u_\varepsilon(x)\right|
 \leq\eta+\varepsilon\log2.
\]
Letting $\varepsilon\downarrow0$ and then $\eta\downarrow0$, and using
$u_\varepsilon\to I$ locally uniformly, proves the assertion at the fixed
time $t$.  Finally, Corollary~\ref{cor:laxoleinik} gives $S_tI=I$.
\end{proof}

\begin{thm}[Exchange of the long-time and zero-noise limits]
\label{thm:exchange}
Assume (A1)--(A3). Then the function $I$ defined in \eqref{eq:log-density-and-rate}
is the unique coercive weak KAM solution of \eqref{eq:HJE0} with the
LDP-selected boundary values
\begin{equation}\label{eq:selected-aubry-boundary-values}
I(x_i)=W_0(x_i)-\min_{x_j\in\mathcal A}W_0(x_j),
\qquad x_i\in\mathcal A.
\end{equation}

Suppose in addition that (A4) holds, and let $\rho_\varepsilon$ be the probability solution of \eqref{FP} with the initial density specified in
Proposition~\ref{prop:finite-time-wkb}.  Then, for every
$x\in\mathbb R^d$, both iterated limits exist and
\begin{equation}
\begin{aligned}
    \lim_{\varepsilon\downarrow0}
    \left(-\varepsilon\log
    \lim_{t\to+\infty}\rho_\varepsilon(x,t)\right)
    =I(x)=\lim_{t\to+\infty}\lim_{\varepsilon\downarrow0}
    \bigl(-\varepsilon\log\rho_\varepsilon(x,t)\bigr).
\end{aligned}
\end{equation}
\end{thm}

\begin{proof}
The graph formula \eqref{eqn:disbdy} determines the LDP-selected values of
$W_0$ on $\mathcal A$, and subtracting their minimum only fixes the additive
constant.  The Peierls-barrier extension in Proposition~\ref{prop:ldp}
therefore gives the weak KAM solution $I$ with the boundary values
\eqref{eq:selected-aubry-boundary-values}.  Its uniqueness within the class
of coercive weak KAM solutions follows from
Corollary~\ref{cor:unique-coercive-weak-kam}.

For each fixed $0<\varepsilon<1/(dN)$,
Proposition~\ref{prop:fixed-epsilon-density-convergence} gives
\[
\lim_{t\to+\infty}\rho_\varepsilon(x,t)=m_\varepsilon(x).
\]
Corollary~\ref{cor:density-convergence-strong-regularity} therefore yields
\[
\begin{aligned}
\lim_{\varepsilon\downarrow0}
\left(-\varepsilon\log
\lim_{t\to+\infty}\rho_\varepsilon(x,t)\right)=\lim_{\varepsilon\downarrow0}
\bigl(-\varepsilon\log m_\varepsilon(x)\bigr)=I(x).
\end{aligned}
\]

Proposition~\ref{prop:finite-time-wkb} and
Corollary~\ref{cor:laxoleinik} imply, for every fixed $t>0$,
\[
\lim_{\varepsilon\downarrow0}
\bigl(-\varepsilon\log\rho_\varepsilon(x,t)\bigr)
=(S_tI)(x)=I(x).
\]
The right-hand side is independent of $t$, and consequently its limit as
$t\to+\infty$ is again $I(x)$.
\end{proof}

\begin{rem}
\label{rem:initial-profile-commutation}
The two iterated limits agree because the initial logarithmic profile is the
selected weak KAM solution $I$, for which $S_tI=I$.
For general initial data, they need not commute. For example, take
$\vb(x)=-x$ and $\rho_\varepsilon(x,0)=(2\pi)^{-d/2}e^{-|x|^2/2}$.
Then $I(x)=|x|^2/2$, and $\rho_\varepsilon(\cdot,t)$ is a centered Gaussian
with covariance $[e^{-2t}+\varepsilon(1-e^{-2t})]I_d$. Hence
\[
 \lim_{\varepsilon\downarrow0}\lim_{t\to\infty}
 [-\varepsilon\log\rho_\varepsilon(x,t)]=\frac{|x|^2}{2},
 \qquad
 \lim_{t\to\infty}\lim_{\varepsilon\downarrow0}
 [-\varepsilon\log\rho_\varepsilon(x,t)]=0.
\]
\end{rem}

\begin{rem}[Comparison with the vanishing discount selection]\label{rem}
On a compact manifold, the vanishing discount limit selects a viscosity solution; see \cite[Theorem~1.1]{DaviniFathiIturriagaZavidovique2016}. The corresponding result in $\mathbb R^d$ holds under additional hypotheses
at infinity; see \cite{ishii2020vanishing}. In general, the solution selected
by vanishing discount can differ from the one selected by vanishing viscosity.
General selection principles for viscosity solutions and Mather measures,
including vanishing-discount and vanishing-viscosity problems, are discussed
in \cite{Gomes2008}.
Moreover, in some settings the vanishing-viscosity family need not converge:
for a mechanical Hamiltonian on the torus, such a failure is exhibited in
\cite{LiuTranYu2026}. For simultaneous vanishing discount and viscosity on
the torus, the selected limit can also depend on the relative scaling of the
two terms; in the critical scaling regime, an explicit limit is obtained in
\cite{MitakeNiTran2026}. The following example illustrates the difference
between the two selections when both limits exist.

Consider the one-dimensional gradient drift $b=-U'$, with
\[
 U(x)=\begin{cases}
 \dfrac{(1-x^2)^2(3-x^2)}8,& |x|\leq1,\\[4pt]
 (|x|-1)^2,& |x|\geq1.
 \end{cases}
\]
Note that $U\in C^3(\mathbb R)$, and (A1)--(A4) hold with
$\ell=0$. The only equilibria are $-1,0,1$, with
$U(\pm1)=0$ and $U(0)=3/8$. The invariant density is
$m_\varepsilon=Z_\varepsilon^{-1}e^{-U/\varepsilon}$, where
$Z_\varepsilon=\int_{\mathbb R}e^{-U/\varepsilon}\,dx$. Hence
\[
 u_\varepsilon=U+\varepsilon\log Z_\varepsilon
 \longrightarrow I=U
 \qquad\text{locally uniformly as }\varepsilon\downarrow0.
\]
On the other hand, for $0<\lambda<4$, let
\[
 \psi_\lambda(x)=
 \begin{cases}
  0,& |x|\leq1,\\[2pt]
  \left(1-\dfrac\lambda4\right)(|x|-1)^2,& |x|\geq1.
 \end{cases}
\]
Note that $\psi_\lambda\in C^1(\mathbb R)$ is coercive and solves
\[
 \lambda\psi_\lambda+(\psi_\lambda')^2+b\psi_\lambda'=0
 \qquad\text{in }\mathbb R.
\]
Consequently,
\[
 \psi_\lambda\longrightarrow\psi_0,\qquad
 \psi_0(x)=
 \begin{cases}
  0,& |x|\leq1,\\[2pt]
  (|x|-1)^2,& |x|>1,
 \end{cases}
 \qquad \psi_0(0)=0\ne\frac38=I(0).
\]
Both limits are coercive and have minimum zero, but they differ on
$(-1,1)$. Thus the vanishing discount limit need not agree with the
vanishing-viscosity limit, even when the discounted solutions are coercive. Moreover, the discounted problem for our Hamiltonian in general does not have the uniqueness even within the coercive class.
\end{rem}

\section*{Acknowledgment}
Yuan Gao was partially supported by NSF under  CAREER award DMS-2440651. Yuxi Han was partially supported by NSF under award DMS-2440651.

\section*{Conflict Of Interest Statement}
The authors have no competing interests to declare that are relevant to the content of this article.

\section*{Data Availability Statement}
No data were used for the research described in this article.

\bibliographystyle{alpha}
\bibliography{kam_path.bib}

\appendix

\bigskip

\section{Some comparison results}

We relate the stationary solution to the time-dependent equation \eqref{eq:limit-HJ}. The following action bound for viscosity subsolutions will provide the upper comparison with $W$.

\begin{lem}\label{lem:lax-oleinik-subsolution}
Assume (A1) and (A3). Let
$u\in\mathrm{USC}(\mathbb R^d\times[0,\infty))$ be a locally bounded
viscosity subsolution of \eqref{eq:limit-HJ} on
$\mathbb R^d\times(0,\infty)$. Then, for every $x,y\in\mathbb R^d$ and $t>0$,
\begin{equation}\label{eq:subsolution-domination}
 u(y,t)-u(x,0)\leq h_t(x,y).
\end{equation}
Consequently,
\begin{equation}\label{eq:lax-oleinik-subsolution}
 u(y,t)\leq\inf_{x\in\mathbb R^d}\{u(x,0)+h_t(x,y)\}.
\end{equation}
\end{lem}

\begin{proof}
Fix $x,y\in\mathbb R^d$ and $t>0$. By
Proposition~\ref{prop:minimizer-regularity}, there exists
$\gamma\in C^2([0,t];\mathbb R^d)$ such that
\[
 \gamma(0)=x,\qquad \gamma(t)=y,\qquad
 h_t(x,y)=\int_0^t L(\dot\gamma(s),\gamma(s))\,ds.
\]
Define
\[
F:\mathbb{R}^d\times[0,t]\to\mathbb{R},
\qquad
F(z,s):=u(z+\gamma(s),s).
\]
Since \(u\) is upper semicontinuous and \(\gamma\) is continuous, $F\in \textrm{USC} \bigl(\mathbb{R}^d\times[0,t]\bigr)$.
Let $(z_0,s_0)\in\mathbb{R}^d\times(0,t)$ and suppose that $\Phi\in C^1\bigl(\mathbb{R}^d\times(0,t)\bigr)$ touches \(F\) from above at \((z_0,s_0)\). Define $\Psi:\mathbb{R}^d\times(0,t)\to\mathbb{R}$ by $\Psi(y,s):=\Phi(y-\gamma(s),s)$.
Then \(\Psi\) touches \(u\) from above at
\[
(y_0,s_0)
:=
(z_0+\gamma(s_0),s_0).
\]
Since \(u\) is a viscosity subsolution of
\begin{equation}
u_t(y,s)+H(Du(y,s),y)=0,
\qquad
(y,s)\in\mathbb{R}^d\times(0,\infty),
\label{eq:HJ-u}
\end{equation}
we have
\[
\Psi_s(y_0,s_0)
+
H(D_y\Psi(y_0,s_0),y_0)
\leq0.
\]
Using
\[
D_y\Psi(y,s)
=
D_z\Phi(y-\gamma(s),s), \qquad
\Psi_s(y,s)
=
\Phi_s(y-\gamma(s),s)
-
\dot\gamma(s)\cdot
D_z\Phi(y-\gamma(s),s),
\]
we obtain
\[
\Phi_s(z_0,s_0)
-
\dot\gamma(s_0)\cdot D_z\Phi(z_0,s_0)
+
H\bigl(D_z\Phi(z_0,s_0),z_0+\gamma(s_0)\bigr)
\leq0.
\]
By the definition of the Legendre transform,
\[
\Phi_s(z_0,s_0)
\leq
L\bigl(\dot\gamma(s_0),z_0+\gamma(s_0)\bigr).
\]
Thus \(F\) is a viscosity subsolution of
\begin{equation}
F_s(z,s)
\leq
L\bigl(\dot\gamma(s),z+\gamma(s)\bigr), \qquad (z, s) \in \mathbb{R}^d\times(0,t).
\label{eq:F-L}
\end{equation}

Now define $w:[0,t]\to\mathbb{R}$ by $w(s):=F(0,s)=u(\gamma(s),s)$. Then $w\in \mathrm{USC}([0,t])$.
We claim that \(w\) is a viscosity subsolution of
\begin{equation}
w'(s)
\leq
L(\dot\gamma(s),\gamma(s)), \qquad s \in (0,t).
\label{eq:w-ode}
\end{equation}
Fix $s_0\in(0,t)$, and let $\phi\in C^1((0,t))$ and suppose \(w-\phi\) has a strict maximum at \(s_0 \in (0,t)\).  Choose $\delta>0$ so that $[s_0-\delta,s_0+\delta]\subset(0,t)$, and consider the compact set
\[
K_\delta
:=
\overline{B_\delta(0)}
\times
[s_0-\delta,s_0+\delta]
\subset
\mathbb{R}^d\times(0,t).
\]
For \(\varepsilon>0\), let $(z_\varepsilon,s_\varepsilon)\in K_\delta$ be a maximum point of
\begin{equation}
(z,s)\mapsto
F(z,s)-\phi(s)-\frac{|z|^2}{\varepsilon}.
\label{eq:penalized-function}
\end{equation}
Since \(F\) is upper semicontinuous, such a maximum exists. Comparing the maximum value with the value at \((0,s_0)\), we obtain
\begin{equation} \label{eq:afcompare}
F(z_\varepsilon,s_\varepsilon)
-
\phi(s_\varepsilon)
-
\frac{|z_\varepsilon|^2}{\varepsilon}
\geq
F(0,s_0)-\phi(s_0).
\end{equation}
Since \(F\) is bounded above on \(K_\delta\), there exists \(C_\delta>0\)
such that
\[
\frac{|z_\varepsilon|^2}{\varepsilon}
\leq C_\delta.
\]
Hence
\begin{equation}
z_\varepsilon\to0
\qquad\text{as }\varepsilon\to0.
\label{eq:z-epsilon}
\end{equation}
By \eqref{eq:afcompare}
\begin{equation}
F(z_\varepsilon,s_\varepsilon)-\phi(s_\varepsilon)
\ge
F(0,s_0)-\phi(s_0).
\label{eq:seps-maximality}
\end{equation}
Let
\(\varepsilon_j\to0\) be any subsequence such that
\[
s_{\varepsilon_j}\to \bar s
\in[s_0-\delta,s_0+\delta].
\]
By upper semicontinuity of \(F\), $\limsup_{j\to\infty}
F(z_{\varepsilon_j},s_{\varepsilon_j})
\le
F(0,\bar s)$.
Taking the \(\limsup\) in \eqref{eq:seps-maximality} and using the
continuity of \(\phi\), we obtain
\[
F(0,\bar s)-\phi(\bar s)
\ge
F(0,s_0)-\phi(s_0).
\]
Since \(F(0,s)=w(s)\), this becomes
\[
w(\bar s)-\phi(\bar s)
\ge
w(s_0)-\phi(s_0).
\]
However, \(w-\phi\) has a strict maximum at \(s_0\). Therefore, $\bar s=s_0$.
Since every convergent subsequence of \(s_\varepsilon\) has the same
limit \(s_0\), we conclude that
\begin{equation}\label{eq:s-epsilon}
 s_\varepsilon\to s_0
\qquad\text{as }\varepsilon\to0.
\end{equation}
For sufficiently small \(\varepsilon\), $(z_\varepsilon,s_\varepsilon)
\in
\mathbb{R}^d\times(0,t)$, so we may apply \eqref{eq:F-L} using the test function
\[
\Phi_\varepsilon(z,s)
:=
\phi(s)+\frac{|z|^2}{\varepsilon},
\qquad
(z,s)\in\mathbb{R}^d\times(0,t),
\]
which gives
\[
\phi'(s_\varepsilon)
\leq
L\bigl(
\dot\gamma(s_\varepsilon),
z_\varepsilon+\gamma(s_\varepsilon)
\bigr).
\]
Letting \(\varepsilon\to0\), and using
\eqref{eq:z-epsilon}, \eqref{eq:s-epsilon}, and the continuity of \(L\),
we obtain
\[
\phi'(s_0)
\leq
L(\dot\gamma(s_0),\gamma(s_0)).
\]
Therefore \eqref{eq:w-ode} holds in the viscosity sense on \((0,t)\).

Define $G:[0,t]\to\mathbb{R}$ by
\begin{equation}
G(s)
:=
w(s)
-
\int_0^s
L(\dot\gamma(r),\gamma(r))\,dr.
\label{eq:def-G}
\end{equation}
Then \(G\) is a viscosity subsolution  of
\[
G'(s)\leq0, \qquad s \in (0,t).
\]
We first show that $G$ is nonincreasing on $(0,t)$.
Suppose, to the contrary, that $G(b)>G(a)$ for some $0<a<b<t$.
Choose $c\in(b,t)$ and $\eta,\varepsilon>0$ sufficiently small that
\[
 G(b)-\eta(b-a)-\frac{\varepsilon}{c-b}
 >G(a)-\frac{\varepsilon}{c-a}.
\]
Since $G$ is upper semicontinuous and bounded above on $[a,c]$,
the function
\[
 G(s)-\eta(s-a)-\frac{\varepsilon}{c-s},\qquad a\leq s<c,
\]
tends to $-\infty$ as $s\uparrow c$ and attains a maximum.
The preceding strict inequality excludes $a$, so a maximum point
$s_*\in(a,c)$ exists. Thus, for an appropriate constant $C$,
\[
 \phi(s)=\eta(s-a)+\frac{\varepsilon}{c-s}+C
\]
touches $G$ from above at $s_*$, whereas
\[
 \phi'(s_*)=\eta+\frac{\varepsilon}{(c-s_*)^2}>0.
\]
This contradicts the viscosity subsolution property of $G$.
Hence $G(b)\leq G(a)$ for every $0<a<b<t$.

To include the terminal time $t$, fix $\rho>0$ and define
\[
 \widetilde\gamma(s)=
 \begin{cases}
 \gamma(s),&0\leq s\leq t,\\
 \gamma(t)+(s-t)\dot\gamma(t),&t<s\leq t+\rho.
 \end{cases}
\]
The values and first derivatives match at $t$, so
$\widetilde\gamma\in C^1([0,t+\rho];\mathbb R^d)$.
Define $\widetilde w$ and $\widetilde G$ by the same formulas as
$w$ and $G$, with $\gamma$ replaced by $\widetilde\gamma$.
The preceding test-function argument applies on $(0,t+\rho)$,
since $u$ is a subsolution for all positive times. The monotonicity
argument just proved then gives
$\widetilde G(t)\leq\widetilde G(a)$ for every $0<a<t$.
Since $\widetilde G=G$ on $[0,t]$, we have $G(t)\leq G(a)$.
The auxiliary extension need not minimize the action beyond $t$.
At the initial time, upper semicontinuity gives
\[
 \limsup_{a\downarrow0}G(a)
 =\limsup_{a\downarrow0}u(\gamma(a),a)
 \leq u(x,0)=G(0).
\]
Therefore
\[
G(t)\leq G(0).
\]
Equivalently,
\[
w(t)-w(0)
\leq
\int_0^t
L(\dot\gamma(s),\gamma(s))\,ds.
\]
Since $w(t)=u(y,t)$, $w(0)=u(x,0)$, and $\gamma$ minimizes the
action on $[0,t]$, this proves \eqref{eq:subsolution-domination}.
Taking the infimum over $x\in\mathbb R^d$ gives
\eqref{eq:lax-oleinik-subsolution}.
\end{proof}

Combining this action bound with $S_tW=W$ controls subsolutions from above. Uniform coercivity of supersolutions supplies the boundary condition at infinity for the opposite inequality.

\begin{prop}
\label{prop:coercive-comparison}
Under the assumptions of Proposition~\ref{prop:calibrated}, let $W$ be the function constructed there.
Let $u\in\mathrm{USC}(\mathbb R^d\times[0,\infty))$ and
$v\in\mathrm{LSC}(\mathbb R^d\times[0,\infty))$ be locally bounded
viscosity sub- and supersolutions, respectively, of \eqref{eq:limit-HJ}.
Suppose that
\[
 u(x,0)\leq W(x)\leq v(x,0),\qquad x\in\mathbb R^d,
\]
and
\begin{equation}
 \lim_{R\to\infty}
 \inf_{\substack{|y|\geq R\\0\leq s\leq T}}v(y,s)=+\infty
 \qquad\text{for every }T>0.
 \label{eq:comparison-uniform-coercivity}
\end{equation}
Then
\[
 u(x,t)\leq W(x)\leq v(x,t)
 \qquad\text{on }\mathbb R^d\times[0,\infty).
\]
In particular, $w(x,t)=W(x)$ is the unique continuous viscosity
solution with initial data $W$ in the class of viscosity solutions coercive
uniformly on finite time intervals.
\end{prop}

\begin{proof}
We first prove $u(x,t)\leq W(x)$ for every $(x,t)\in \mathbb{R}^d\times[0,\infty)$.

Fix
$(x,t)\in\mathbb R^d\times(0,\infty)$.
By Lemma~\ref{lem:lax-oleinik-subsolution},
\[
u(x,t)
\leq
\inf_{y\in\mathbb{R}^d}
\{u(y,0)+h_t(y,x)\}.
\]
Since \(u(y,0)\leq W(y)\) for every $y\in\mathbb R^d$,
\[
u(x,t)
\leq
\inf_{y\in\mathbb{R}^d}
\{W(y)+h_t(y,x)\}.
\]
By Corollary~\ref{cor:laxoleinik},
\[
W(x)
=
\inf_{y\in\mathbb{R}^d}
\{W(y)+h_t(y,x)\}.
\]
Therefore, $u(x,t)\leq W(x)$.

Next, we prove that \(W(x)\leq v(x,t)\).

Fix $(x,t)\in\mathbb R^d\times(0,\infty)$ and choose $T>t$. Choose $M>W(x)$ and define
\[
 W_M:\mathbb{R}^d\to\mathbb{R}, \qquad W_M(y):=\min\{W(y),M\}.
\]

We first prove that \(W_M\) is dominated by \(L\). Let $\tau>0$ and let $\eta:[0,\tau]\to\mathbb{R}^d$ be absolutely continuous. Since \(W\prec L\),
\begin{equation}
 W(\eta(\tau))-W(\eta(0)) \leq \int_0^\tau L(\dot\eta(s),\eta(s))\,ds. \label{eq:W-dominated}
\end{equation}
Moreover,
\begin{equation}
 W_M(\eta(\tau))-W_M(\eta(0)) \leq \max\left\{ W(\eta(\tau))-W(\eta(0)),0 \right\}. \label{eq:min-inequality}
\end{equation}
Since $L(q,y)\geq0$ for all $(q,y)\in\mathbb R^d\times\mathbb R^d$,
it follows from \eqref{eq:W-dominated},\eqref{eq:min-inequality} that
\begin{equation}
 W_M(\eta(\tau))-W_M(\eta(0)) \leq \int_0^\tau L(\dot\eta(s),\eta(s))\,ds. \label{eq:WM-dominated}
\end{equation}

We next show that $W_M\in \mathrm{USC}(\mathbb{R}^d)$ is a viscosity subsolution on \(\mathbb{R}^d\) of
\begin{equation}
 H(DW_M(y),y)\leq0, \qquad y\in\mathbb{R}^d. \label{eq:WM-stationary}
\end{equation}
Let $y_0\in\mathbb{R}^d$ and suppose that $\varphi\in C^1(\mathbb{R}^d)$ touches \(W_M\) from above at \(y_0\). Fix $q\in\mathbb{R}^d$ and \(h>0\), and define $\eta_h:[0,h]\to\mathbb{R}^d$ by $\eta_h(s) = y_0-hq+sq$. Then
\[
 \eta_h(0)=y_0-hq, \qquad \eta_h(h)=y_0, \qquad \dot\eta_h(s)=q.
\]
Applying \eqref{eq:WM-dominated} gives
\[
 W_M(y_0)-W_M(y_0-hq) \leq \int_0^h L(q,y_0-hq+sq)\,ds.
\]
Since \(\varphi\) touches \(W_M\) from above at \(y_0\),
\[
 W_M(y)-\varphi(y) \leq W_M(y_0)-\varphi(y_0)
\]
for \(y\) sufficiently close to \(y_0\). In particular, for sufficiently
small \(h>0\),
\[
 \varphi(y_0)-\varphi(y_0-hq) \leq W_M(y_0)-W_M(y_0-hq).
\]
Hence
\[
 \frac{ \varphi(y_0)-\varphi(y_0-hq) }{h} \leq \frac1h \int_0^h L(q,y_0-hq+sq)\,ds.
\]
Letting \(h\downarrow0\) yields $D\varphi(y_0)\cdot q\leq L(q,y_0)$.
Since \(q\in\mathbb{R}^d\) is arbitrary,
\[
 H(D\varphi(y_0),y_0)
 =\sup_{q\in\mathbb R^d}\{D\varphi(y_0)\cdot q-L(q,y_0)\}\leq0.
\]
Thus \eqref{eq:WM-stationary} holds in the viscosity sense.

Equivalently, defining $\widetilde W_M: \mathbb{R}^d\times[0,T]\to\mathbb{R}$ by $\widetilde W_M(y,s):=W_M(y)$, we have that \(\widetilde W_M\) is a viscosity subsolution on $\mathbb{R}^d\times(0,T)$ of
\begin{equation}
 (\widetilde W_M)_s + H(D_y\widetilde W_M,y) \leq0. \label{eq:WM-evolution}
\end{equation}

Since
\begin{equation}
 \lim_{R\to\infty} \inf_{\substack{|y|\geq R\\ s\in[0,T]}} v(y,s) = +\infty, \label{eq:v-coercivity}
\end{equation}
there exists $R>|x|$ such that
\begin{equation}
 v(y,s)\geq M, \qquad (y,s)\in\partial B_R\times[0,T]. \label{eq:lateral-v}
\end{equation}
Since $W_M(y)\leq M$ for every $y\in\mathbb R^d$, we obtain
\begin{equation}
 W_M(y)\leq v(y,s), \qquad (y,s)\in\partial B_R\times[0,T]. \label{eq:lateral-order}
\end{equation}

At the initial time,
\[
 W_M(y)\leq W(y)\leq v(y,0), \qquad y\in\overline{B_R},
\]
so
\begin{equation}
 W_M(y)\leq v(y,0), \qquad y\in\overline{B_R}. \label{eq:initial-order}
\end{equation}

We now apply the viscosity comparison principle on the bounded
space-time cylinder $Q_{R,T} := B_R\times(0,T)$. Its parabolic boundary is
\[
 \partial_p Q_{R,T} = \bigl(\overline{B_R}\times\{0\}\bigr) \cup \bigl(\partial B_R\times [0,T) \bigr).
\]
By \eqref{eq:lateral-order} and \eqref{eq:initial-order},
\[
 W_M(y)\leq v(y,s) \qquad \text{on }\partial_p Q_{R,T}.
\]

On the bounded spatial domain $B_R$, the Hamiltonian
\[
 H:\mathbb{R}^d\times\overline{B_R}\to\mathbb{R}, \qquad H(p,y)=|p|^2+\vb(y)\cdot p,
\]
satisfies the spatial continuity condition for viscosity comparison. Indeed, by (A3), $\vb$ is Lipschitz on $\overline{B_R}$ with some constant $K_R<\infty$, so
\[
 |\vb(y)-\vb(z)| \leq K_R|y-z|, \qquad y,z\in\overline{B_R}.
\]
Consequently, for $y,z\in\overline{B_R}$, $\alpha>0$, and
$p=\alpha(y-z)$, we have
\begin{align*}
|H(p,y)-H(p,z)|
&=
|(\vb(y)-\vb(z))\cdot p|
\leq K_R|y-z|\,|p|
=K_R\alpha|y-z|^2.
\end{align*}

Hence the comparison principle on \(Q_{R,T}\) yields (\cite[Theorem~8.2]{Crandall_Ishii_Lions_1992})
\[
 W_M(y)\leq v(y,s), \qquad (y,s)\in\overline{B_R}\times[0,T).
\]
In particular, since $|x|<R$ and $t<T$, we obtain $W_M(x)\leq v(x,t)$. Because \(M>W(x)\), $W_M(x)=W(x)$. Together with the initial ordering, this proves
\begin{equation}
 W(x)\leq v(x,t), \qquad (x,t)\in\mathbb{R}^d\times[0,\infty). \label{eq:W-le-v}
\end{equation}

Combining the upper bound proved above with \eqref{eq:W-le-v}, we conclude that
\[
 u(x,t)\leq W(x)\leq v(x,t) \qquad \text{for every } (x,t)\in\mathbb{R}^d\times[0,\infty).
\]
\end{proof}

\end{document}